\def\psd{\mathbb{S}_{+}}
\newcommand{\off}[1]{\text{off}{#1}}

\documentclass[preprint,12pt]{elsarticle}
\usepackage{algorithm}
\usepackage{algorithmic}

\journal{}

\usepackage{amsfonts}
\usepackage{amssymb}
\usepackage{amsmath}
\usepackage{amsthm}
\usepackage{latexsym}
\usepackage{verbatim}
\usepackage{epsfig}
\usepackage{color}
										
\usepackage{graphicx}
\usepackage{amsbsy}
\usepackage{amscd}
\usepackage{bm}

\newcommand{\R}{{\mathbb R}}

\DeclareMathOperator*{\cov}{Cov}

\DeclareMathOperator*{\tr}{tr}
\DeclareMathOperator*{\diag}{diag}
\DeclareMathOperator*{\E}{\mathbf{E}}
\DeclareMathOperator*{\argmin}{argmin}
\DeclareMathOperator*{\argmax}{argmax}

\DeclareMathOperator*{\p}{\mathbf{P}}
\providecommand{\wt}[1]{\widetilde{#1}}
\providecommand{\wh}[1]{\widehat{#1}}

\providecommand{\norm}[1]{\left \lVert#1 \right \rVert}
\providecommand{\nnorm}[1]{ \lVert#1 \rVert}

\providecommand{\mc}[1]{\mathcal#1}
\providecommand{\T}{\top}

\newcommand{\indep}{\rotatebox[origin=c]{90}{$\models$}}

\newcommand{\blanco}[1]{  }

\newcommand{\deriv}[3]{%
\ifthenelse{#1 = 1}{\frac{d\,#2}{d\,#3}}{\frac{d^{{#1}} #2}{d{#3}^{{#1}}}}
}

\newcommand{\partials}[3]{%
\ifthenelse{#1 = 1}{\frac{\partial\,#2}{\partial\,#3}}{\frac{\partial^{#1}
    #2}{\partial#3^{#1}}}
} 

\def\su{\sum_{i=1}^n}

\def \coloneq{\mathrel{\mathop:}=}
\def \invcoloneq{=\mathrel{\mathop:}}
\def \eps{\varepsilon} 
\def \gec{\succeq}

\newtheorem{theo}{Theorem}
\newtheorem{propo}{Theorem}

\newtheorem{definitioA}{Theorem}[section]
\newtheorem{propoA}{Theorem}[section]
\newtheorem{lemmaA}{Theorem}[section]
\newtheorem{theoA}{Theorem}[section]
\newtheorem{remA}{Theorem}[section]

\newtheorem{theorem}[theo]{Theorem}
  
  \newtheorem{defnApp}[definitioA]{Definition}
  \newtheorem{prop}[propo]{Proposition}
\newtheorem{propApp}[propoA]{Proposition}

\newtheorem{lemmaApp}[lemmaA]{Lemma}
\newtheorem{theoremApp}[theoA]{Theorem}

\newtheorem{remApp}[remA]{Remark}

\newenvironment{bew}{\begin{proof}[Proof]}{\end{proof}}

\begin{document}

\begin{frontmatter}

  

\title{Estimation of positive definite $M$-matrices and structure learning for
attractive Gaussian Markov Random fields}


\author{Martin Slawski and Matthias Hein\\ {\footnotesize \textsf{\{ms,hein\}@cs.uni-saarland.de}}}

\address{\begin{tabular}{c} Saarland University \\ {\scriptsize
        Department of Computer Science, 
        Campus E 1.1, 66041 Saarbr\"ucken}
    \end{tabular}  
}

\begin{abstract}
Consider a random vector with finite second moments. If its  
precision matrix is an $M$-matrix, then all partial correlations
are non-negative. If that random vector is additionally Gaussian, the
corresponding Markov random field (GMRF) is called \emph{attractive}.\\
We study estimation of $M$-matrices taking the role of inverse second moment
or precision matrices using sign-constrained log-determinant divergence minimization.
We also treat the high-dimensional case with the number of variables exceeding the
sample size. The additional sign-constraints turn out to greatly
simplify the estimation problem: we provide evidence that explicit regularization is no longer
required. To solve the resulting convex optimization problem, we propose an
algorithm based on block coordinate descent, in which each sub-problem can 
be recast as non-negative least squares problem.\\ 
Illustrations on both simulated and real world data are provided.        
\end{abstract}


\end{frontmatter}


\section{Introduction}\label{sec:introduction}

The covariance matrix of a random vector and its inverse play
an important role in multivariate statistical analysis due to
their presence, for example, in PCA, discriminant analysis, tests of
hypotheses and confidence intervals for the mean. The sparsity pattern of the inverse
covariance, or synonymously precision matrix, reveals all pairwise conditional
independence relations in a Gaussian Markov random field (GMRF), which
correspond to missing edges in the resulting conditional independence graph 
\cite{Whittaker1990, Lauritzen1996, Rue2001}. The latter is central to 
Gaussian graphical modelling \cite{Dempster1972}, where one aims at a
parameter-parsimonious models in terms of a conditional independence graph
consisting of few edges. In recent years, precision matrix 
estimation and Gaussian graphical models (GGMs) have received considerable
attention in statistics, machine learning and optimization due to the prevalence    
of high-dimensional datasets in areas such as genomics, finance and neuroscience for which the
number of variables $p$ is comparable or even larger than the sample size
$n$. This setup has stimulated the development of various new inferential
procedures typically hinging on sparsity assumptions on the precision
matrix. One class of approaches tries to infer only its sparsity pattern, or
equivalently in the multivariate Gaussian case, the edges of the conditional
independence graph. For this purpose, procedures based on conditional
independence tests \cite{Kalisch2007, Anandkumar2012} and nodewise sparse
regression (\emph{neighbourhood selection} \cite{Mei2006, Zhou2011}) have been suggested. Given the graph structure, the precision matrix can be estimated subject to additional
constraints \cite{Drton2007, Dahl2005}. A second line of research    
is concerned with estimation of the precision matrix with the help of sparsity-promoting regularization schemes. 
In the references \cite{Yuan2007, Ban2008, Friedman2008, Rothman2008, Foygel2010,
  Ravikumar2011} $\ell_1$-regularized log-determinant divergence minimization,
which amounts to $\ell_1$-penalized maximum likelihood estimation in the
Gaussian case, is investigated. Related regularization schemes
enforcing sparsity of the off-diagonal elements of the precision matrix are proposed in \cite{Fan2009,
  Yuan2010, Cai2011, Shen2012}. In \cite{Duchi2009} and \cite{Honorio2009}, regularization schemes to enforce different forms of structured
sparsity are considered. In the present paper, we adopt the high-dimensional
setting that is addressed in the cited references, while considering elements
from a subcone of the positive semidefinite cone as target in
precision matrix estimation. Specifically, we consider positive definite
matrices that are symmetric $M$-matrices \cite{Ostrowski1937, BermanPlemmons}, i.e.~elements of the set 
\begin{equation}\label{eq:pdM}
\mc{M}^p = \{\Omega = (\omega_{jk}) \in \R^{p \times p}:\; \Omega \in \psd^p,
\; \; \,  \omega_{jk} \leq 0, \; \, j,k=1,\ldots,p, \; \, j \neq k \}, 
\end{equation}
where $\psd^p = \{\Omega \in \R^{p \times p}:\; \Omega = \Omega^{\T}, \; \,
\Omega \succ 0 \}$ denotes the set of symmetric, positive definite
matrices. In the statistics literature, B\o{}lviken \cite{Bolviken1982}
appears to be the first to consider matrices \eqref{eq:pdM} as precision
matrices whose partial correlations $-\omega_{jk}/\{\omega_{jj}
\omega_{kk} \}^{1/2}$, $j \neq k$, are all non-negative. Karlin and Rinott
\cite{Karlin1983} studied elements from \eqref{eq:pdM} as covariance or
precision matrices of a multivariate Gaussian distribution with a focus on 
total positivity \cite{Karlin1980}; see also \cite{Rinott2004}. In
\cite{Malioutov2006, Anandkumar2012}, the GMRF corresponding to a precision matrix of the form
\eqref{eq:pdM} is referred to as \emph{attractive} GMRF. In \cite{Malioutov2006},
attractive GMRFs are shown to be a sub-class of of \emph{non-frustrated}
GMRFs, which in turn form a sub-class of \emph{walk-summable} GMRFs. Statistical inference specifically for the class
\eqref{eq:pdM}, has, to the best of our knowledge, not been studied in
the literature. In \cite{LakeTenenbaum2010}, the authors consider
MAP estimation for the case that the precision matrix of a Gaussian random
vector belongs to the following subset of \eqref{eq:pdM}:
\begin{align}\label{eq:LpI}
\begin{split}
\mc{L}_I^p =  \Bigg\{\Omega \in &\R^{p \times p}:\; \Omega = \kappa I
- W + \text{diag}\left(\sum_{k = 1}^p w_{1p},\ldots, \sum_{k = 1}^p w_{pk}
\right), \\
&\kappa > 0, \; w_{jk} = w_{kj}, \; w_{jk} \geq 0, \; j,k=1,\ldots,p, \; \, j \neq k
\Bigg\} \subset \mc{M}^p, 
\end{split}
\end{align}
where the containment in $\mc{M}^p$ holds because by construction, all elements
of $\mc{L}_I^p$ are diagonally dominant and thus positive definite. We note
that \eqref{eq:LpI} equals the set of matrices that can be written as a
positive multiple of the identity plus the combinatorial Laplacian of an
undirected graph on $p$ vertices and positive edge weights $w_{jk}$,
$j,k=1,\ldots,p$. In \cite{LakeTenenbaum2010}, an exponential prior for the 
weights is proposed, so that MAP estimation amounts to $\ell_1$-penalized
maximum likelihood estimation. As discussed in more detail below, restricting
the class of admissible precision matrices by imposing sign-constraints on the
off-diagonals as in \eqref{eq:pdM} and \eqref{eq:LpI} can be a blessing and a curse                                 
at the same time. On the negative side, the requirement that all partial
non-negative correlations be non-negative is realistically not fulfilled in
most contemporary datasets. For example, in gene expression analysis, genes
may have both up- and down-regulatory effects on other genes. It is
a priori unclear what the consequences of estimation under model
misspecification are (see Section \ref{sec:misspecification} below). On the positive side, we show
that the presence of the additional sign constraints suffices to establish
existence and uniqueness of maximum likelihood estimation in the Gaussian case
even in a high-dimensional regime ($n < p$), which is unlike the unconstrained
case. Furthermore, we present empirical evidence that explicit regularization is not required and 
that subsequent thresholding of the off-diagonal entries of the constrained estimate yields a simple yet effective
procedure to recover the sparsity pattern of an underlying sparse target from
the class \eqref{eq:pdM} and hence also the structure of the associated
graph. This is akin to recent work on (thresholded) non-negative least squares in high-dimensional 
sparse regression \cite{SlawskiHein2011nips, Meinshausen2013}. Absence of tuning parameters and the tendency to produce sparse
solutions make the approach attractive in exploratory data analysis when the
goal is to find a sparse graph depicting positive dependence relations among
variables. In \cite{LakeTenenbaum2010}, learning taxonomies is presented as an
example where only positive dependence are of interest. In the present paper,
we also discuss a possible application to the analysis of landmark data
similar in spirit to \cite{Gu07}.
\paragraph{Outline} In Section \ref{sec:main}, we study central properties of 
sign-constrained log-determinant divergence minimization and positive
definite $M$-matrices. Sparse estimation based on thresholding is subsequently
discussed in Section \ref{sec:sparsification}. In Section \ref{sec:optimization}, we develop a block coordinate descent
algorithm to solve the resulting convex optimization problem and prove
its convergence. An extensive empirical study including the analysis of
real world datasets is presented in Section \ref{sec:experiments}. We conclude with a short summary. The appendix
contains all proofs.  

\paragraph{Notation}

Matrices are denoted by uppercase Latin or Greek letters and its elements
by the corresponding lowercase letters. For matrices starting with the letter $\Omega$,
the letter $\Sigma$ is used for their inverses. We use double subscripts to denote submatrices, i.e.~$A_{IJ}$ is the submatrix of some matrix $A$ 
with row indices $I$ and column indices $J$.  We write $A_{II}^{-1} = (A_{II})^{-1}$ for the inverse of a square invertible sub-matrix $A_{II}$ of $A$. A superscript $^c$ denotes the set complement.  
We will frequently arrange a symmetric and invertible matrix $A \in \R^{p
  \times p}$ in the following way. For $j \in
\{1,\ldots,p \}$ arbitrary, let $a_j \in \R^{p-1}$ be the vector with
components $a_{jk}, \; k=1,\ldots,p, \, k \neq j$, and $A_{jj} \in \R^{(p-1) \times (p-1)}$ the square submatrix of $A$ having
entries $\{a_{lm}, \; l,m \neq j \}$. After row and column permutations, $A$
and accordingly its inverse $B$ can be partitioned as
\begin{equation}\label{eq:partitioning}
 \left[ \begin{array}{cc}
           a_{jj}    & a_j^{\T} \\
           a_j      &  A_{jj}
          \end{array} \right], \quad \text{respectively} \quad \left[ \begin{array}{cc}
b_{jj} & b_j^{\T} \\
b_j   & B_{jj} 
\end{array} \right],
\end{equation}
where $b_{jj}$, $B_{jj}$ and $b_j$ are given by
\begin{align}\label{eq:partitioning_schur}
b_{jj} = \frac{1}{a_{jj} - a_j^{\T} A_{jj}^{-1} a_j}, \quad B_{jj} =
\left(A_{jj} - \frac{a_j a_j^{\T}}{a_{jj}} \right)^{-1}, \quad
b_j = -b_{jj} \cdot A_{jj}^{-1} a_j. 
\end{align}
For a square matrix $A$, $D(A)$ denotes the matrix resulting from $A$ after
setting all off-diagonal entries to zero. Likewise, $\off(A)$ denotes the matrix
resulting from $A$ after setting all diagonal elements to zero. Moreover, $\tr(A)$ denotes the trace of $A$. For square matrices $A_1,\ldots,A_K$, $\text{bdiag}(A_1,\ldots,A_K)$ denotes the block diagonal matrix
composed of these matrices. We write $A \succ 0$ and $A \gec 0$ for a positive definite respectively positive semidefinite matrix $A$,
whereas $\leq, \geq$, $<,>$ are used to denote component-wise inequalities, e.g.~$A \leq B$ means
that $a_{jk} \leq b_{jk}$ for all $j$ and $k$. The
symbols $I$ and $\bm{1}$ are used to denote identity matrices and vectors of ones, respectively.

\section{Positive definite $M$-matrices and sign-constrained log-determinant divergence minimization}\label{sec:main}

\subsection{Problem formulation}\label{sec:problem}
Let $x_1, \ldots, x_n$ be a sample of $n$ i.i.d.~realizations from a
multivariate Gaussian random vector $X$ with mean $\mu_* \in \R^p$ and covariance $\Sigma_* =
(\sigma_{jk}^*) \in \psd^p$ and precision matrix $\Sigma_*^{-1} = \Omega_* =
(\omega_{jk}^*)$. Assuming that $\mu_*$ is known and that $\Omega_* \in \mc{M}^p$ as
defined in \eqref{eq:pdM}, constrained maximum likelihood estimation of $\Omega_*$      
leads to the minimization problem
\begin{equation}\label{eq:logdet}
\min_{\Omega \in \mc{M}^p} -\log \det(\Omega) + \tr(\Omega S), \; \; \;
\text{where} \; \, S = \frac{1}{n} \su (x_i - \mu_*)(x_i - \mu_*)^{\T}.
\end{equation}
In case that $\Omega_* \notin \mc{M}^p$ or $X$ is non-Gaussian,
\eqref{eq:logdet} can be understood as M-estimation based on minimizing the
Bregman divergence between positive definite matrices that is induced by the
function $\Omega \mapsto -\log \det(\Omega)$,
cf.~\cite{Dhillon2007, Ravikumar2011} and \eqref{eq:bregmandivergence} below. Accordingly, we will henceforth refer to  
\eqref{eq:logdet} and related problems as (constrained or regularized)
log-determinant divergence minimization. It is well known that if $\mc{M}^p$ in
\eqref{eq:logdet} is replaced by $\psd^p$, i.e.~if the additional
sign-constraints on the off-diagonal elements are omitted, and $p > n$,
a minimizer of \eqref{eq:logdet} in general does not exist since the minimum in       
\eqref{eq:logdet} is not finite. Hence, it is a priori unclear whether the 
minimization problem \eqref{eq:logdet} is well-defined in the case $p > n$. As
stated in the following theorem, the additional constraint $\Omega \in \mc{M}^p$
makes a drastic difference.
\begin{theorem}\label{theo:existenceanduniqueness}
Consider the optimization problem \eqref{eq:logdet} and suppose that $S =
(s_{jk})$ has strictly positive diagonal elements. Then, unless there exists $(j,k), j \neq k$, such that  
$s_{jk} = \sqrt{s_{jj} s_{kk}}$, a minimizer of \eqref{eq:logdet} exists and is unique. 
\end{theorem}
In other words, unless there exists a pair of variables of perfect positive
sample correlation, which can easily be checked in practice, the constrained log-determinant divergence minimization
problem \eqref{eq:logdet} is well-posed, even though one may have $n < p$ and
no additional regularization is employed. We note that the conditions of
Theorem \ref{theo:existenceanduniqueness} are mild, because they are fulfilled
with probability one provided $n > 1$ and the random vector $X$ has a
distribution that is absolutely continuous w.r.t.~the Lebesgue measure.

\subsection{Optimality conditions and dual problem}\label{sec:dual}   

Within the present subsection, we study problem \eqref{eq:logdet} from the
point of view of convex optimization. It is standard to extend 
the negative log-determinant to the entire positive semidefinite cone
$\overline{\psd^p} = \{\Omega \in \R^{p \times p}:\;\Omega = \Omega^{\T}, \;
\,\Omega \gec 0  \}$ by setting $-\log \det(\Omega) = +\infty$ if
$\Omega \notin \psd^p $. Accordingly, we define $\overline{\mc{M}^p}$ as the
subset of matrices in $\overline{\psd^p}$ having only non-positive off-diagonal
elements. We may then re-write \eqref{eq:logdet} as  
\begin{equation}\label{eq:logdetext}
\min_{\Omega \in \overline{\mc{M}^p}} -\log \det(\Omega) + \tr(\Omega S),
\end{equation}
which constitutes a convex optimization problem. In fact, the constraint set, as the intersection of two convex cones, is a convex
cone, and the negative log-determinant is convex on $\overline{\psd^p}$, 
cf.~\cite{BoydVandenberghe2004}. The Lagrangian for \eqref{eq:logdetext} is 
given by 
\begin{equation}\label{eq:lagrangian}
L(\Omega, \Gamma) = -\log \det(\Omega) + \tr(\Omega S) + \tr(\Omega \Gamma),
\end{equation}  
where $\Gamma = (\gamma_{jk})$ is a symmetric, non-negative matrix of
Lagrangian multipliers with all diagonal entries being zero. By the
Karush-Kuhn-Tucker (KKT) optimality conditions, $(\wh{\Omega}, \wh{\Gamma})$
is an optimal solution if and only if
\begin{align}\label{eq:kkt}
S + \wh{\Gamma} = \wh{\Omega}^{-1}, \;\;\, \tr(\wh{\Omega} \wh{\Gamma}) = 0,
\;\;\, \wh{\Omega} \in \mc{M}^p, \; \; \, \wh{\Gamma} \in \R_+^{p \times p}, \,
\wh{\Gamma} = \wh{\Gamma}^{\T}, \,\text{diag}(\wh{\Gamma}) = 0.  
\end{align}
Note that under the stated conditions
\begin{equation}\label{eq:compslackness_componentwise}
\tr(\wh{\Omega} \wh{\Gamma}) = 0 \; \Longleftrightarrow \, \wh{\omega}_{jk} \wh{\gamma}_{jk} = 0, \; j,k=1,\ldots,p. 
\end{equation}
Convex duality yields 
\begin{align}\label{eq:dualproblem}
\min_{\Omega \in \overline{\mc{M}^p}} -\log \det(\Omega) + \tr(\Omega S) \notag
&= \max_{\Gamma \geq 0, \, \text{diag}(\Gamma) = 0} \; \min_{\Omega \in \overline{\mc{M^p}}} L(\Omega, \Gamma) \notag \\
&= \max_{S + \Gamma \in \overline{\psd^p}, \, \Gamma \geq 0, \, \text{diag}(\Gamma) = 0} \log \det(S + \Gamma) + p \notag \\
&= \max\limits_{\Sigma \in \overline{\psd^p}, \, \Sigma \geq S, \;
  \text{diag}(\Sigma) = \text{diag}(S)} \log \det(\Sigma) + p.
\end{align}
The second identity follows after taking the derivative of 
$L$ w.r.t.~$\Omega$, setting the result to zero and substituting this relation back into $L$ (cf.~the first condition in \eqref{eq:kkt}), while the third equality is by a change of variables. In other words, in the problem dual to \eqref{eq:lagrangian}, one seeks for
a positive definite matrix of maximum determinant, which dominates $S$
entry-wise and has the same diagonal entries. As corollary, we obtain the 
following characterization of \emph{inverse} positive definite $M$-matrices.
\begin{theo}\label{theo:characterization} 
$\Sigma \in \psd^p$ is an inverse $M$-matrix if and only if 
\begin{equation*}
\argmax\limits_{\Sigma' \in \overline{\psd^p}, \; \; \Sigma' \geq \Sigma, \;
  \diag(\Sigma') = \diag(\Sigma)} \log \det(\Sigma') = \Sigma.
\end{equation*}
\end{theo}
From \eqref{eq:dualproblem}, we can also read off necessity of the condition in Theorem
\ref{theo:existenceanduniqueness}: if there is a pair of variables of perfect
positive sample correlation, it is not possible to find $\Sigma \geq S$ with $\text{diag}(\Sigma) = \text{diag}(S)$ that
is strictly positive definite, so that \eqref{eq:dualproblem}
is unbounded from below.

\subsection{The class $\mc{M}^p$ as a model of multivariate dependence:
  restrictions and consequences of mis-specification}\label{sec:misspecification}
Recapitulating facts from \cite{BermanPlemmons, Karlin1983}, we will see that the
constraint $\Omega \in \mc{M}^p$ induces a rather specific model of multivariate
dependence for an underlying random vector $X = (X_j)_{j=1}^p$. Consequently, the target $\Omega_*$ can in general not
be expected to satisfy the given constraint. It is therefore of interest to know how
sign-constrained log-determinant divergence minimization \eqref{eq:logdetext} behaves
under model mis-specification, and we will investigate this issue for selected
examples.
\paragraph{$M$-matrices as precision matrices} Let $\Omega \in \mc{M}^p$. Then,
it is not hard to see that $\Omega = \delta I - B$ for symmetric $B \in \R_+^{p \times
  p}$ and $\delta \in (\lambda_1(B), \infty)$, where $\lambda_1$ denotes the
largest eigenvalue of $B$ (cf.~Appendix A). Expressing the inverse by a
Neumann series, one obtains $\Sigma = \Omega^{-1} = \delta^{-1} \sum_{k = 0}^{\infty}
(B/\delta)^k$ and hence $\Sigma \in \R_+^{p \times p}$. That is, a precision
matrix with non-positive off-diagonal entries implies non-negative marginal
correlations, i.e.~$\cov(X_{j}, X_k) \geq 0$ for all $j,k$. More generally,
for any pair of variables $(X_j, X_k)$ and any set of conditioning variables
$(X_l)_{l \in L}, \; L \subseteq \{1,\ldots,p \} \setminus \{j,k \}$,
the partial correlation of $(X_j, X_k)$ conditional on $X_L$ is
non-negative. This follows from the fact that covariances conditional on $L$
are given by the Schur complement of $\Sigma$ w.r.t.~$L$, that is
\begin{equation}\label{eq:Sigma_Schur}
\Sigma_{L^c L^c} - \Sigma_{L^c L} \Sigma_{LL}^{-1} \Sigma_{L L^c} =
\Omega_{L^c L^c}^{-1},     
\end{equation}
where the right hand side results by using partitioned inverses. Since
$\Omega$ is an $M$-matrix, so must be the sub-matrix $\Omega_{L^c L^c}$
and the claim follows from the same argument as above. Exchanging roles
of $\Sigma$ and $\Omega$ in \eqref{eq:Sigma_Schur}, we find that 
\begin{equation}\label{eq:Omega_Schur}
\Omega_{L^c L^c} - \Omega_{L^c L} \Omega_{LL}^{-1} \Omega_{L L^c} =
\Sigma_{L^c L^c}^{-1},
\end{equation}  
the Schur complement of $\Omega$ w.r.t.~$L$, is an $M$-matrix, or equivalently,
that the principal sub-matrix $\Sigma_{L^c L^c}$ of $\Sigma$ is an inverse $M$-matrix. This observation
implies that the sign of the partial correlations remain unchanged when confining oneself to
any subset of variables. To verify this, note that $\Omega_{LL} \in \mc{M}^{|L|}$ implies that $\Omega_{LL}^{-1} \geq 0$.
Combining this with $\Omega_{L^c L} \leq 0$, we have that $\Omega_{L^c L} \Omega_{LL}^{-1} \Omega_{L L^c} \geq 0$,
and in turn that the off-diagonal entries of the Schur complement are non-positive.\\
Finally, we remark that the regression coefficients of a linear regression for any variable $X_j$ on the
remaining variables are non-negative. To see this, first partition $\Sigma$ and $\Omega$ 
as in \eqref{eq:partitioning}:
\begin{equation*}
\left[ \begin{array}{cc}
           \sigma_{jj}    & \sigma_j^{\T} \\
           \sigma_j      &  \Sigma_{jj}
          \end{array} \right], \quad \text{and} \quad \left[ \begin{array}{cc}
\omega_{jj} & \omega_j^{\T} \\
\omega_j   & \Omega_{jj} 
\end{array} \right].
\end{equation*}
The regression coefficients equal $\Sigma_{jj}^{-1} \sigma_j =
-\omega_j/\omega_{jj} \in \R_{+}^{p-1}$ in view of \eqref{eq:partitioning_schur}, which are non-negative
as $\Omega$ has only non-positive off-diagonal elements. 


\begin{figure}
\normalsize
\begin{minipage}{0.78\textwidth}
\begin{tiny}
\begin{equation*}
\bordermatrix{S & \text{mec} & \text{vec} & \text{alg} & \text{ana} & \text{stat} \cr
                \text{mec}&  1   & .553   & .547  & .409 & .389\cr
                \text{vec} & .553 & 1      &  .610 & .485 & .436 \cr
                \text{alg} & .547 & .610   &  1    & .711 & .665 \cr 
                \text{ana} & .409 & .485   &  .711 & 1    & .607 \cr 
                \text{stat}& .389 & .436   &  .665 & .607 & 1} 
\bordermatrix{S^{-1} & \text{mec} & \text{vec} & \text{alg} & \text{ana} & \text{stat} \cr
                \text{mec}&   1.604   & -.560   & -.509  & \bm{.003} & -.043 \cr
                \text{vec} &  -.560 &   1.802   & -.657   & -.155 & -.038 \cr
                \text{alg} &  -.509 &   -.657 &   3.043   & -1.112 &  -.862 \cr 
                \text{ana} &  \bm{.003} & -.155   & -1.112  &  2.178    & -.517 \cr 
                \text{stat}&   -.043 &  -.038  &  -.862 &  -.517 & 1.921} 
\end{equation*}
\begin{equation*}
\bordermatrix{\wh{\Sigma} & \text{mec} & \text{vec} & \text{alg} & \text{ana} & \text{stat} \cr
                \text{mec}&  1   & .553   & .547  & .410 & .389\cr
                \text{vec} & .553 & 1      &  .610 & .485 & .436 \cr
                \text{alg} & .547 & .610   &  1    & .711 & .665 \cr 
                \text{ana} & .410 & .485   &  .711 & 1    & .607 \cr 
                \text{stat}& .389 & .436   &  .665 & .607 & 1} 
\bordermatrix{\wh{\Omega} & \text{mec} & \text{vec} & \text{alg} & \text{ana} & \text{stat} \cr
                \text{mec}&   1.604   & -.559   & -.508  & \bm{0} & -.042 \cr
                \text{vec} &  -.559 &   1.802   & -.658   & -.154 & -.038 \cr
                \text{alg} &  -.508 &   -.658 &   3.042   & -1.111 &  -.862 \cr 
                \text{ana} &  \bm{0} & -.154   & -1.111  &  2.178    & -.517 \cr 
                \text{stat}&   -.042 &  -.038  &  -.862 &  -.517 & 1.920} 
\end{equation*}
\end{tiny}
\end{minipage}
\begin{minipage}{0.20\textwidth}
\begin{flushright}
\includegraphics[height =
  0.14\textheight]{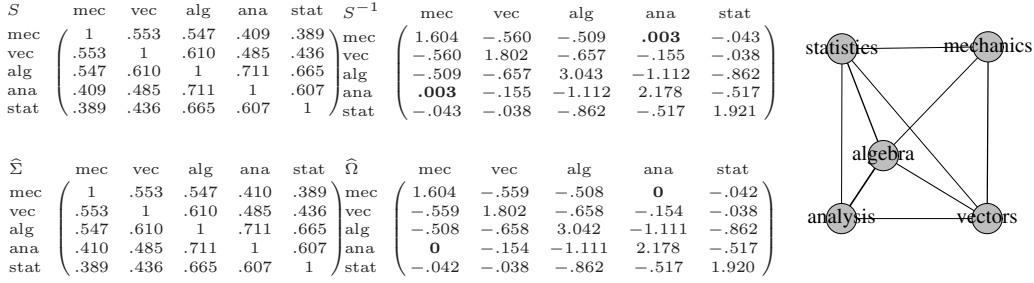}
\end{flushright}
\end{minipage}
{\normalsize
\caption{{\normalsize Solution pair $(\wh{\Sigma}, \wh{\Omega})$ for the math marks dataset
with correlation matrix $S$. The edge widths of the graph correspond to
partial correlations inferred from $\wh{\Omega}$.}}\label{fig:mathmarks}}
\end{figure}

\paragraph{Example: Math marks data}
The math marks dataset \cite{Mardia1979, Whittaker1990} contains the marks of
$n = 88$ students achieved in $p = 5$ subjects of mathematics ('mechanics',
'vectors', 'algebra', 'analysis', 'statistics'). Figure \ref{fig:mathmarks}
shows the resulting sample correlation matrix, its inverse and the solution 
pair $(\wh{\Sigma}, \wh{\Omega})$ of the minimization problem
\eqref{eq:logdetext} obtained with the computational approach described in
Section \ref{sec:optimization} below. The $M$-matrix model for the precision matrix appears to be
adequate, with $(\wh{\Sigma}, \wh{\Omega})$ closely matching $(S,
S^{-1})$. There is only one pair of variables that yields a tiny positive
off-diagonal entry in $S^{-1}$, which equals exactly zero in
$\wh{\Omega}$. From a practical point of view, the good fit of the $M$-matrix model 
suggests that if a student's performance in a subset of disciplines increase,
then so does the performance the remaining disciplines, i.e.~it least remains constant, but
there is no drop in performance. 

\begin{figure}
\begin{flushleft}
\begin{tabular}{cccc}
\hspace{-0.6cm} (1) & (2) & (3) & (4) \\
\hspace{-0.6cm}\includegraphics[height = 0.07\textheight]{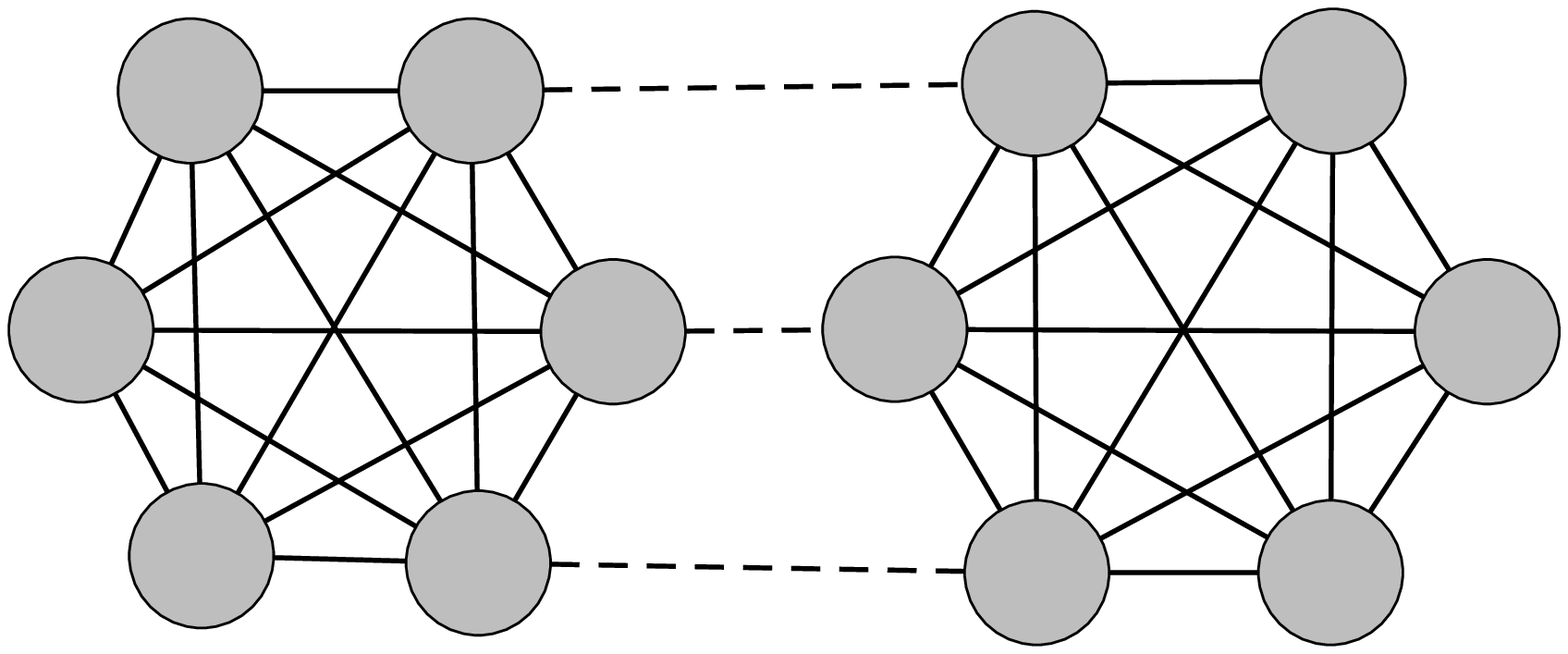} & \hspace{-0.1cm}\includegraphics[height = 0.022\textheight]{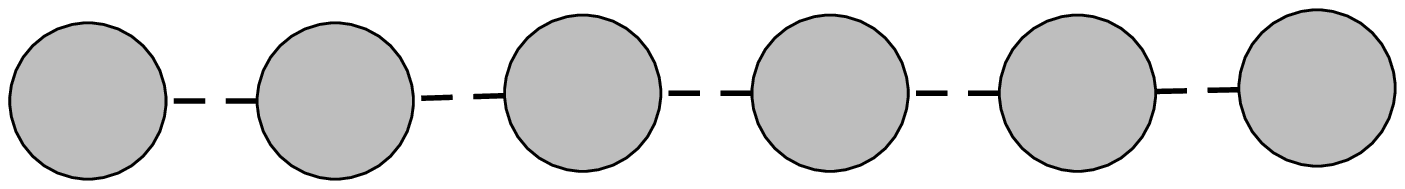}   & \hspace{-0.3cm} \includegraphics[height = 0.047\textheight]{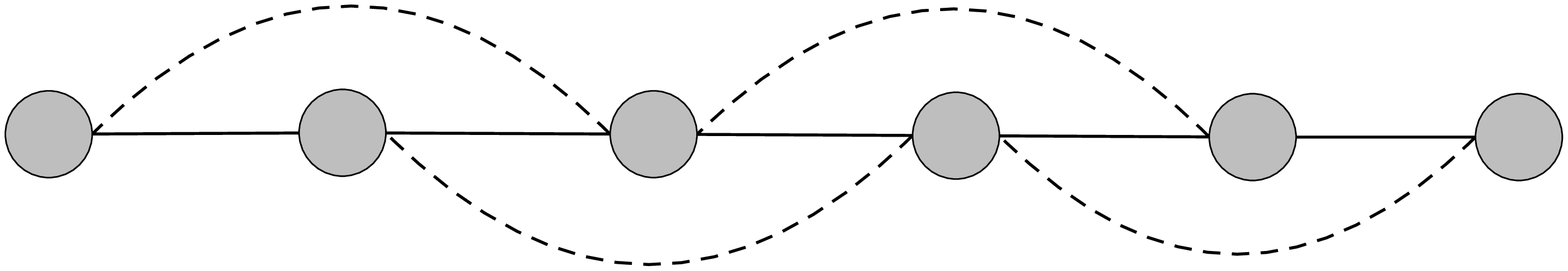}  &  \hspace{-0.3cm}\includegraphics[height = 0.088\textheight]{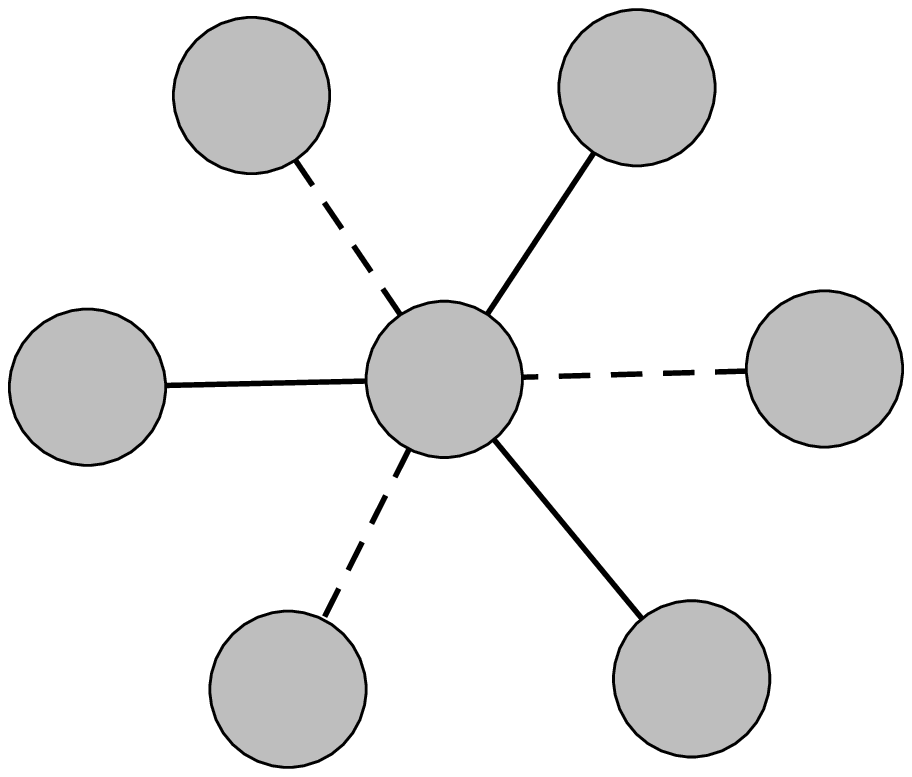}  \\
\hspace{-0.6cm}\includegraphics[height = 0.07\textheight]{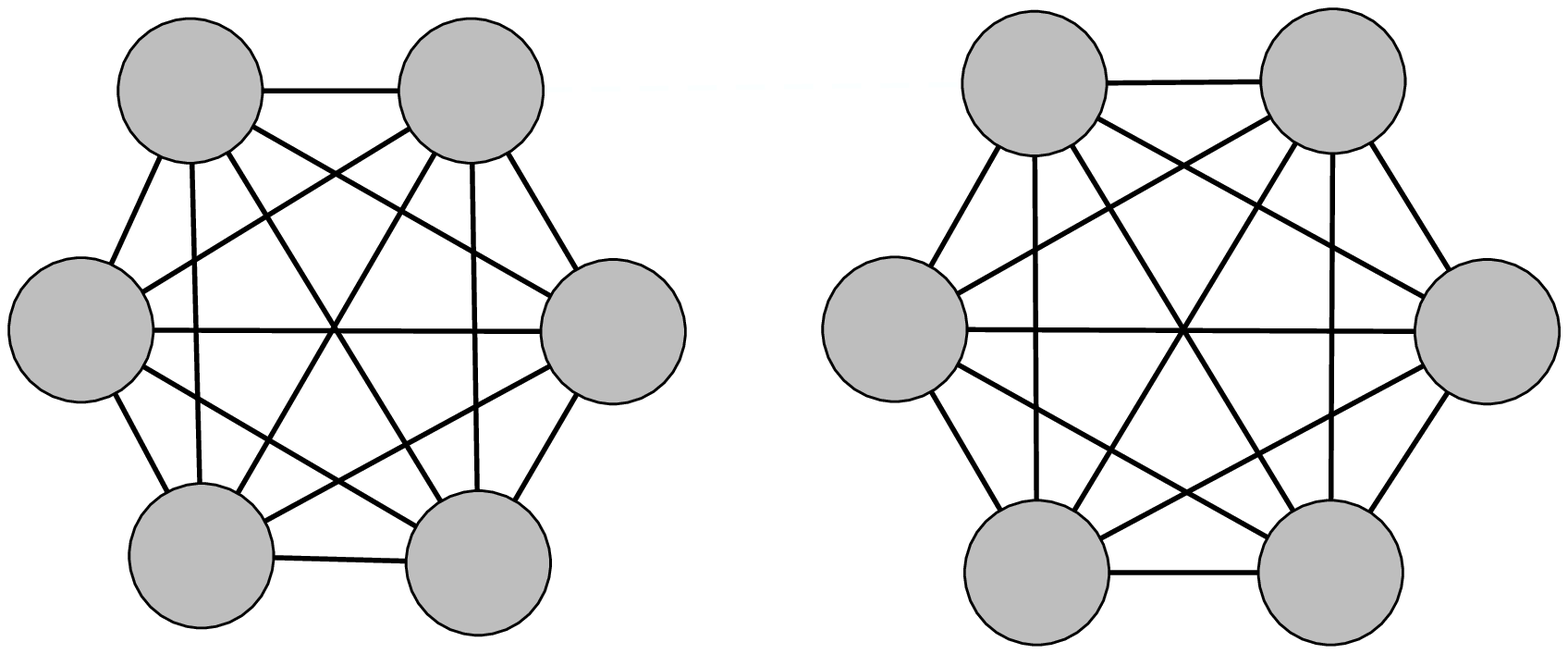}  & \hspace{-0.1cm}\includegraphics[height = 0.054\textheight]{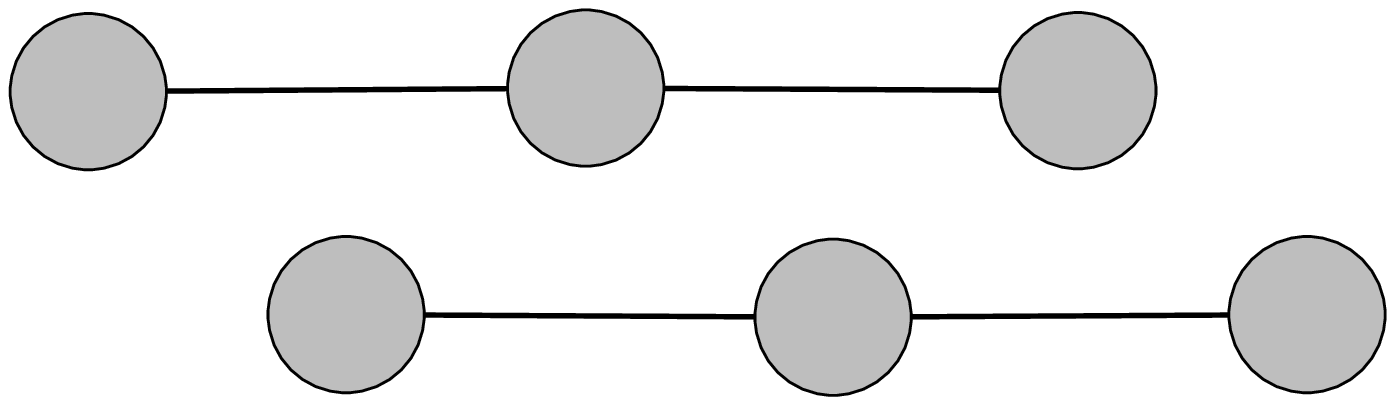}  & \hspace{-0.3cm}\includegraphics[height = 0.047\textheight]{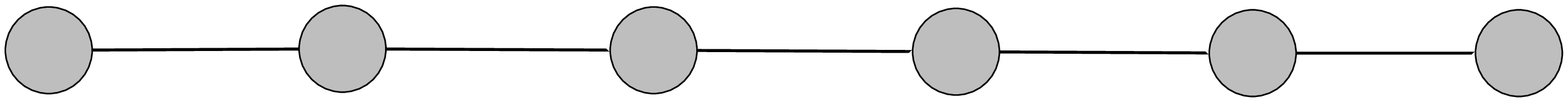}  &   \hspace{-0.3cm}\includegraphics[height = 0.088\textheight]{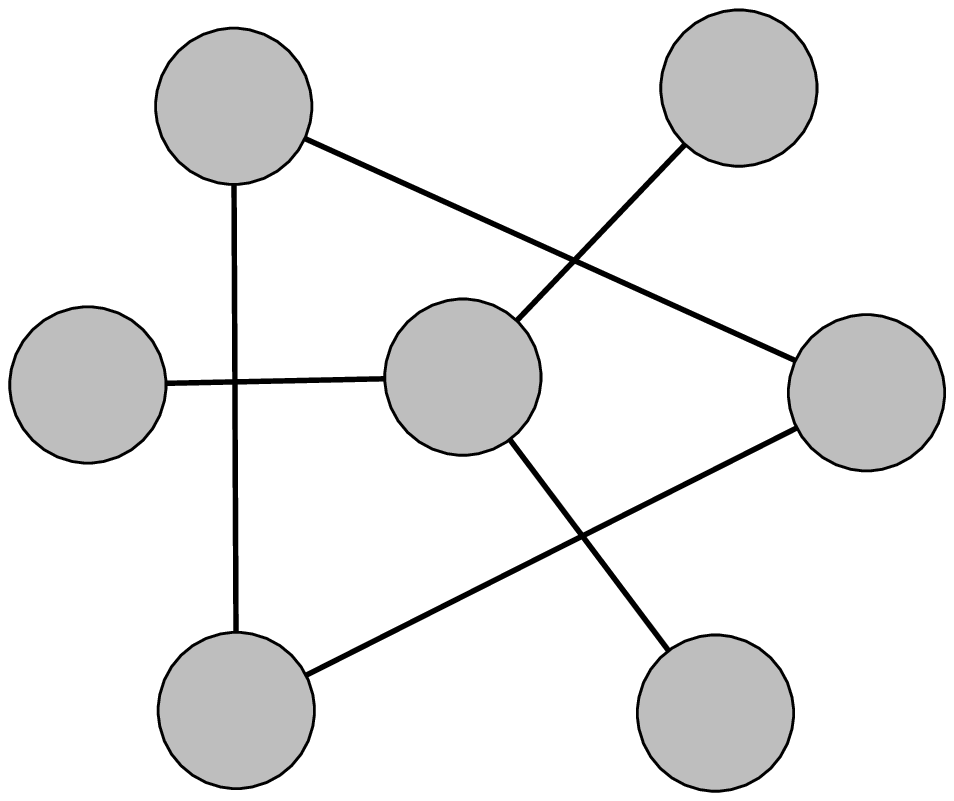} \\
\end{tabular}
\end{flushleft}              
\caption{Graph-based description of the examples (1) to (4). Dashed edges correspond to positive off-diagonal entries and solid edges to negative ones. Top: graph structure of $\Omega_*$. Bottom: graph structure of $\Omega_{\bullet}$.}\label{fig:misspec}
\end{figure}

\paragraph{Model mis-specification} Since the constraint $\Omega \in \mc{M}^p$ is rather strong, it is important 
to have some understanding about how sign-constrained log-determinant divergence minimization \eqref{eq:logdetext} behaves
under mis-specification, i.e.~the population precision matrix $\Omega_* \notin \mc{M}^p$. In this case, the estimator
$\wh{\Omega}$ may be subject to a substantial bias. In the following, we discuss this issue at
the population level. To this end, we define
\begin{align}
\Omega_{\bullet} &= \argmin_{\Omega \in \overline{\mc{M}^p}} -\log \det(\Omega) + \tr(\Omega \Sigma_*) \label{eq:logdet_population_Omega} \\
             &= \argmin_{\Omega \in \overline{\mc{M}^p}} -\log \det(\Omega) + \log \det(\Omega_*) + \tr((\Omega - \Omega_*) \Sigma_*) \notag \\ 
             &\invcoloneq \argmin_{\Omega \in \overline{\mc{M}^p}} D(\Omega_* \parallel \Omega), \label{eq:bregmandivergence} 
\end{align} 
where $D(\Omega_* \parallel \Omega)$ denotes the Bregman divergence of $\Omega$ from $\Omega_*$, which coincides
(apart from a factor of $1/2$) with the Kullback-Leibler divergence of two zero-mean Gaussian distributions
with precision matrices $\Omega$ and $\Omega_*$, respectively. Note that according to the dual problem 
\eqref{eq:dualproblem}, $\Sigma_{\bullet} = \Omega_{\bullet}^{-1}$ satisfies 
\begin{equation}\label{eq:logdet_population_Sigma}
\Sigma_{\bullet} = \argmax\limits_{\Sigma \in \overline{\psd^p}, \; \; \Sigma \geq \Sigma_*, \; \text{diag}(\Sigma) = \text{diag}(\Sigma_*)} \log \det(\Sigma). 
\end{equation}  
Depending on the degree of mis-specification, it may be possible that $\Omega_{\bullet}$ preserves positive partial correlations in $\Omega_*$. Ideally,       
one has for all $j \neq k$
\begin{equation}\label{eq:misspec_ideal}
\omega_{jk}^{\bullet} = \omega_{jk}^* \quad \text{if} \; \, \omega_{jk}^* < 0, \quad \text{and} \;\,  \omega_{jk}^{\bullet} = 0 \; \; \, \text{otherwise.}  
\end{equation}
That is, entries matching the sign constraints are maintained, while negative
partial correlations are zeroed out. If interest is only in the identification of pairs of variables of positive partial correlation and 
if \eqref{eq:misspec_ideal} holds, then moving from $\Omega_{*}$ to $\Omega_{\bullet}$ does not result in a loss of information. For example, in recommender systems, one is 
interested in finding pairs of items where the purchase of one item increases the chances of purchasing
the other one. Several more examples are presented in Section \ref{sec:experiments}. Beyond the ideal case \eqref{eq:misspec_ideal}, one 
can ask whether $\Omega_{\bullet}$ at least preserves the pairs of negative sign, i.e.~whether it holds that
\begin{equation}\label{eq:misspec_signpreserve}
 \mc{E}^* = \{(j,k):\; \omega_{jk}^* < 0 \}  = \mc{E}^{\bullet} = \{(j,k):\; \omega_{jk}^{\bullet} < 0 \}.      
\end{equation}
Below, we shed some light on this question for specific choices of $(\Omega_*, \Sigma_*)$ for which it is possible to compute the corresponding 
solution pair $(\Omega_{\bullet}, \Sigma_{\bullet})$ in closed form.\\
%
\emph{(1) Block structure.} Let $\Omega_* \in \psd^p$ be partitioned as
\begin{equation*}
\Omega_* = \left[ \begin{array}{cc}
                  \Omega_{*,11} & \Omega_{*,12} \\
                  \Omega_{*,21} & \Omega_{*,22} \\ 
                \end{array} \right], \; \Omega_{*,11} \in \psd^{p_1}, \; \, \Omega_{*,22} \in \psd^{p_2}, \; \, \Omega_{*,12} =  \Omega_{*,21}^{\T} \in \R^{p_2 \times p_1},
\end{equation*}
and let $\Sigma_{*,11}$ etc.~be defined accordingly.\\ 
\emph{(a).} If $\Omega_{*,11} \in \mc{M}^{p_1}$, $\Sigma_{*,22} \in \mc{M}^{p_2}$ and $\Omega_{*,12} = 0$, then 
$\Omega_{\bullet}$ in \eqref{eq:logdet_population_Omega} is given by $\Omega_{\bullet} = \text{bdiag} \left(\Omega_{*,11}, \{D(\Sigma_{*,22})\}^{-1} \right)$, where $D(\Sigma_{*,22})$ is the restriction of $\Sigma_{*,22}$ to its diagonal. 
In order to verify this, let us consider the dual in \eqref{eq:logdet_population_Sigma}. We have 
$\det(\Sigma_*) = \det(\Sigma_{*,11}) \cdot \det(\Sigma_{*,22}) \leq \det(\Sigma_{*,11}) \det(D(\Sigma_{*,22}))$ by 
Hadamard's inequality. Since we have $D(\Sigma_{*,22}) \geq \Sigma_{*,22} \in \mc{M}^{p_2}$,  $\Sigma_{\bullet} = \text{bdiag}(\Sigma_{*,11}, D(\Sigma_{*,22}))$
is dual feasible. Given the upper bound on $\det(\Sigma_*)$ and 
Theorem \ref{theo:characterization}, it must be the solution \eqref{eq:logdet_population_Sigma}. Observe that $\Omega_{\bullet}$ fulfills \eqref{eq:misspec_ideal}.\\
\emph{(b).} Let now $\Omega_{*,11} \in \mc{M}^{p_1}, \Omega_{*,22} \in \mc{M}^{p_2}$, $\Omega_{*,12} \geq 0$.  Then we have \\
$\Omega_{\bullet} = \text{bdiag}(\Sigma_{*,11}^{-1}, \Sigma_{*,22}^{-1})$. To see this, note that the partitioned inverse formula yields 
\begin{equation*}
\Sigma_{*,12} = -\Omega_{*,11}^{-1} \Omega_{*,12} (\Omega_{*,22} - \Omega_{*,21} \Omega_{*,11}^{-1} \Omega_{*,12})^{-1} = -\Omega_{*,11}^{-1} \Omega_{*,12} \Sigma_{*,22} \leq 0,             
\end{equation*} 
because $\Omega_{*,11}^{-1}$, $\Sigma_{*,22}$ and $\Omega_{*,12}$ have only non-negative entries. Consequently, $\Sigma_{\bullet} = \Omega_{\bullet}^{-1}$ is 
feasible for the dual in \eqref{eq:logdet_population_Sigma}. Feasibility of $\Omega_{\bullet}$ for the primal in \eqref{eq:logdet_population_Sigma}
follows from 
\begin{equation*}
\off(\Sigma_{*,11}^{-1}) = \off(\Omega_{*,11} - \Omega_{*,12} \Omega_{*,22}^{-1} \Omega_{*,21}) \leq 0  
\end{equation*}
and an according argument for $\off(\Sigma_{*,22}^{-1})$. From the Hadamard-Fischer inequality, we obtain
that $\det(\Sigma_{*}) \leq \det(\Sigma_{*,11}) \cdot \det(\Sigma_{*,22})$, and the claim follows from Theorem            
\ref{theo:characterization}. Observe that $\Omega_{\bullet}$ preserves signs according to \eqref{eq:misspec_signpreserve}.\\ 
\\
The next two examples deal with precision matrices corresponding to stationary autoregressive (AR) processes of orders 1 and 2.\\ 
\emph{(2) AR(1)-structure.} Let $\Sigma_*$ have entries $\sigma_{jk}^* = \rho^{|j-k|}$, $j,k=1,\ldots,p$ for $p$ even
and $-1 < \rho < 0$. Then $\sigma_{jk}^* < 0$ if one of $(j, k)$ is even and the other one is odd, and $\sigma_{jk}^* > 0$ if both 
of $(j,k)$ are even/odd. The inverse $\Omega_*$ is a band matrix of bandwidth one, the non-zero off-diagonal entries being all 
equal to $-\rho/(1 - \rho^2)$ (see e.g.~\cite{Rue2001}, p.2). Set $\Sigma_{\bullet} = \Sigma_* + \Gamma_{\bullet}$, where the entries of $\Gamma_{\bullet}$ are given by
$\gamma_{jk}^{\bullet} = 0$ if $(j,k)$ are both even/odd and $\gamma_{jk}^{\bullet} = -\sigma_{jk}^*$ otherwise. We will show that 
$\Sigma_{\bullet}$ and its inverse $\Omega_{\bullet}$ are the solutions \eqref{eq:logdet_population_Omega}/\eqref{eq:logdet_population_Sigma} by 
verifying the KKT optimality conditions as given in \eqref{eq:kkt}. First note that $\Sigma_{\bullet} \geq \Sigma_*$ and
$\Sigma_{\bullet} \in \psd^p$, where second claim can be seen from $\Pi \Sigma_{\bullet} \Pi^{\T}= \text{bdiag}(\wt{\Sigma}, \wt{\Sigma})$, where $\Pi$ is a permutation matrix permuting odd rows $1,3,5,\ldots$ on 
the first $1,\ldots,p/2$ rows and the even rows on the rows $p/2+1,p/2+2,\ldots,p$,
and $\wt{\Sigma} \in \R^{p/2 \times p/2}$ has entries $\wt{\sigma}_{lm} = \rho^{2|l-m|} \invcoloneq \phi^{|l-m|}$, $l,m=1,\ldots,p/2$, i.e.~$\wt{\Sigma}$ 
has AR(1)-structure with parameter $0 < \phi < 1$. The non-zero off-diagonal entries of $\wt{\Sigma}^{-1}$ are all equal 
to $-\phi/(1 - \phi^2)$ , hence $\Omega_{\bullet} \in \mc{M}^{p}$. Finally observe that $\Omega_{\bullet}$ satisfies 
the complementarity slackness condition $\tr(\Omega_{\bullet} \Gamma_{\bullet}) = 0$.\\
Note that \eqref{eq:misspec_signpreserve} is violated.\\ 
\emph{(3) AR(2)-structure.} Let $\Sigma_*$ have entries $\sigma_{jk}^* = \rho_{|j-k|}$, where
$\rho_{\ell}$, $\ell = 1,\ldots,p-1$, is defined by the recursion
\begin{equation}\label{eq:yulewalkerAR2}
\rho_0 = 1, \qquad \rho_1 = \frac{\phi_1}{1 - \phi_2}, \qquad \rho_{\ell} =
\phi_1 \rho_{\ell-1} + \phi_2 \rho_{\ell - 2}, \; \; \ell \geq 2. 
\end{equation}
for parameters $\phi_1$ and $\phi_2$ satisfying the stationarity condition
\begin{equation*}
\left|\frac{1}{2 \phi_2} \left(-\phi_1 + \sqrt{\phi_1^2 + 4 \phi_2} \right) \right| > 1, \quad \left|\frac{1}{2 \phi_2} \left(-\phi_1 - \sqrt{\phi_1^2 + 4 \phi_2} \right) \right| > 1,
\end{equation*}   
cf.~\cite{Chatfield2003}. The inverse $\Omega_*$ has bandwidth two \cite{Rue2001}. In the appendix, we 
prove that if $4 |\phi_2| < \phi_1$, the minimizer $\Omega_{\bullet}$ in \eqref{eq:logdet_population_Omega} preserves
signs according to \eqref{eq:misspec_signpreserve}. The condition $4 |\phi_2| < \phi_1$ is found to be tight in
the sense that if it fails, \eqref{eq:misspec_signpreserve} in general does not hold.\\   
\emph{(4) Star structure.} Let $\Sigma_* = \begin{pmatrix} 1 & -\rho^{\T} \\ -\rho  & I + \rho \rho^{\T} \end{pmatrix}$, where $\rho \in \R^{p-1}$. The inverse
results as $\Omega_* = \begin{pmatrix} 1 + \norm{\rho}_2^2 & \rho^{\T} \\ \rho  & I \end{pmatrix}$. In the appendix, it is
shown that the solutions \eqref{eq:logdet_population_Omega}/\eqref{eq:logdet_population_Sigma} are given by
\begin{align*}
\begin{split}
&\Omega_{\bullet} = \begin{pmatrix}
              1 + \nnorm{\wt{\rho}}_2^2  & \wt{\rho}^{\T} \\
              \wt{\rho} &  I - \frac{(\rho - \wt{\rho}) (\rho - \wt{\rho})^{\T}}{1 + \nnorm{\rho - \wt{\rho}}_2^2}     
              \end{pmatrix},
              \; \Sigma_{\bullet} = \begin{pmatrix}
              1 & -\wt{\rho}^{\T} \\
              -\wt{\rho} &  I + \rho \rho^{\T} - \rho \wt{\rho}^{\T} -
              \wt{\rho} \rho^{\T} + 2 \wt{\rho} \wt{\rho}^{\T} 
              \end{pmatrix},\\
&\text{where} \; \; \, \wt{\rho} = (\min(\rho_j, 0)).  
\end{split}
\end{align*}  
It follows that in general, sign preservation \eqref{eq:misspec_signpreserve} does not hold: if $\rho$ has at
least two positive components $j,k$, $j \neq k$, then $\rho_j - \wt{\rho}_j > 0$, $\rho_k - \wt{\rho}_k > 0$ and 
hence $(\rho_j - \wt{\rho}_j)(\rho_k - \wt{\rho}_k) > 0$, i.e.~the corresponding entry in the bottom right block 
of $\Omega_{\bullet}$ is negative.\\ 
Figure \ref{fig:misspec} provides a graphical description of the above four settings. In summary, it is not
guaranteed that $\Omega_{\bullet}$ recovers the set $\mc{E}^*$ \eqref{eq:misspec_signpreserve}. While there are instances of $\Omega_*$ where this is fulfilled, there are cases like the AR(1) example for which $\Omega_{\bullet}$ even satisfies an opposite sign constraint. It is important to bear in mind that for simplicity, we have limited our discussion to the population setting. For a complete treatment, one would need to additionally take into account the effect resulting from the replacement of $\Sigma_*$ by
the sample covariance $S$.     

\section{Sparsification}\label{sec:sparsification}

Sparsity has been connected to precision matrix estimation since
the seminal work of Dempster \cite{Dempster1972} on covariance selection. With the
advent of high-dimensional data analysis, sparsity has become a key concept both
to obtain interpretable results and to establish guarantees of various concrete
sparsity-promoting estimation techniques. Among these, penalty-based approaches
are most prominent \cite{Ban2008, Rothman2008, Fan2009}. It hence
appears natural to complement sign-constrained log-determinant divergence minimization \eqref{eq:logdetext} by such
a penalty when sparsity is desired. This prompts the following modification of \eqref{eq:logdetext}.   
\begin{equation}\label{eq:logdetext_penalty}
\min_{\Omega \in \overline{\mc{M}^p}} -\log \det(\Omega) + \tr(\Omega S) + \lambda \sum_{(j,k): j \neq k} \textsf{pen}(\omega_{jk}), \quad \lambda  \geq 0. 
\end{equation} 
Choosing the negative identity for the penalty $\textsf{pen}$, one ends up with a sign-constrained version of the graphical lasso. 
We here refrain from the penalty approach. Instead, we argue for a post-processing procedure
combining thresholding and re-fitting.

\subsection{Hard thresholding and re-fitting}
Suppose that the population precision matrix $\Omega_* \in \mc{M}^p$ has few non-zero entries, i.e.~the
set $\mc{E}^* = \{(j,k):\,\omega_{jk}^* < 0 \}$ as already defined in \eqref{eq:misspec_signpreserve} has
small cardinality. We may interpret $\mc{E}^*$ as the edge set of the graph of positive partial
correlations associated with $\Omega_*$. Our aims are recovery of $\mc{E}^*$ and accurate
estimation of $\Omega_*$ under sparsity given a finite sample of i.i.d.~observations. We suggest the following scheme.
\begin{enumerate}
\item We first compute the minimizer $\wh{\Omega}$ of \eqref{eq:logdetext} as initial estimate. 
\item We apply hard thresholding to the off-diagonal entries, i.e.~for some
      threshold $t \geq 0$, we let
      \begin{equation}\label{eq:Omega_thres}
       \wh{\Omega}(t) = (\wh{\omega}_{jk}(t)) = \begin{cases}
      \wh{\omega}_{jk} \quad &\text{if} \; j=k, \\
      I(-\wh{\omega}_{jk} > t) \, \wh{\omega}_{jk} \quad &\text{if} \; j \neq k.
\end{cases}
\end{equation}
\item We estimate $\mc{E}^*$ by $\wh{\mc{E}}(t) = \{(j,k):\,\wh{\omega}_{jk}(t) < 0\}$.
\end{enumerate}
Since $\wh{\Omega}(t)$ is in general no longer positive definite, we perform 
a re-fit subject to additional zero constraints on the off-diagonal entries
as represented by $\wh{\mc{E}}(t)$, that is we compute
\begin{equation}\label{eq:refitting}
\wh{\wh{\Omega}}(t) = \argmin_{\Omega \in \overline{\mc{M}^p}, \; \omega_{jk} = 0 \;
  \forall (j,k) \notin \wh{\mc{E}}(t), \, j \neq k} -\log \det(\Omega) + \tr(\Omega
S).
\end{equation}            
In addition to being positive definite, the final estimator $\wh{\wh{\Omega}}(t)$ potentially improves over $\wh{\Omega}$ 
regarding estimation of $\Omega_*$ if $\wh{\mc{E}}(t) \supseteq \mc{E}^*$ and $\wh{\mc{E}}(t)$ is of small cardinality; see
\cite{Zhou2011} for an analysis of precision matrix estimation subject to zero constraints on off-diagonal elements.    

\subsection{Justification}
In the sequel, we provide some theoretical underpinning for the suggested approach. Successful identification 
of $\mc{E}^*$ via thresholding entails that the initial estimate $\wh{\Omega}$ obeys the condition
\begin{equation}\label{eq:cond_suc_thresholding}
\max_{(j,k), \, j \neq k} |\wh{\omega}_{jk} - \omega_{jk}^*| \leq t, \; \quad \; \min_{(j,k) \in \mc{E}^*} |\omega_{jk}^*| > 2t.
\end{equation} 
We wish to take the threshold $t$ as small as possible so that it is possible to detect even non-zero off-diagonal elements in $\Omega_*$ 
of small absolute magnitude, but still large enough to filter out all pairs not in $\mc{E}^*$. Performance
hence depends on the $\ell_{\infty}$-distance of off-diagonal entries of $\wh{\Omega}$ to those of $\Omega_*$. In the classic
setting with $p$ fixed as $n \rightarrow \infty$, consistency of $\wh{\Omega}$ can be established by using standard arguments. 
\begin{prop}\label{prop:fixedp} Let $x_1,\ldots,x_n$ be i.i.d.~realizations from a $p$-dimensional random
vector $X$ with precision matrix $\Omega_*$. Furthermore, suppose that
$X$ has finite fourth moments. Then, with $\Omega_{\bullet}$ as defined in \eqref{eq:logdet_population_Omega}, the minimizer of the sign-constrained log-determinant divergence \eqref{eq:logdetext}
satisfies $\wh{\Omega} = \Omega_{\bullet} + o_{\p}(1)$ as $n \rightarrow \infty$ and $p$ stays fixed. 
\end{prop}
Note that if $\Omega_{*} \in \mc{M}^p$, we have $\wh{\Omega} = \Omega_* + o_{\p}(1)$. In that case, for $n$ large enough, there exists a threshold $t$ such that \eqref{eq:cond_suc_thresholding} holds. In the modern setting 
with $p$ being of the order of $n$ or even $p > n$, the situation is less clear. Here, the use of the thresholding
procedure is mainly justified by its strong empirical performance (cf.~Section \ref{sec:experiments}). In addition, we provide
the very first step towards an understanding of what one observes empirically. Namely, we can show the following.
\begin{prop}\label{prop:singleedge_recovery} For $\rho \in (0,1)$, let 
\begin{equation*}
\Sigma_* = \emph{\text{bdiag}} \left(\begin{pmatrix} 1 & \rho  \\
                                             \rho & 1  \end{pmatrix}, I_{p-2} \right),\;\;\, 
\Omega_* = \emph{\text{bdiag}}\left(\begin{pmatrix} \omega_{11}^* &  \omega_{12}^*  \\
                                   \omega_{12}^* & \omega_{22}^*  \end{pmatrix}, I_{p-2} \right),
\end{equation*}
with $\omega_{12}^* = -\rho/(1 - \rho^2)$ so that $\mc{E}^* = \{(1,2),(2,1) \}$ and $\omega_{11}^* = \omega_{22}^* = 1/(1 - \rho^2)$. Denote 
\begin{equation}\label{eq:supbound}
B \coloneq \max_{1 \leq j,k \leq p} |s_{jk} - \sigma_{jk}^*|,
\end{equation}
where $S = (s_{jk})$ is the sample covariance matrix \eqref{eq:logdet} based
on an i.i.d.~sample whose population covariance equals $\Sigma_*$. 
Suppose that sign-constrained log-determinant divergence minimization \eqref{eq:logdetext} has a 
unique minimizer $\wh{\Omega}$ and that $B = o_{\p}(1)$ as $n \rightarrow \infty$. Then there exists constants $c_1, c_2, c_3 > 0$ such that 
\begin{align*}
\max_{(j,k) \notin \mc{E}^*, \, j \neq k} (-\wh{\omega}_{jk}) \leq c_1B. \; \text{Moreover, if} \; \rho > c_2B, \; \, (-\wh{\omega}_{12}) \geq (-\omega_{12}^*) - c_3B > 0.
\end{align*}
as $n \rightarrow \infty$. Consequently, if $(-\omega_{12}^*)  - (c_3 + c_1)B > 0$, when applying thresholding with $t = c_1B$, it holds that $\wh{\mc{E}}(t) = \mc{E}^*$. 
\end{prop}
Proposition \ref{prop:singleedge_recovery} studies the special case in which the graph associated with $\Omega_*$ consists of a single edge. 
The ingredients of Proposition \ref{prop:singleedge_recovery} are existence of $\wh{\Omega}$ (which
may be concluded from Theorem \ref{theo:existenceanduniqueness}), and a uniform bound $B$ on the
differences of the entries of $S$ and those of $\Sigma_*$. Finite sample, high probability
upper bounds on $B$ in dependence on the decay of the tails of the underlying distribution can be found in \cite{Ravikumar2011}.
In particular, for sub-Gaussian tails, one has $B = O(\sqrt{\log(p)/n})$ with high probability, so that 
Proposition \ref{prop:singleedge_recovery} asserts identification of $\mc{E}^*$ even in the high-dimensional 
case, provided $|\omega_{12}^*|$ is not too small in the sense that it stands out of the effective noise level.\\
A result for general $\mc{E}^*$ is out of the scope of the present paper, though we conjecture that the rates
of convergence for the elementwise $\ell_{\infty}$-error \eqref{eq:cond_suc_thresholding} are at least comparable
to those available for $\ell_1$-regularization-based methods, e.g.~\cite{Ravikumar2011, Cai2011}. The proof
techniques employed therein exploit the presence of explicit regularization, whereas regularization induced by the constraint $\Omega \in \mc{M}^p$ is implicit and requires substantial extra work in order to be leveraged; see the proof of Proposition \ref{prop:singleedge_recovery} in \ref{app:singleedge_recovery}.        

\paragraph{Thresholding vs.~penalization} We now briefly explain why we prefer thresholding 
over the commonly used penalization approach \eqref{eq:logdetext_penalty}. First, sparsity-promoting 
penalty terms induce a bias, which adversely affects the identification of small, yet non-zero 
off-diagonal entries of $\Omega_*$. Second, thresholding is a direct way to achieve a desired
level of sparsity: when one aims at a sparsity level of $q$ percent, one simply keeps the
$q$ percent off-diagonal entries of $\wh{\Omega}$ that are largest in absolute magnitude, and
sets the rest equal to zero. The sparsity level is not as easy to control when using penalization,
because it cannot be read off directly what the resulting sparsity level for a specific choice 
of the regularization parameter will be. It is hence common to compute multiple solutions along
a grid of specific values \cite{Ban2008, Friedman2008}. However, this is unfavorable from
a computational point of view, because one has to solve multiple instances of a challenging 
convex optimization problem.

\section{Computational approach}\label{sec:optimization}

Within this section, we present our computational approach for solving the sign-
constrained log-determinant divergence minimization problem \eqref{eq:logdetext} numerically. As a convex
optimization problem, it can be handled by general purpose solvers like \textsf{CVX} \cite{cvx}. However,   
\textsf{CVX} becomes rather slow once the dimension $p$ crosses $30$. Note that the number of variables and 
constraints is quadratic in $p$, i.e.~with $p$ in the hundreds, the number of variables is in the ten thousands, and with 
$p$ in the thousands, the number of variables is in the millions. It is clear that standard 'off-the-shelf' interior point
methods as used in \textsf{CVX} are not suitable both with regard to runtime and memory requirements. We hence devise
a customized solver for the problem.

\subsection{Block coordinate descent}
 
The algorithm that we propose follows the one pioneered in \cite{Ban2008} to solve the graphical
lasso problem, and gives rise to an analogous interpretation: just like the approach in \cite{Ban2008}
amounts to recursively solving $\ell_1$-penalized regression (lasso) problems, our algorithm amounts
to recursively solving non-negative least squares regression (NNLS) problems. In view of this connection, 
optimization can be delegated to an arbitrary NNLS solver. Apart from conceptual simplicity and ease of
implementation, the algorithm has a solid theoretical foundation as block coordinate descent scheme, so
that existing theory can be leveraged to establish convergence. The approach solves problem up to 
$p=1000$ still reasonably fast, but comes with a rather sharp increase in runtime as $p$ increases, with
a complexity of $O(p^4)$.

\begin{algorithm}[b!]
\caption{Algorithm for problem \eqref{eq:logdetext}}\label{alg:blockdescent}
\begin{algorithmic}
\STATE \textbf{Input}: sample covariance matrix $S$
\STATE \textbf{Initialization}: $t \leftarrow 0$, $\; \; \, $ $\Sigma^{t} \gets D(S)$,  $\; \; \, $ $\Omega^{t}
\leftarrow \left\{ \Sigma^{t} \right\}^{-1}$.
\WHILE{stopping criterion not fulfilled}
\FOR{$j=1,\ldots,p$} 
\STATE $\bullet \,$ Call routine \textsc{SOLVEBLOCK} below to obtain \begin{align}\label{eq:blockproblem} 
(\wt{\omega}_{jj}, \wt{\omega}_j) \gets \argmin_{\omega_{jj},
  \omega_j} -&\log \det \wt{\Omega}^t(\omega_{jj}, \omega_j ) +
\tr(\wt{\Omega}^t(\omega_{jj}, \omega_j) S)\\
&\text{sb.~to} \; \wt{\omega}_j \leq 0, \; \wt{\Omega}^t(\omega_{jj}, \omega_j) \in \overline{\psd^p} \notag,
\end{align}
where, with the partitioning scheme \eqref{eq:partitioning}, $\wt{\Omega}^t(\omega_{jj}, \omega_j) = {\scriptsize \left( \begin{array}{cc}
                          \omega_{jj} &  \omega_j^{\T}     \\
                          \omega_j  & \left\{ \Omega_{jj}^t \right \}    
                         \end{array} \right)}$.
 
  \STATE $\bullet \,$ $\Omega^{t+1} \gets  \wt{\Omega}^t(\wt{\omega}_{jj},
  \wt{\omega}_j)$ and $\Sigma^{t+1} \gets {\scriptsize \begin{pmatrix}
  s_{jj} & (\sigma_j^{t+1})^{\T} \\  
  \sigma_j^{t+1} & \Sigma_{jj}^t
\end{pmatrix}}
  $,\\ 
\vspace{0.15cm}
 $\qquad \qquad \qquad \qquad \;\;\;\,$ with $\sigma_j^{t+1}  = s_{jj} (\Omega_{jj}^t)^{-1}(-\wt{\omega}_j)$. 
\ENDFOR
\RETURN$(\wh{\Omega}, \wh{\Sigma})$.
\ENDWHILE
\end{algorithmic} 
\vspace{-.3cm}
\rule{\textwidth}{0.01cm}
\textsc{SOLVEBLOCK}
\begin{algorithmic}
\STATE \textbf{Input}: $\Sigma^t$, $S$, $j$. 
\STATE Solve the following linear complementarity problem in $(\eta, \lambda)$ (recall \eqref{eq:partitioning}):
       \begin{equation}\label{eq:lcp}
       \left(\Sigma_{jj} - \frac{\sigma_j \sigma_j^{\T}}{s_{jj}}\right) \eta = \frac{s_j}{s_{jj}} + \lambda, \; \; \; \eta \geq 0, \; \, \lambda \geq 0, \; \; \eta^{\T} \lambda = 0.
       \end{equation}
\RETURN $\wt{\omega}_j \gets -\eta$, $\; \; \, \wt{\omega}_{jj} \gets (1 + s_{jj} \wt{\omega}_j^{\T} \Omega_{jj}^{-1} \wt{\omega}_j)/(1 + s_{jj})$. 
\end{algorithmic}
\end{algorithm}

\subsection{Reduction to a linear complementarity/non-negative least squares problem}     
Algorithm \ref{alg:blockdescent} is a block coordinate descent scheme in which one variable block consisting
of a single column/row is optimized at a time, while the remaining entries are kept fixed. This is cyclically
repeated for all $p$ blocks until a suitable stopping criterion is satisfied. The approach is appealing 
because it turns out that the sub-problems \eqref{eq:blockproblem} are particularly easy to solve by means of a
conversion to linear complementarity problems \eqref{eq:lcp}, for which efficient solvers exist. In the sequel, we show that the routine  
\textsc{SOLVEBLOCK} indeed provides the solution of \eqref{eq:blockproblem}. We start by decomposing
the determinant part as
\begin{equation*}
\det \wt{\Omega}^t(\omega_{jj}, \omega_j) = \det(\Omega_{jj}^t) \cdot \left( \omega_{jj} - \omega_j^{\T} \left\{ \Omega_{jj}^t
\right\} ^{-1} \omega_j \right),
\end{equation*} 
assuming for a moment that $\Omega_{jj}^t \in \psd^{p-1}$ as will be shown below. After taking logarithms, the first
factor becomes a constant not depending on the optimization variables $(\omega_{jj},  \omega_j)$ and can hence be omitted. Similarly, 
the trace term in \eqref{eq:blockproblem} can be decomposed as
\begin{equation*}
\tr(\wt{\Omega}^t(\omega_{jj}, \omega_j) S) = 2 \omega_j^{\T} s_j + s_{jj} \omega_{jj} + \tr(S_{jj} \Omega_{jj}^t),
\end{equation*}
where $s_{jj}, s_j, S_{jj}$ are the components of a partitioning of $S$ analogous to that of $\Omega^t$. The
term  $\tr(S_{jj} \Omega_{jj}^t)$ does not depend on the optimization variables and can be dropped as well. Altogether, we find that
\eqref{eq:blockproblem} is equivalent to the following optimization problem:
\begin{align}
\begin{split}\label{eq:innerproblemequivalent}
  &\min_{\omega_{jj}, \omega_j}  -\log(\omega_{jj} - \omega_j^{\T} \left\{
  \Omega_{jj}^t \right\} ^{-1} \omega_j) + 2 \omega_j^{\T} s_j + s_{jj}
\omega_{jj},\\
&\text{sb.~to} \; \; \omega_j \leq 0, \quad \omega_{jj} - \omega_j^{\T} \left\{
  \Omega_{jj}^t \right\} ^{-1} \omega_j \geq 0, \quad \omega_{jj} \geq 0, 
\end{split}
\end{align}
Observe that the constraint $\wt{\Omega}^t(\omega_{jj}, \omega_j) \in \overline{\psd^p}$ is equivalent to the second
and third constraint in \eqref{eq:innerproblemequivalent} given $\Omega_{jj}^t \in \psd^{p-1}$, recalling
that a symmetric matrix is positive semi-definite if and only if all its principal minors are non-negative.
Setting $\log(x) = -\infty$ if $x \leq 0$, the second and third constraint can be dropped as long as 
$\Omega_{jj}^t \in \psd^{p-1}$. 
The objective function in \eqref{eq:innerproblemequivalent} consists of a linear part and the
composition of the negative logarithm (which is convex and non-increasing) and
a concave function (assuming again $\Omega_{jj}^t \in \psd^{p-1}$). Such a function is again
convex (\cite{BoydVandenberghe2004}, p.84). We conclude that $(\wt{\omega}_{jj}, \wt{\omega}_j)$ is
a minimizer of \eqref{eq:blockproblem} if and only if it satisfies the KKT optimality conditions of 
\eqref{eq:innerproblemequivalent} given by
\begin{align}\label{eq:kkt_innerproblem}
&\frac{1}{\wt{\omega}_{jj} - \wt{\omega}_{j}^{\T} \left\{ \Omega_{jj}^t \right
  \}^{-1}\wt{\omega}_j} = s_{jj}, \quad \frac{ \left\{ \Omega_{jj}^t \right
  \}^{-1} \omega_j}{\wt{\omega}_{jj}  - \wt{\omega}_{j}^{\T} \left\{ \Omega_{jj}^t \right
  \}^{-1}\wt{\omega}_j} = -(s_{j} + \kappa),\\
& \wt{\omega}_j \leq 0, \quad \kappa \geq 0, \quad \wt{\omega}_j^{\T} \kappa = 0. \notag
\end{align}     
where $\kappa \geq 0$ is a Lagrangian multiplier. Substituting the first equation in \eqref{eq:kkt_innerproblem} into
the second one, we obtain 
\begin{equation*}
\left\{ \Omega_{jj}^t \right
  \}^{-1}  \wt{\omega}_j = \frac{-(s_j + \kappa)}{s_{jj}}.  
\end{equation*} 
It is hence possible to solve the KKT system \eqref{eq:kkt_innerproblem} for $\wt{\omega}_j$ first and then resolve
for $\wt{\omega}_{jj}$. This automatically ensures that  $\wt{\omega}_{jj} - \wt{\omega}_{j}^{\T} \left\{ \Omega_{jj}^t \right
  \}^{-1}\wt{\omega}_j > 0$ and in turn that $\Omega^{t + 1} \in \psd^p$. Applying this argument recursively, it follows that all 
iterates must be strictly positive definite, provided the initial iterate $\Omega^0 = \{ D(S) \}^{-1}$ is. Re-parameterizing
$\eta = -\wt{\omega}_j$ and $\lambda = \kappa/s_{jj}$, solution of \eqref{eq:innerproblemequivalent} respectively
\eqref{eq:kkt_innerproblem} boils down to solving 
\begin{equation}\label{eq:lcp_b}
\left\{ \Omega_{jj}^t \right
  \}^{-1} \eta = \frac{s_{j}}{s_{jj}}, \quad \eta^{\T} \lambda = 0,  \quad \eta \geq 0, \quad \lambda
\geq 0,
\end{equation}
We finally recover \eqref{eq:lcp} in Algorithm \ref{alg:blockdescent} from \eqref{eq:lcp_b} by re-writing $\left\{ \Omega_{jj}^t \right \}^{-1}$ using Schur complements. 
Problem \eqref{eq:lcp_b} is a so-called \emph{monotone linear complementarity problem} (\cite{BermanPlemmons}, Ch.11). It is not hard 
to handle \eqref{eq:lcp_b}, because it is equivalent to the quadratic programming problem 
\begin{equation*}
\min_{\eta \geq 0} \frac{1}{2} \eta^{\T} \left\{ \Omega_{jj}^t \right \}^{-1} \eta - \eta^{\T} s_j/s_{jj},
\end{equation*}
which  is in turn equivalent to the non-negative least squares (NNLS) problem   
\begin{equation*}
\min_{\eta \geq 0} \frac{1}{2}\norm{b - A \eta}_2^2, \quad A = \left\{ \Omega_{jj}^t \right
  \}^{-1/2}, \quad b = \left\{ \Omega_{jj}^t \right
  \}^{1/2} \frac{s_j}{s_{jj}}.
\end{equation*}
Consequently, \eqref{eq:lcp_b} can be solved by one of the many existing solvers for the previous two problems. We use the block principal pivoting algorithm of
\cite{Portugal1994}, which operates directly on the linear complementarity
problem and is experimentally the fastest method for strictly convex NNLS problems where the number 
of variables does not exceed a few thousands \cite{Slawski2013nnlsalg}.

\subsection{Properties}

\paragraph{Convergence} We now state that Algorithm \ref{alg:blockdescent} converges to the unique minimizer $\wh{\Omega}$ 
of \eqref{eq:logdetext}, existence provided (cf.~Theorem \ref{theo:existenceanduniqueness}). 
\begin{theo}\label{theo:convergence} Under the conditions of Theorem \ref{theo:existenceanduniqueness}, the sequence of iterates of Algorithm \ref{alg:blockdescent} satisfies $\lim_{t \rightarrow \infty} \Omega^t = \wh{\Omega}$.  
\end{theo} 
The statement can be derived as a consequence of general result regarding block coordinate descent in \cite{Bertsekas1999}. Regarding
the speed of convergence, our experiments suggest a linear rate (see Figure \ref{fig:runtimes}). This is observation is supported by theoretical work in \cite{LuoTseng1992}. 
\paragraph{Computational complexity} One complete cycle of block updates requires $O(p^4)$ operations, with each call to \textsc{SOLVEBLOCK} amounting to
$O(p^3)$ operations. The workhorse in \textsc{SOLVEBLOCK} is block principal pivoting, in which a linear system of dimension $p$ has to be solved
per iteration. \textsc{SOLVEBLOCK} typically terminates after few iterations. 
\paragraph{Stopping criterion} We suggest to stop Algorithm \ref{alg:blockdescent} once one comes close to KKT optimality. In view of \eqref{eq:kkt}, this can
be quantified by means of the criterion 
\begin{align}\label{eq:kkt_practice}
\eps = \max\{\eps_1, \eps_2 \}, \; \text{where}& \; \eps_1 = \max_{(j,k) \in \mc{E}^t} |\sigma_{jk}^t - s_{jk}|, \\
   &\eps_2 = \max_{(j,k) \notin \mc{E}^t} \max\{s_{jk} - \sigma_{jk}^t,0\},\, \;  \mc{E}^t = \{(j,k):\, \omega_{jk}^t < 0 \}. \notag
\end{align}        
\section{Experiments}\label{sec:experiments}

In the first part, we use synthetic datasets to study systematically the performance of
our thresholding approach (Section \ref{sec:sparsification}) in the high-dimensional, sparse regime  as compared to various competing methods proposed 
in the literature. In the second part, we present possible applications of precision matrix
estimation under non-positivity constraints on the off-diagonal elements. Specifically, we consider
learning of taxonomies and analysis of planar landmark data.
\begin{figure}
\begin{tabular}{cccc}
\hspace{-.4cm}\textsf{chain} &  \textsf{grid} & \textsf{random} & \textsf{star} \\ 
\hspace{-0.6cm}\includegraphics[height = 0.19\textheight]{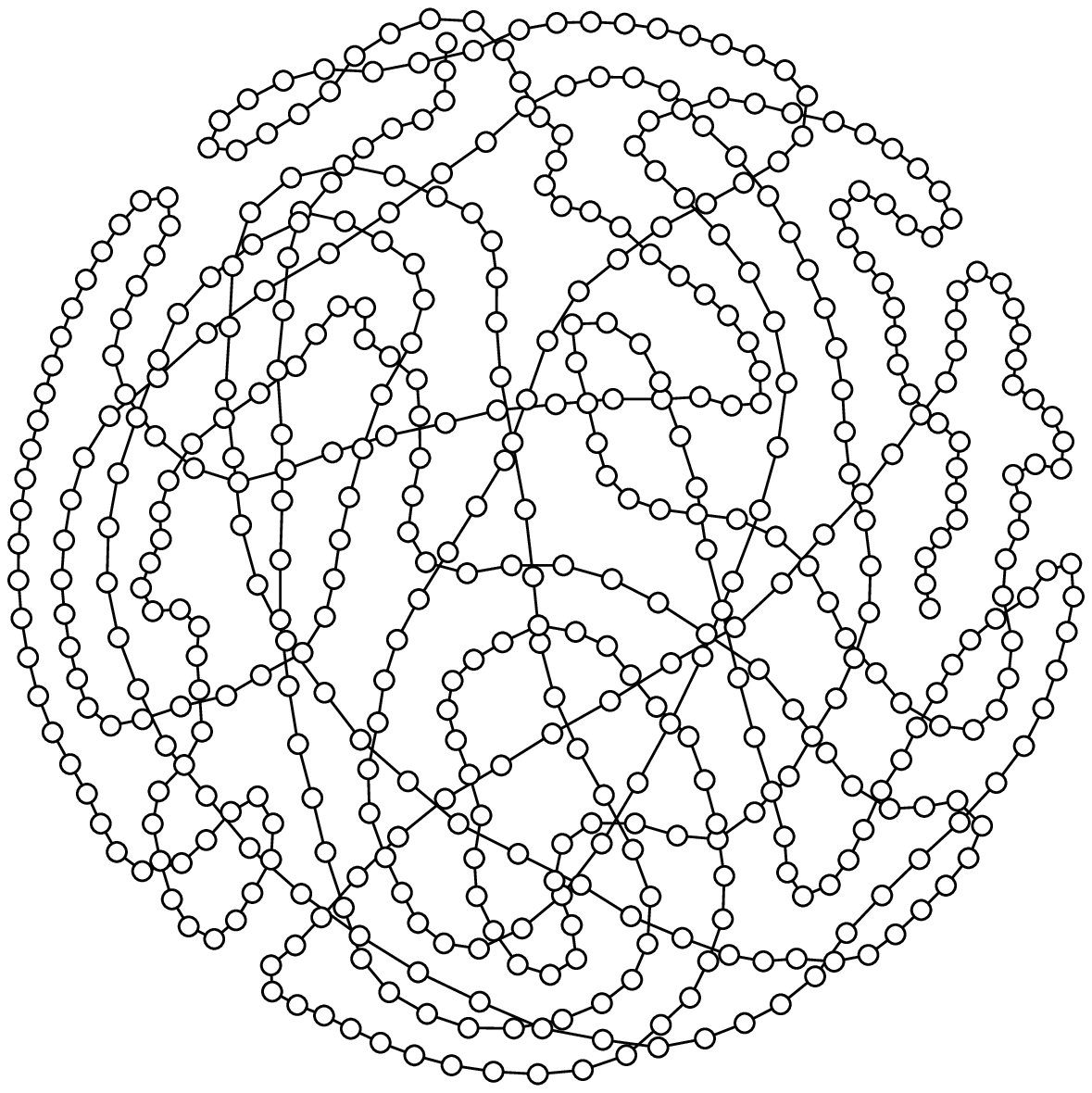} & 
\hspace{-.4cm}\includegraphics[height = 0.19\textheight]{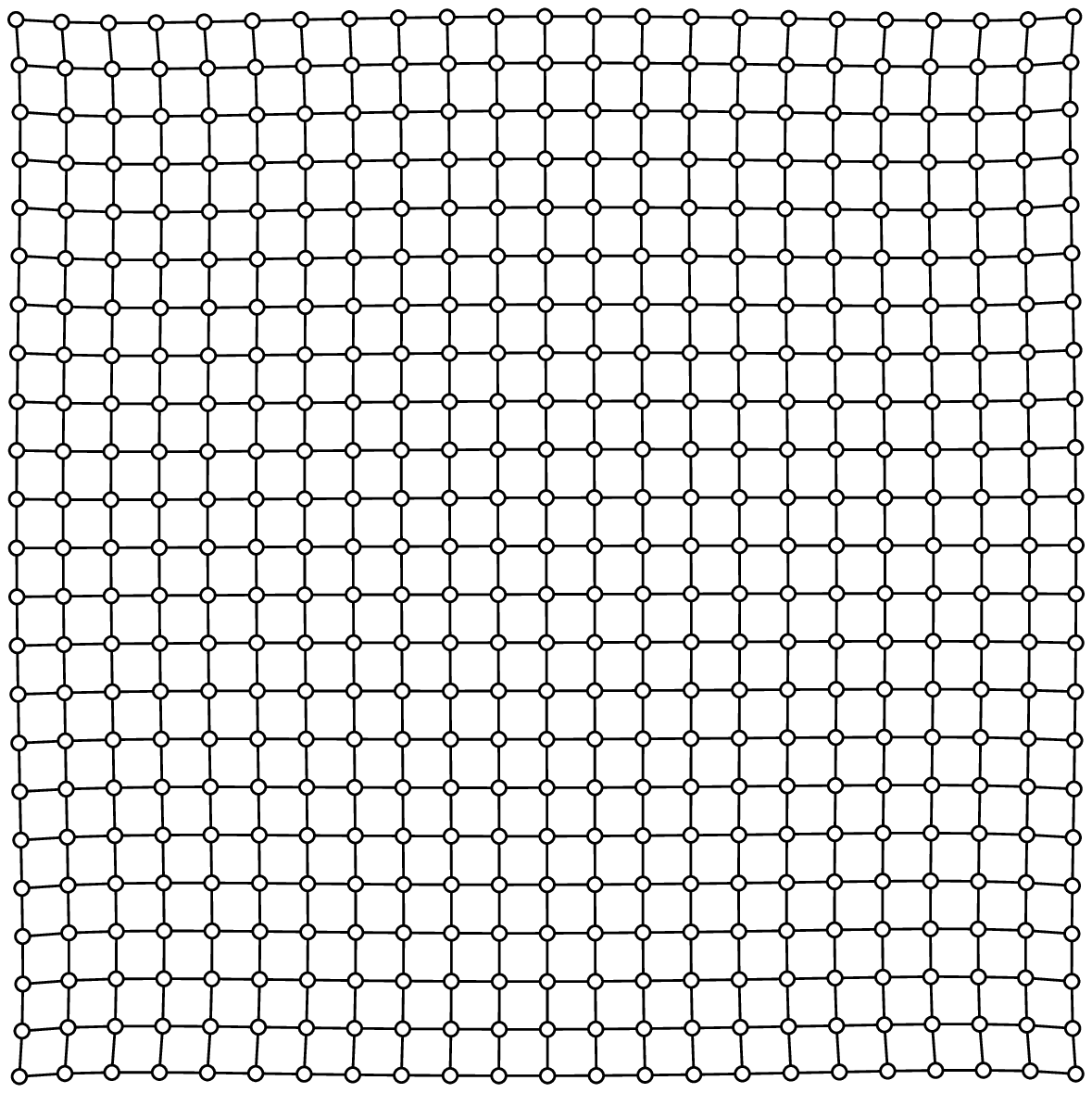} & 
\hspace{-.4cm}\includegraphics[height = 0.19\textheight]{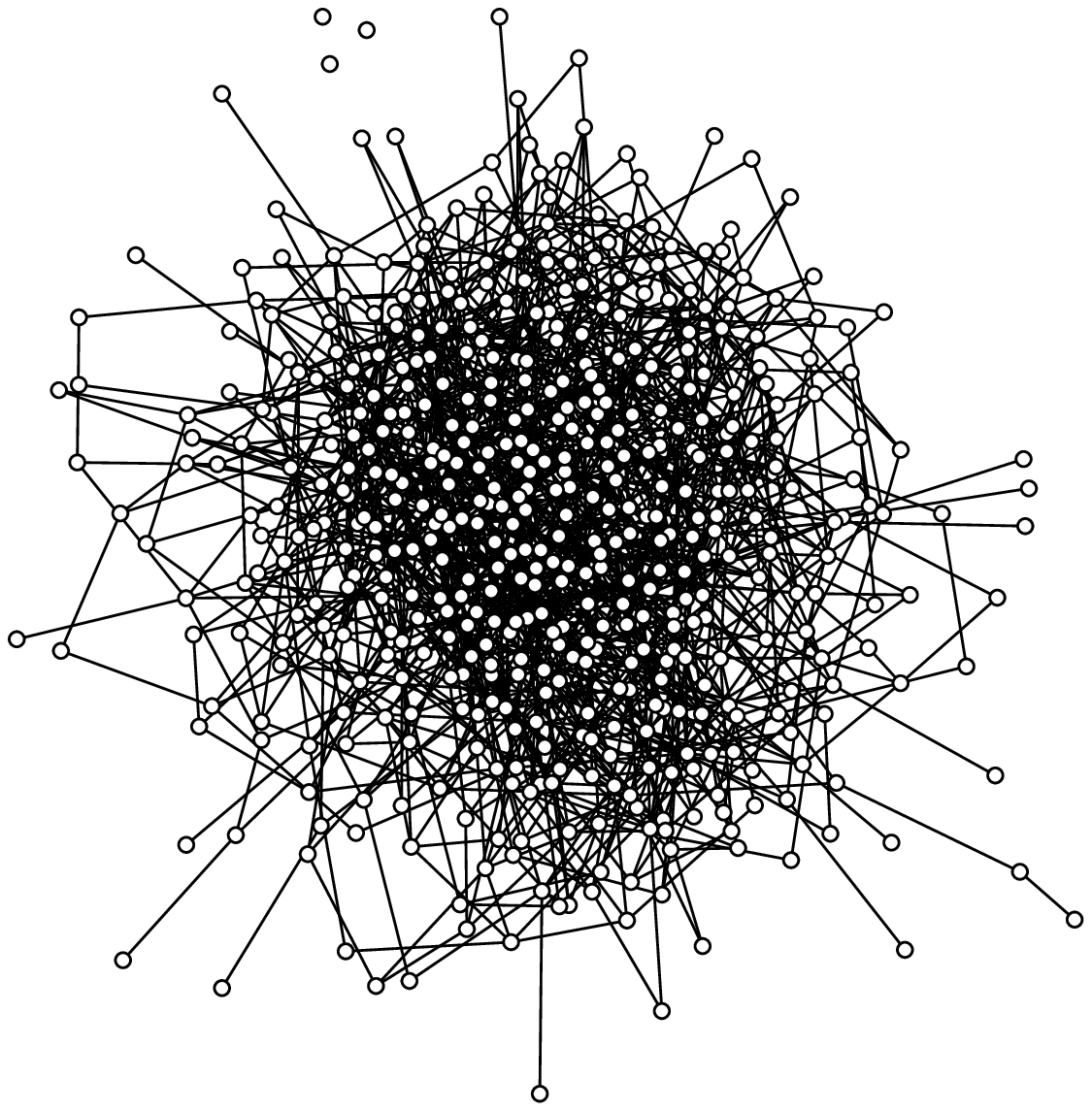} &
\hspace{-.4cm}\includegraphics[height = 0.15\textheight]{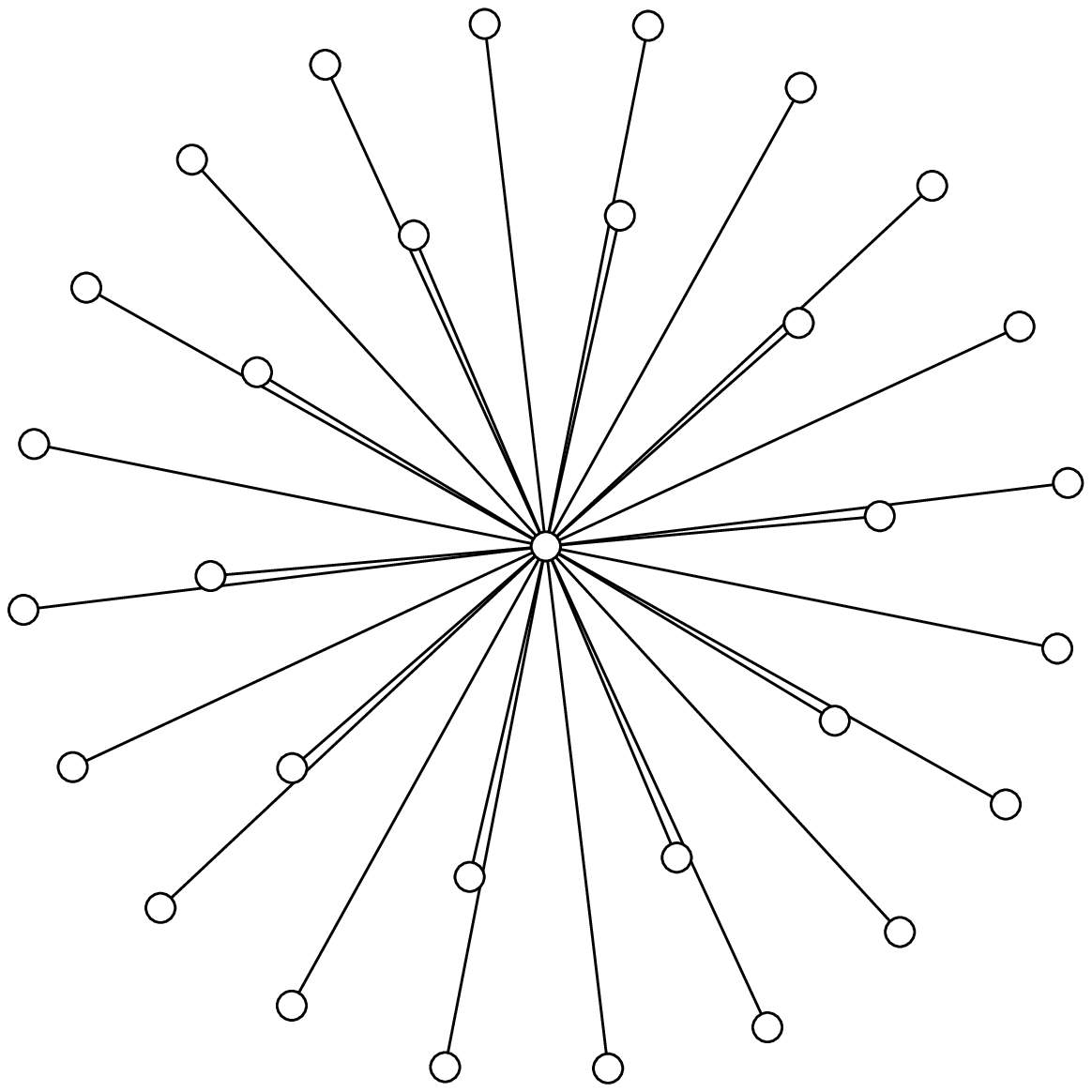} \\ 
\end{tabular}
\caption{Graphs underlying a subset of the experimental setups.}\label{fig:graphstructures}
\end{figure}
\subsection{Synthetic data}
\begin{figure}
\vspace{-0.125\textheight}
\begin{tabular}{ccc}


\hspace{-.7cm}\includegraphics[height = 0.20\textheight]{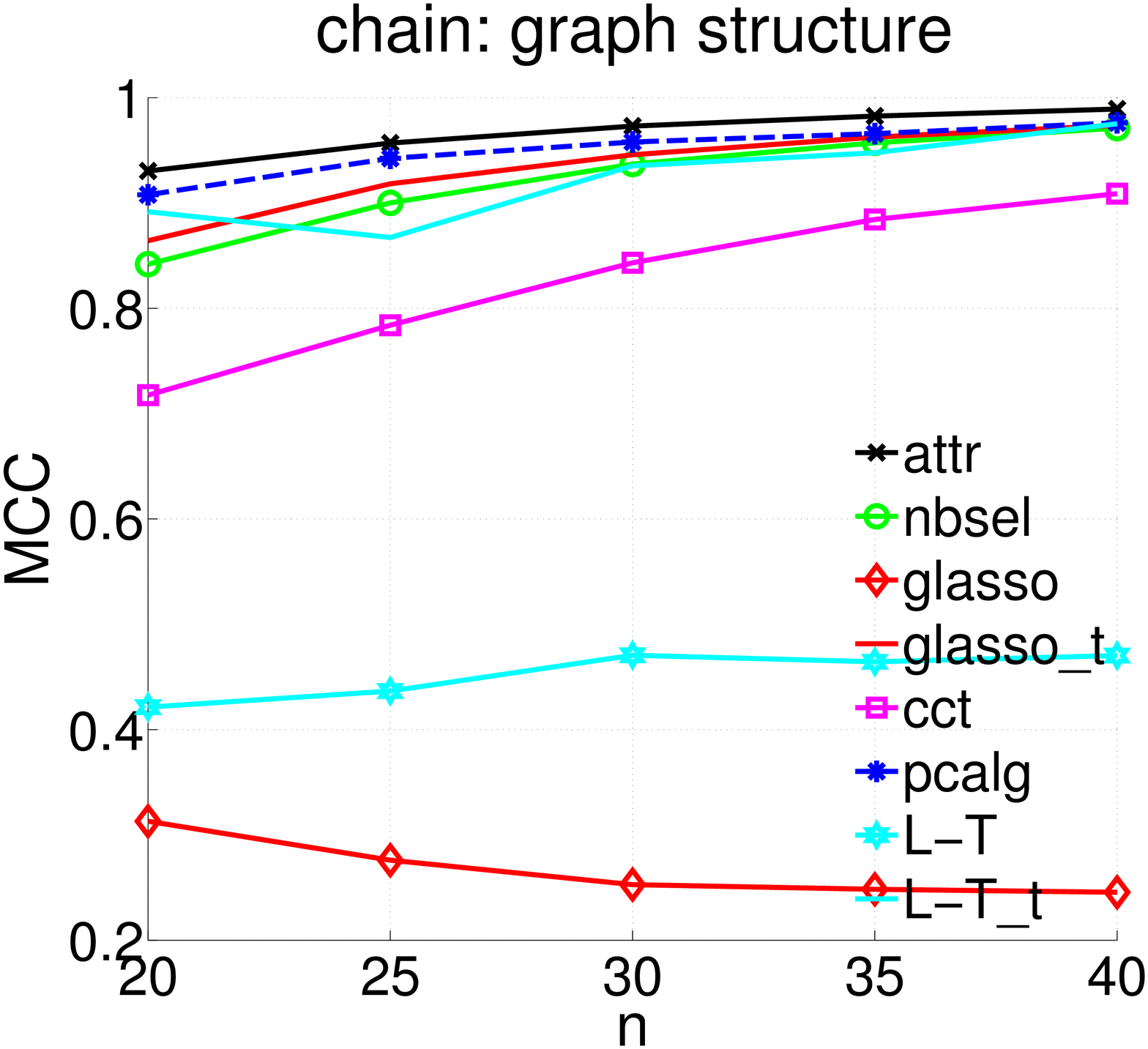} & 
\hspace{-.1cm}\includegraphics[height = 0.20\textheight]{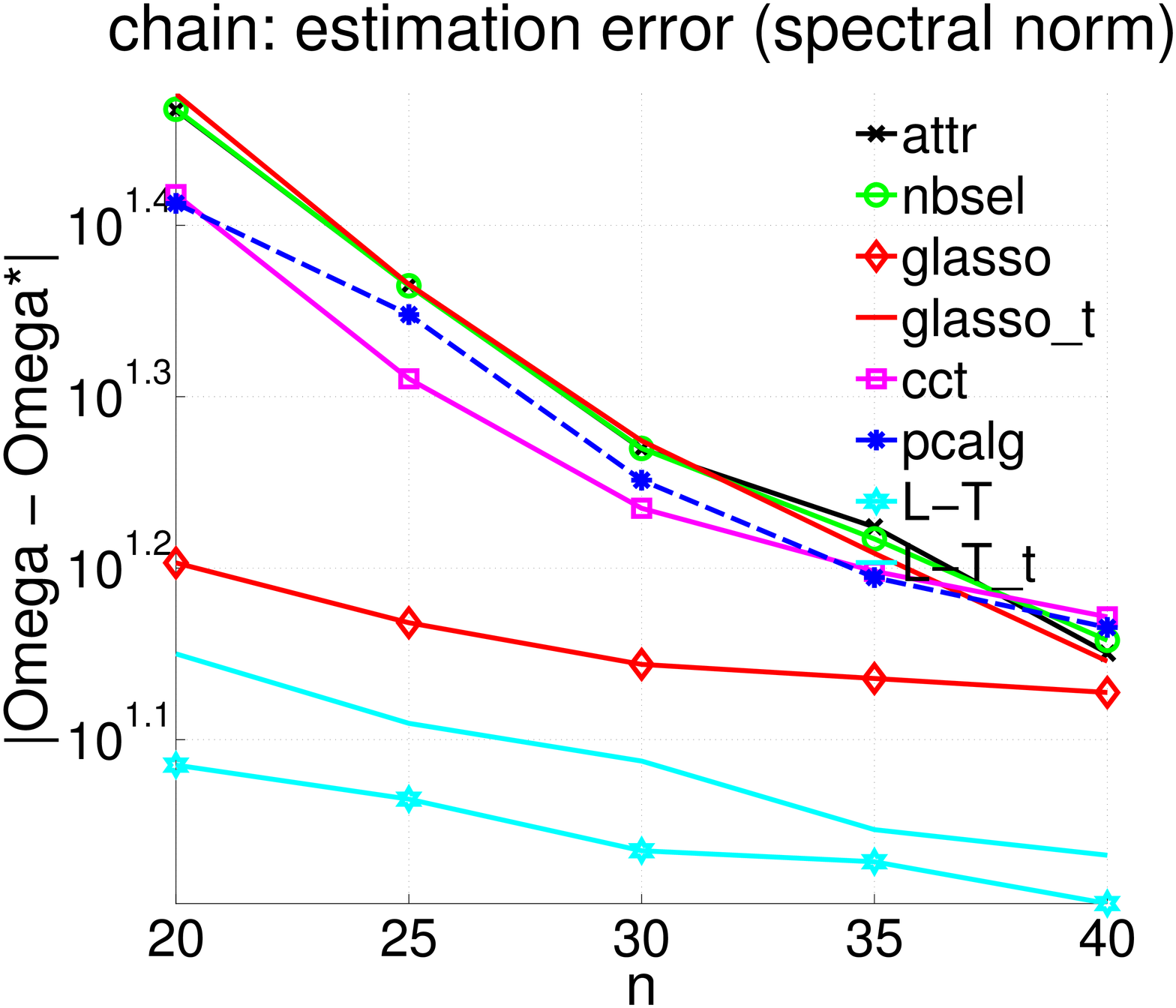} &
\hspace{-.1cm}\includegraphics[height = 0.20\textheight]{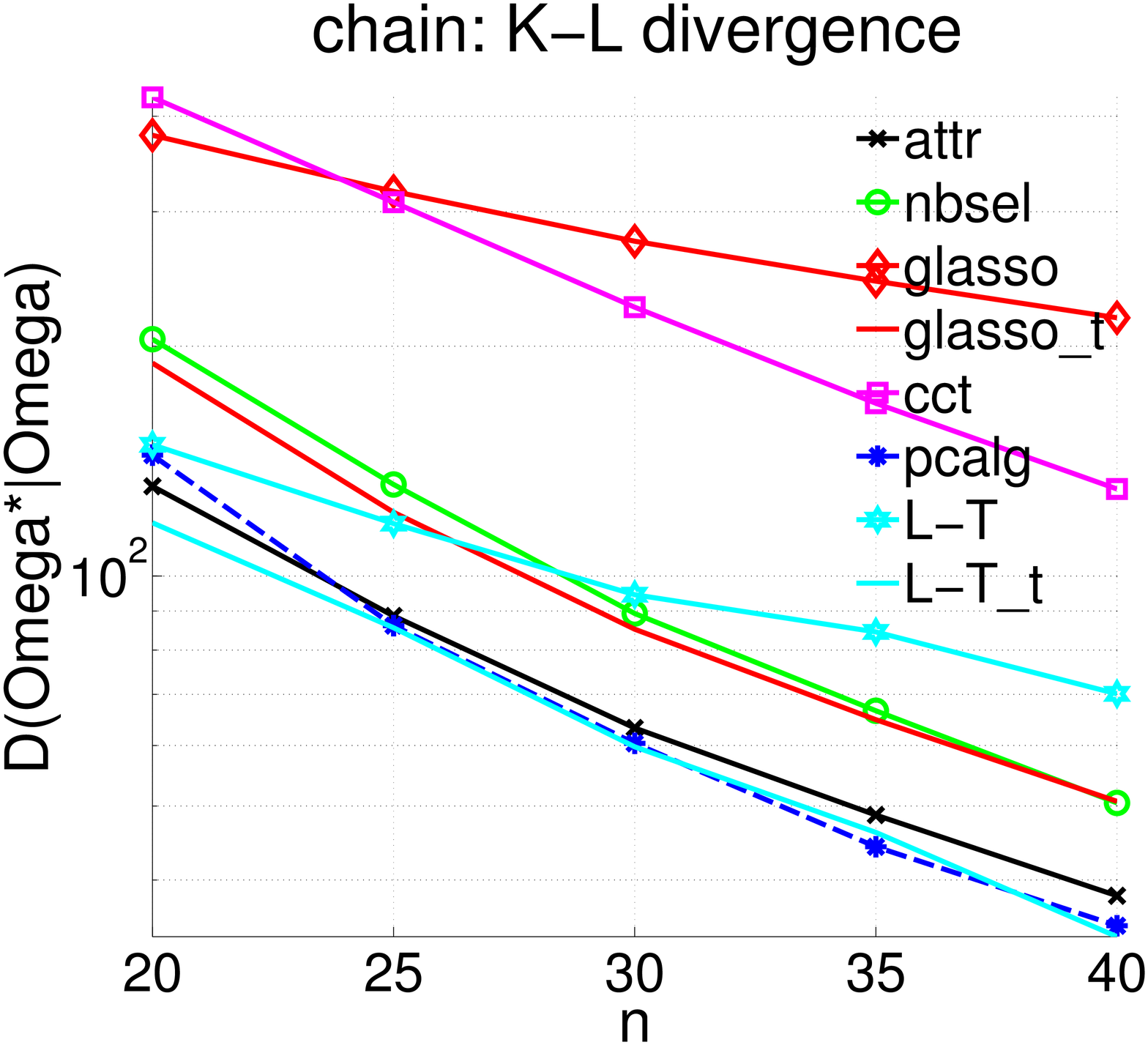} \\ 
\hspace{-.7cm}\includegraphics[height = 0.20\textheight]{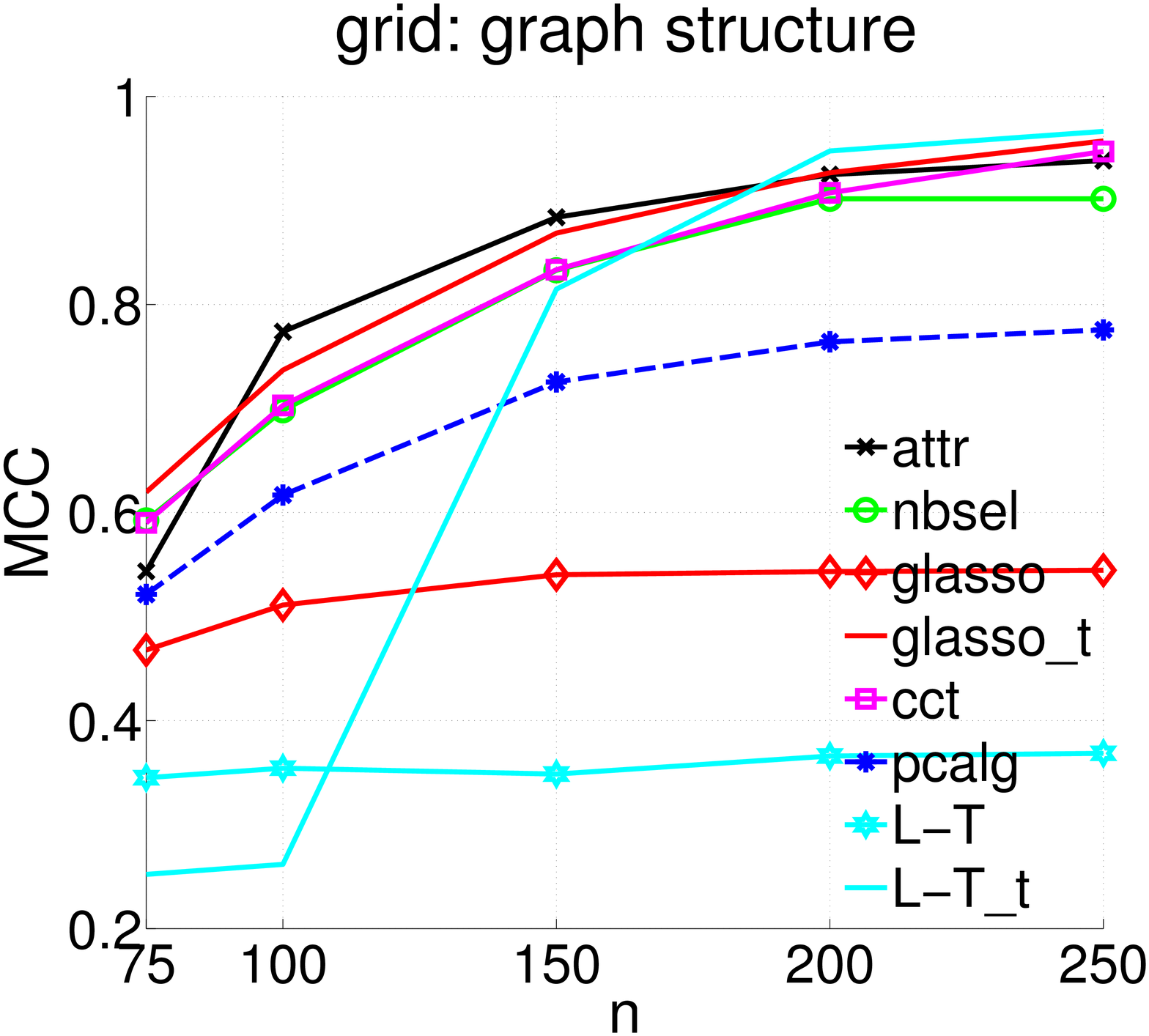} & 
\hspace{-.1cm}\includegraphics[height = 0.20\textheight]{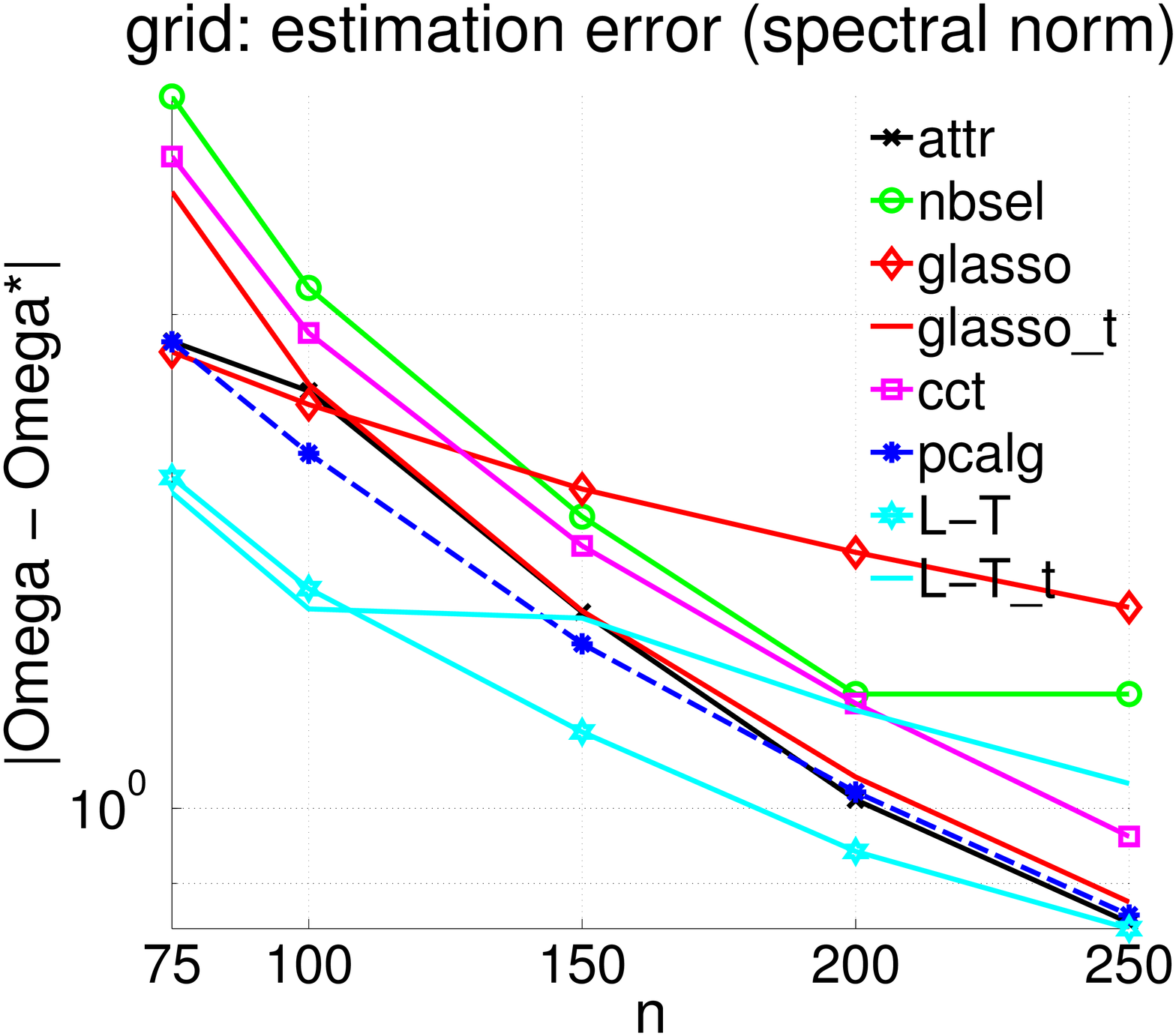} &
\hspace{-.1cm}\includegraphics[height = 0.20\textheight]{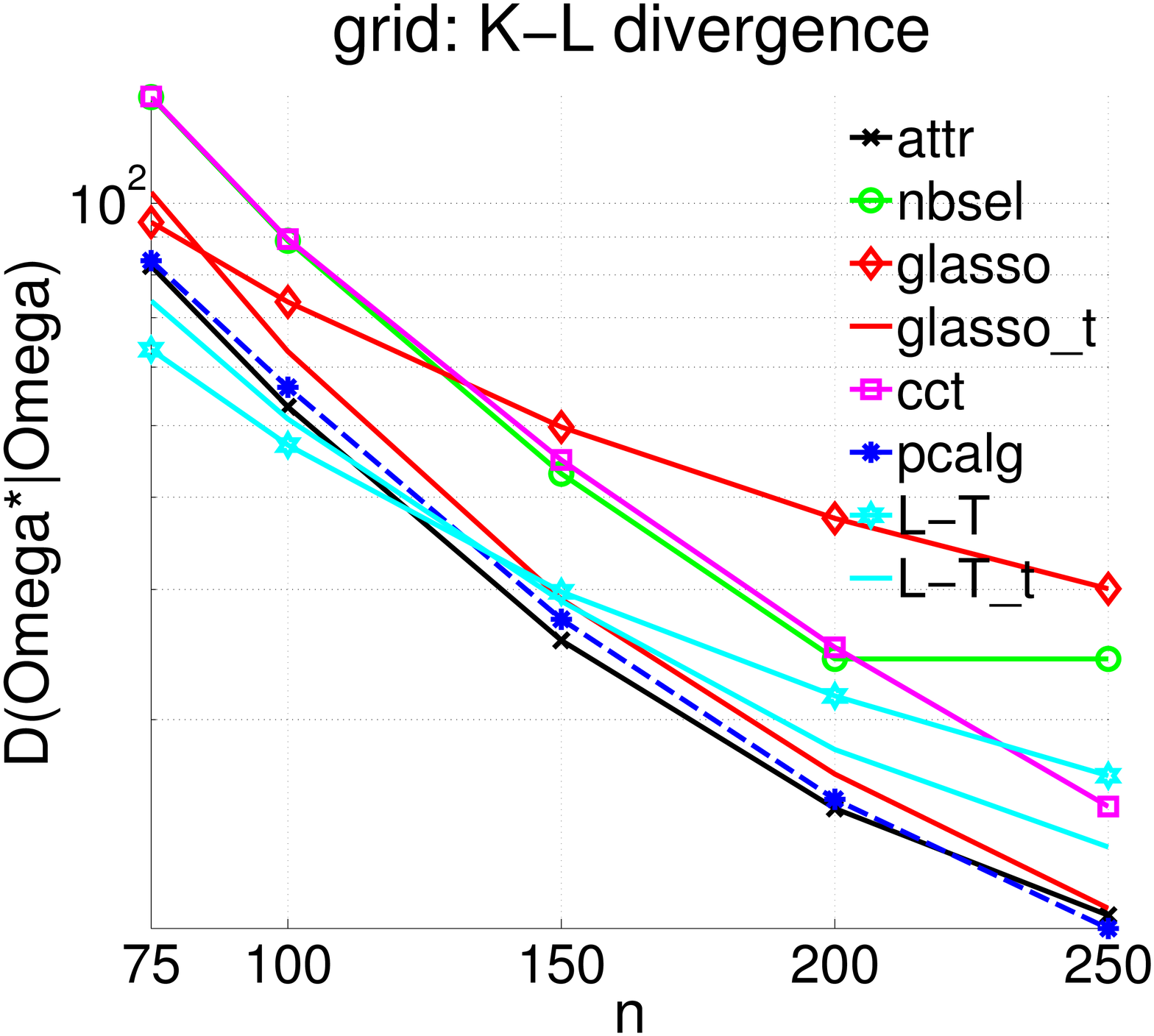} \\ 
\hspace{-.7cm}\includegraphics[height = 0.20\textheight]{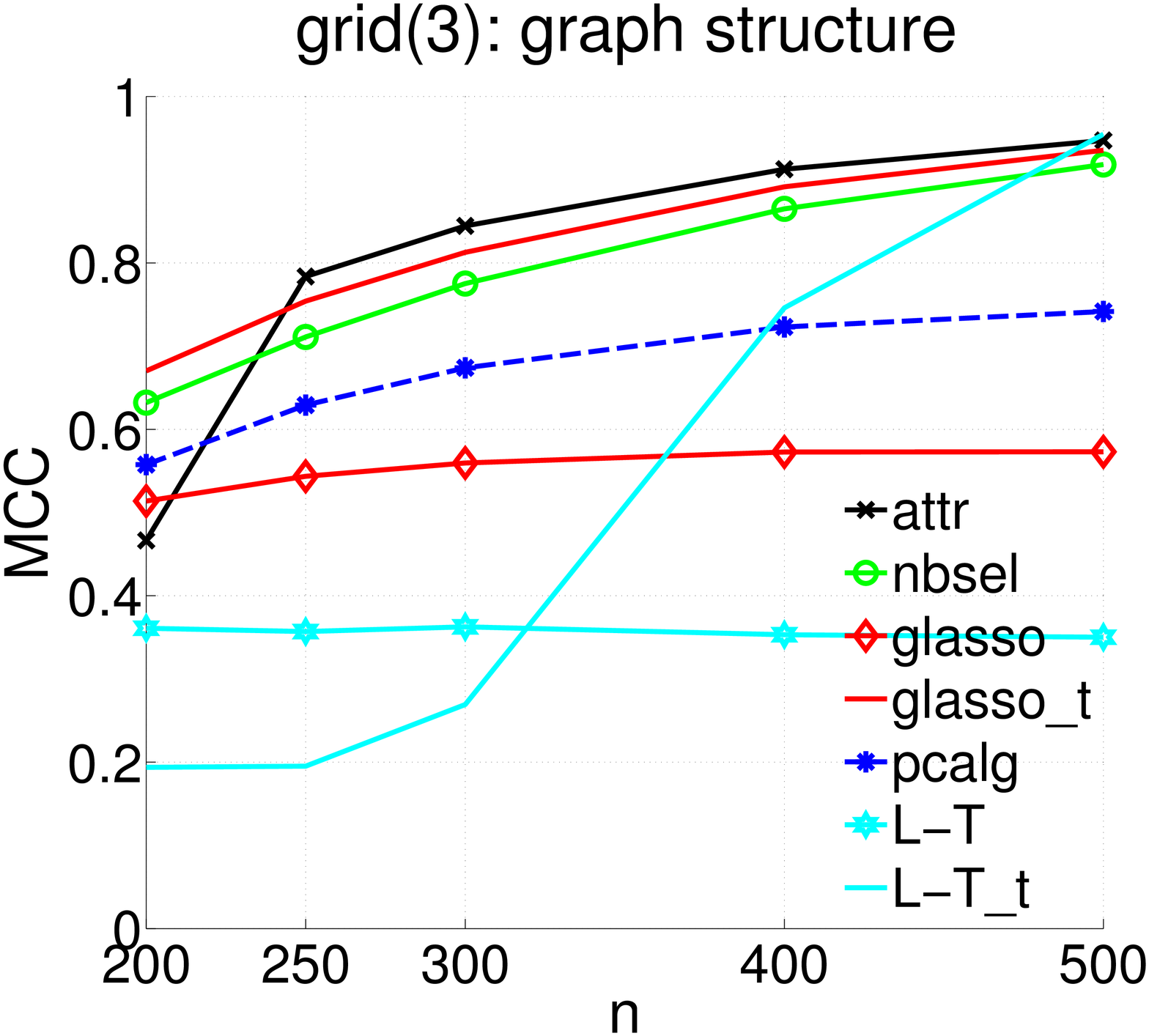} & 
\hspace{-.1cm}\includegraphics[height = 0.20\textheight]{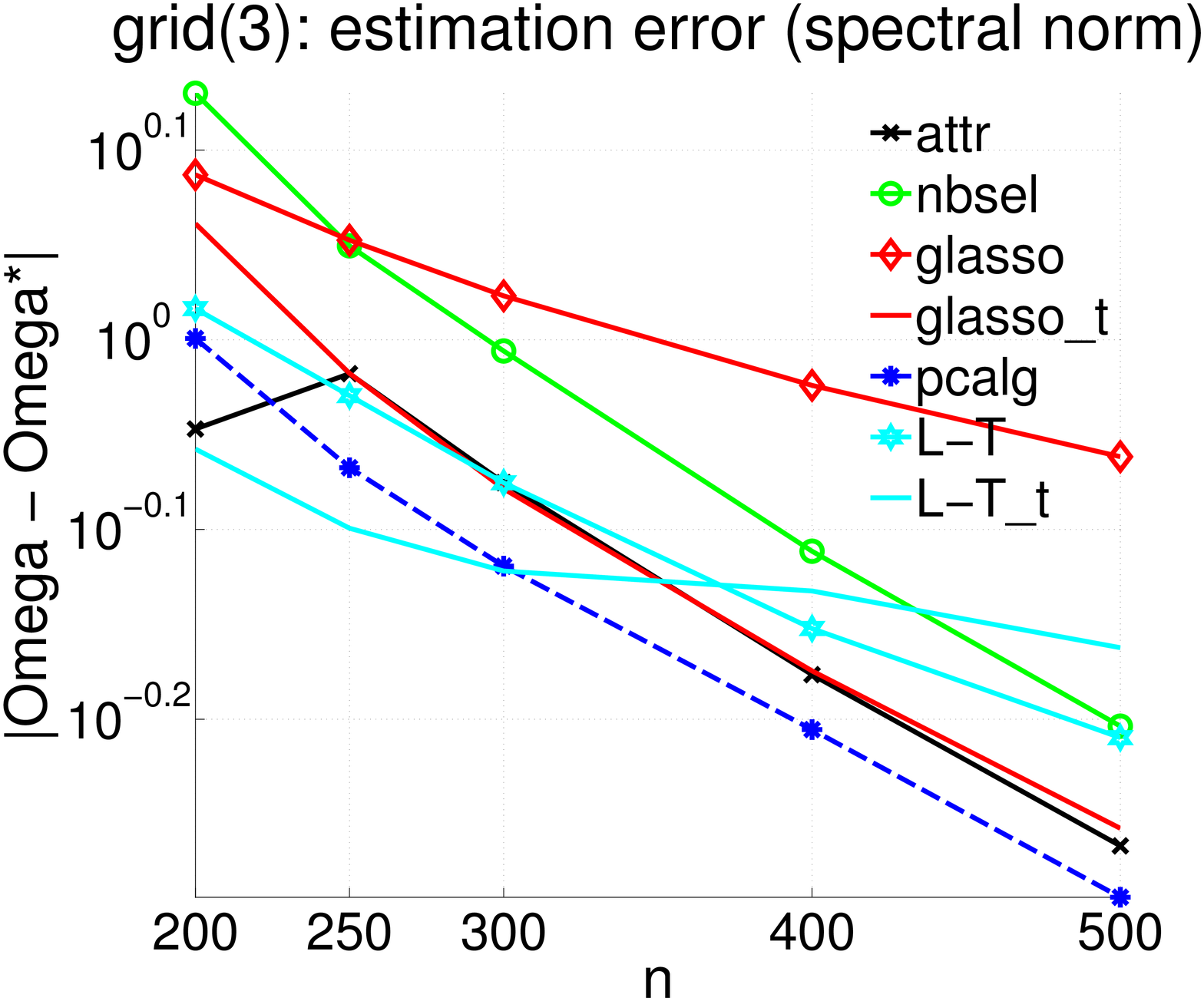} &
\hspace{-.1cm}\includegraphics[height = 0.20\textheight]{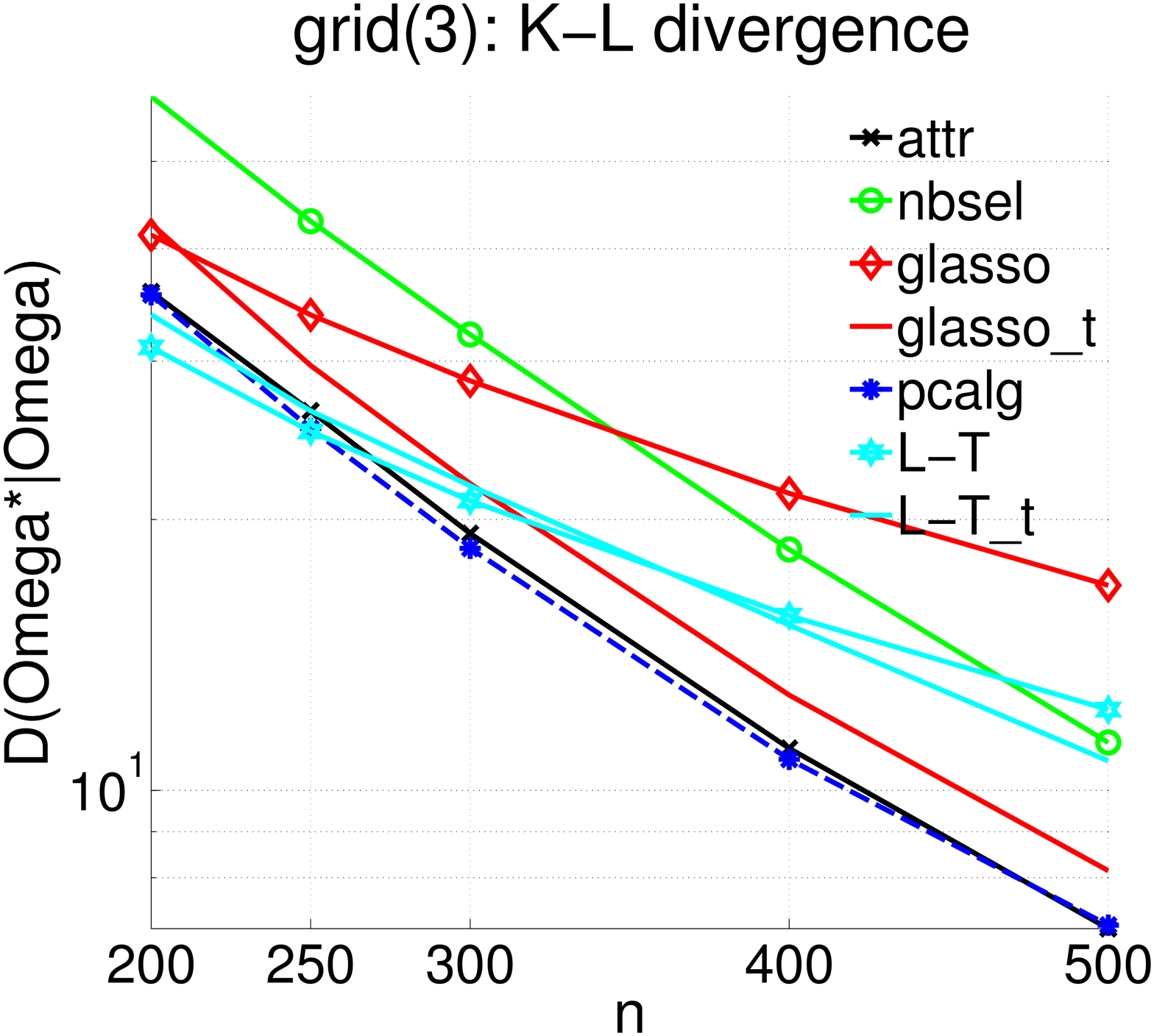} \\ 
\hspace{-.7cm}\includegraphics[height = 0.20\textheight]{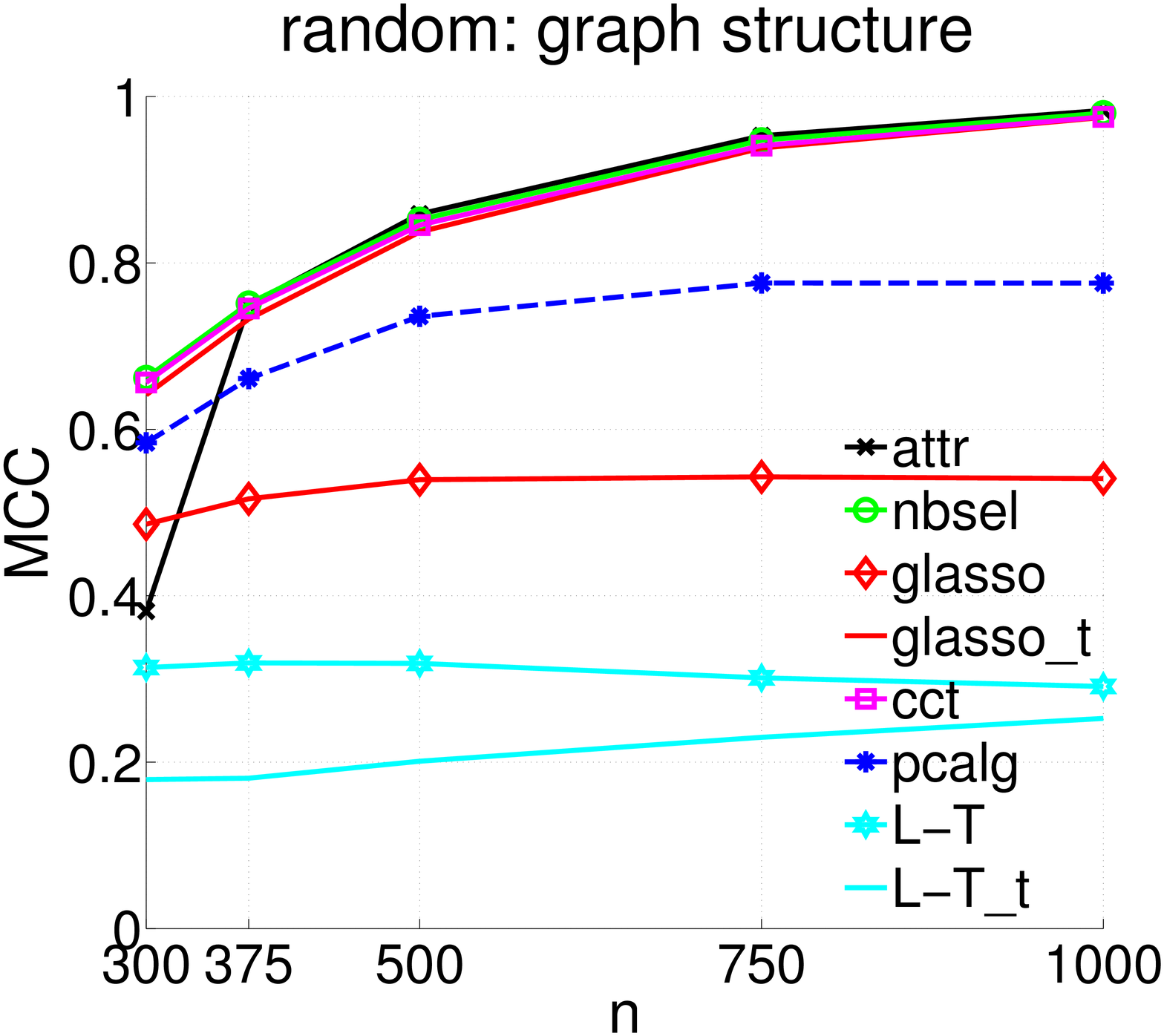} & 
\hspace{-.1cm}\includegraphics[height = 0.20\textheight]{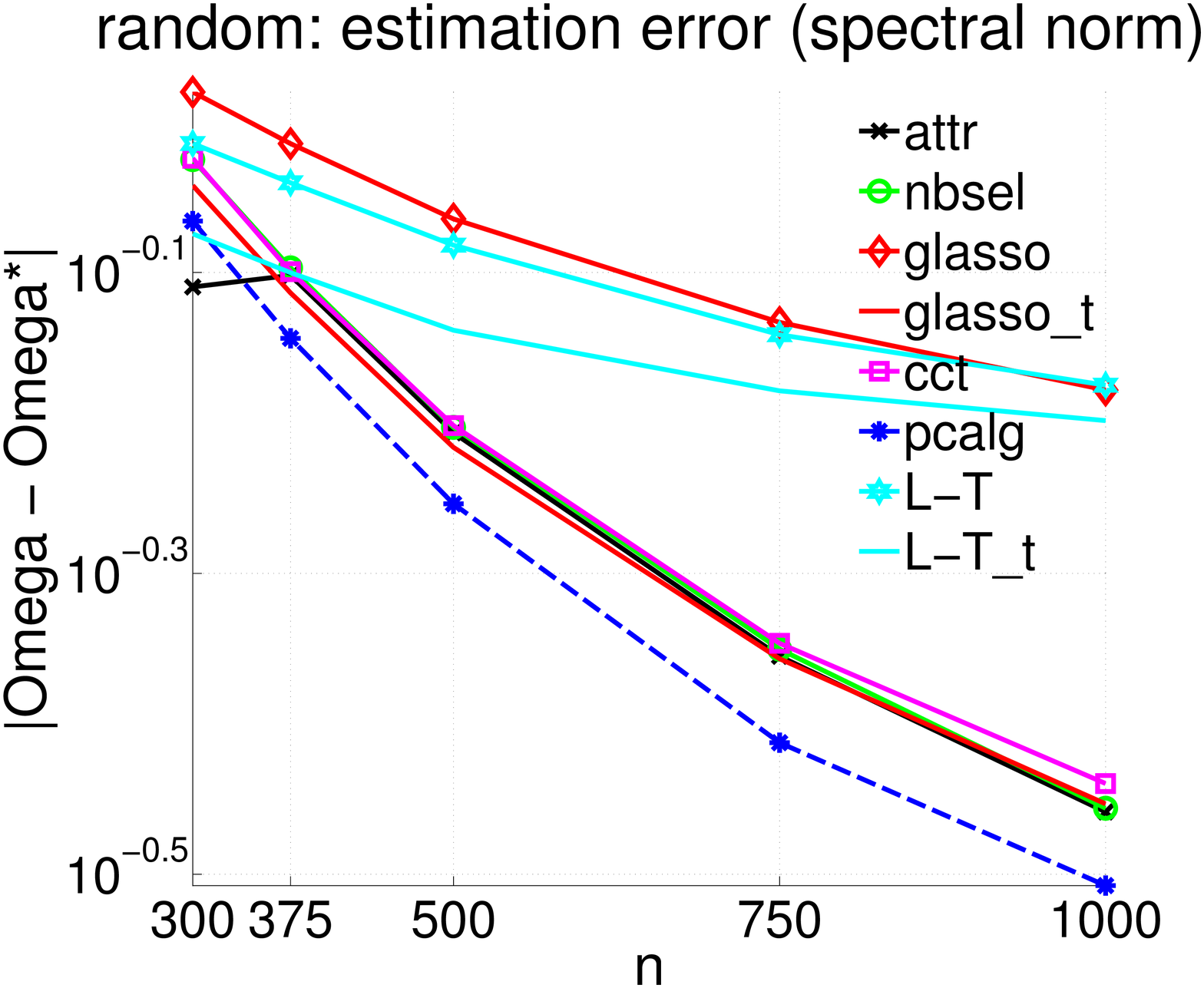} &
\hspace{-.1cm}\includegraphics[height = 0.20\textheight]{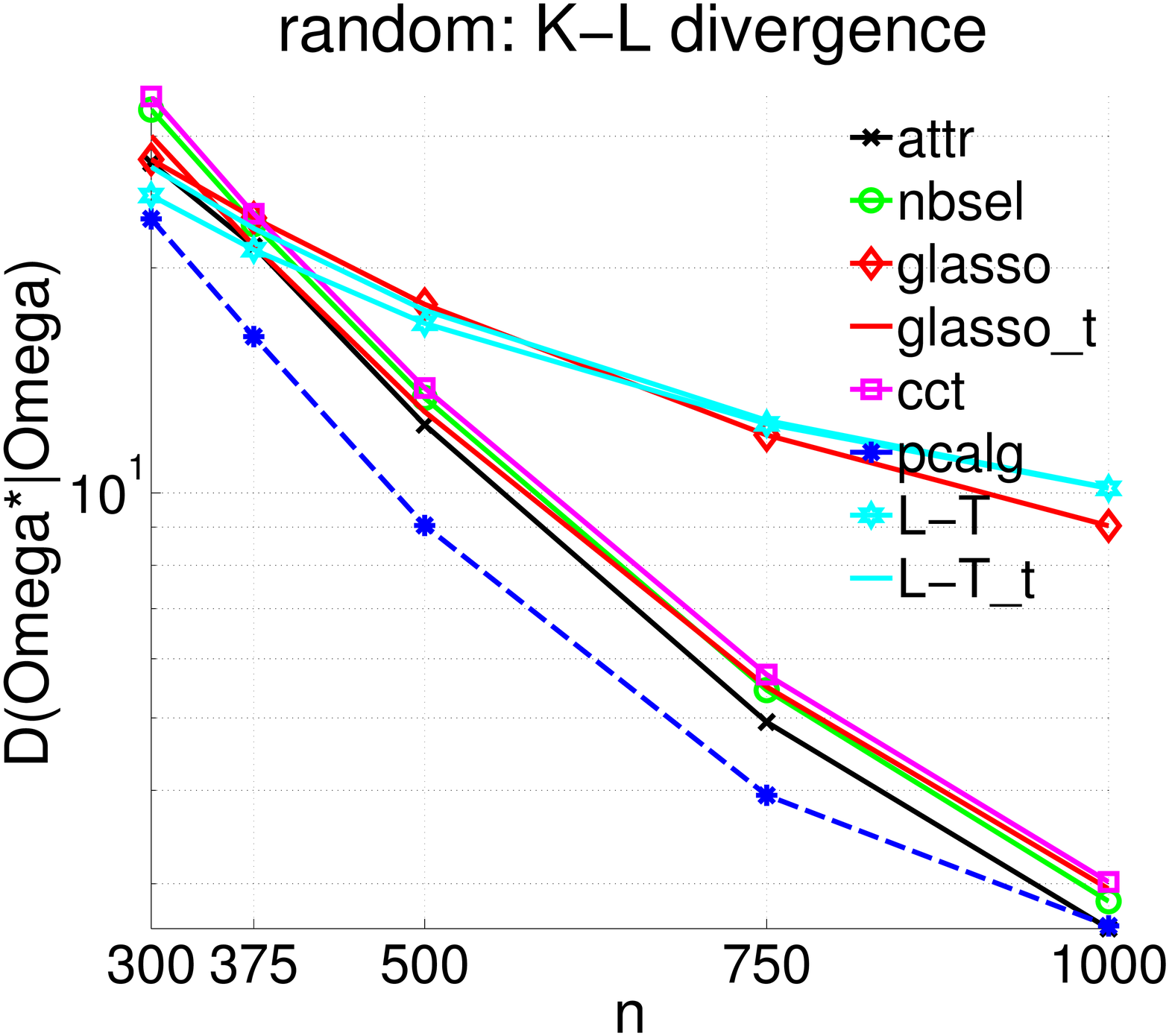} \\ 
\hspace{-.7cm}\includegraphics[height = 0.20\textheight]{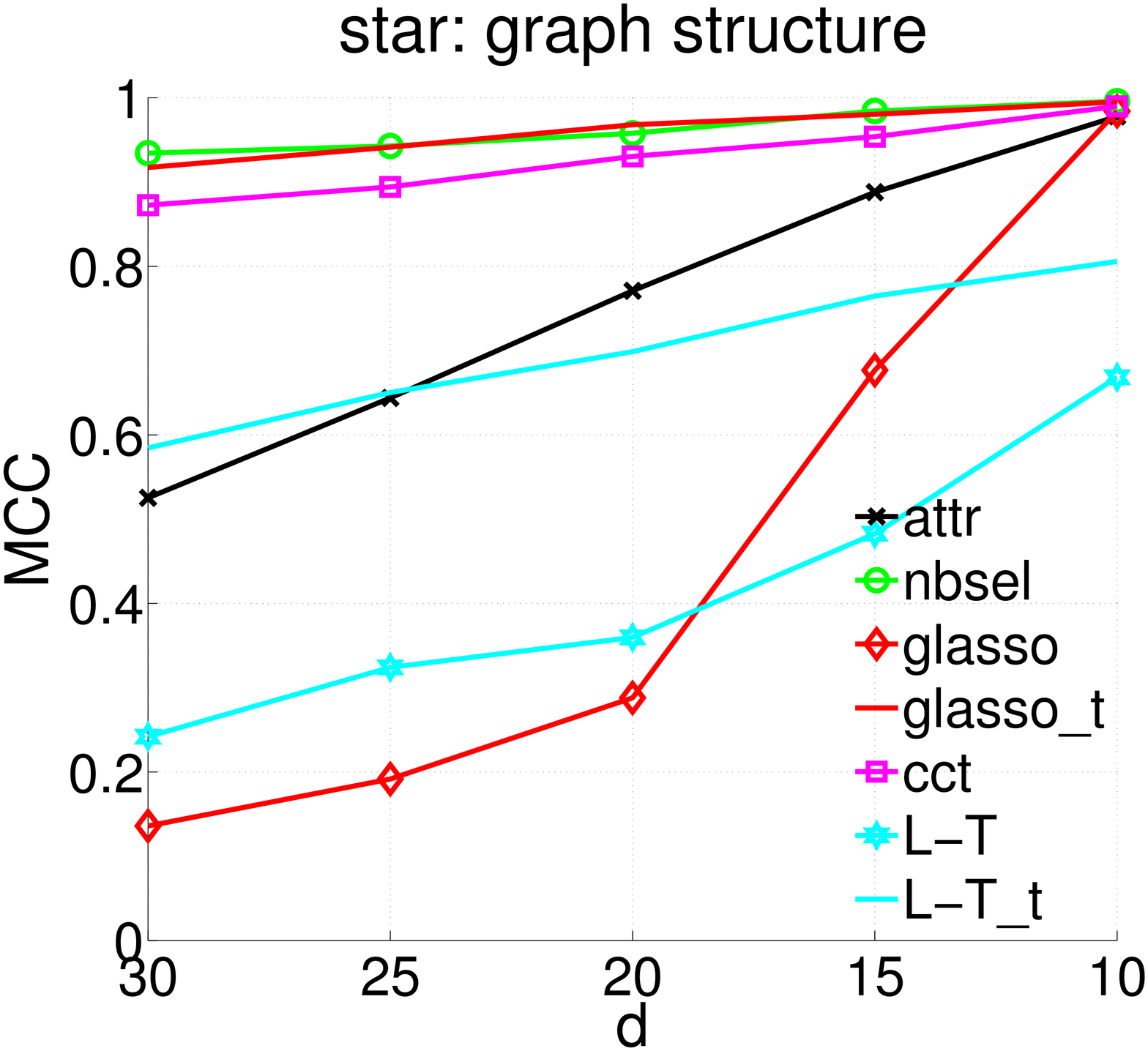} & 
\hspace{-.1cm}\includegraphics[height = 0.20\textheight]{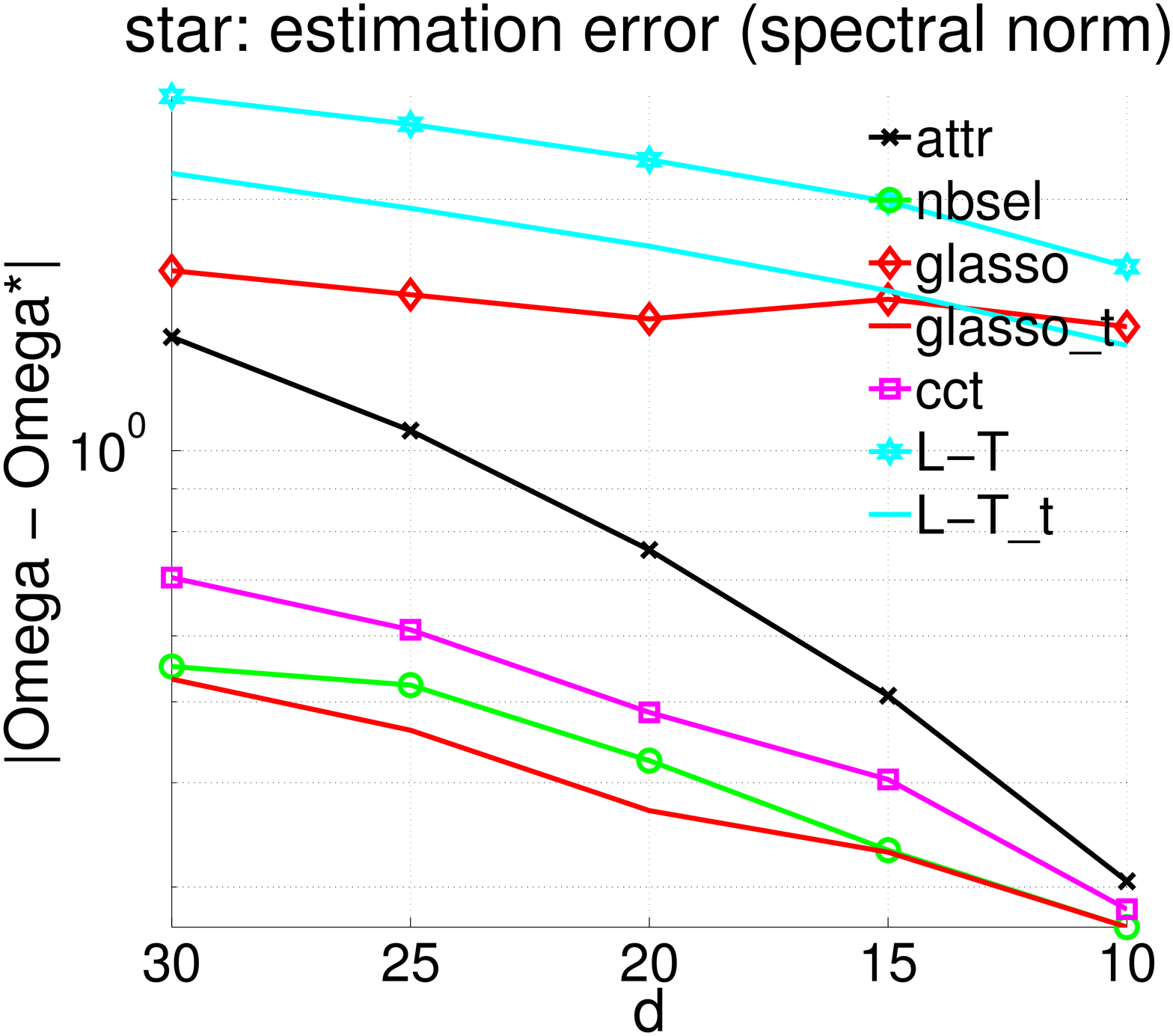} &
\hspace{-.1cm}\includegraphics[height = 0.20\textheight]{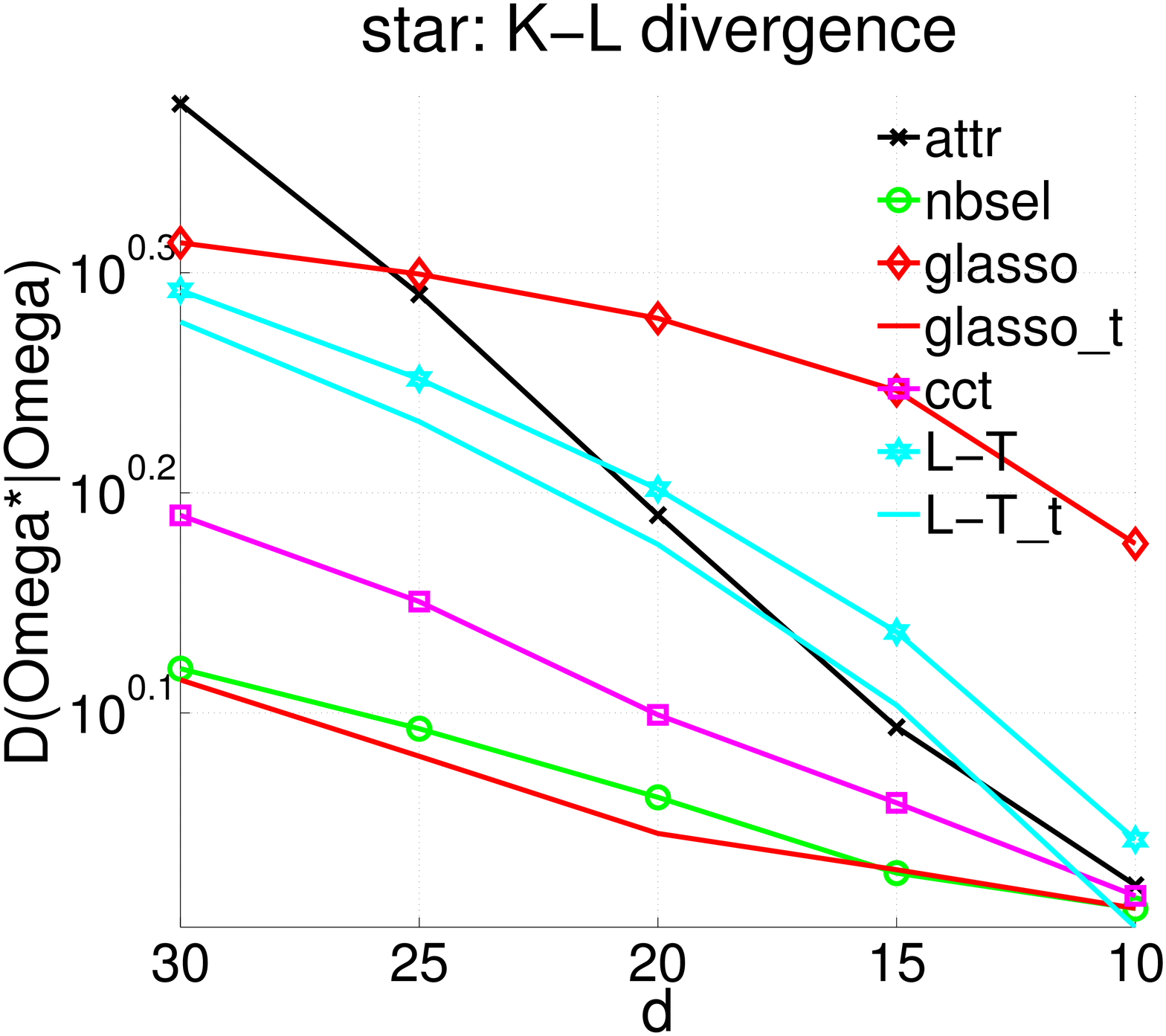} \\
\end{tabular}
\caption{Average performance over 50 replications. Left to right: performance measures; top to bottom: settings under consideration. For spectral norm error and
K-L divergence, a log-scale is used for both axes.}\label{fig:synthetic_benchmark}
\end{figure}

\paragraph{Data generation} We generate two sets of i.i.d.~samples $x_1,\ldots,x_n$ respectively 
$\wt{x}_1,\ldots,\wt{x}_n$ from a multivariate Gaussian distribution with mean zero and precision
matrix $\Omega_* \in \mc{M}^p$ according to one of the setups below. The first sample $\{x_i\}_{i=1}^n$ is used 
to obtain the sample covariance matrix $S = \frac{1}{n} \su x_i x_i^{\T}$ (assuming the mean to be known), 
which is the input for all methods under comparison. The second sample $\{ \wt{x}_i \}_{i=1}^n$ is kept
aside and used only for hyperparameter selection. The parameter $n$ is chosen in a setup-specific manner (see below).\\
\\    
\emph{chain.} We use an AR(1)-model (cf.~example (2) in Section \ref{sec:misspecification}) with positive parameter,
setting $\Sigma_* = (\sigma_{jk}^* ) = (0.9^{|j-k|})$, $j,k=1,\ldots,p=500$. The conditional
independence graph encoded by $\Omega_*$ is a chain.\\ 
\emph{grid.} We set $\wt{\Omega}_* = \delta I - B$, where $B$ is the adjacency matrix
of a 2d-grid (cf.~Figure \ref{fig:graphstructures}) of size $p = 23 \cdot 23 = 529$, and $\delta = 1.05 \lambda_1(B)$. We then 
set $\Omega_* = D \wt{\Omega}_* D$, where $D$ is a diagonal matrix chosen such that $\Sigma_* = \Omega_*^{-1}$ has unit
diagonal entries.\\
\emph{grid(3).} As for 'grid', but with $B$ replaced by the adjacency matrix of a 
3d-grid of size $p = 8 \cdot 8 \cdot 8 = 524$.\\    
\emph{random.} As for 'grid', but with $B$ replaced by a binary symmetric matrix of dimension $p = 500$ having 
zero diagonal and one percent non-zero off-diagonal entries, generated uniformly at random.\\
\emph{star.} We set $\Sigma_* = \begin{pmatrix}
1 & \rho^{\T} \\
\rho & I + \rho \rho^{\T}
\end{pmatrix}
\in \psd^{500}$ , where $\rho = 0.6 \cdot (\underbrace{1,\ldots,1}_{d \, \text{times}},0,\ldots,0) / d^{1/4}$ ($d \in \{10,15,20,25,30\}$), 
so that $\Omega_* = \begin{pmatrix}
D & -\rho^{\T} \\
-\rho & I
\end{pmatrix}
$, with $D = 1 + \norm{\rho}_2^2$.\\       
\emph{decay.} Unlike the previous setups, $\Omega_*$ is no longer sparse. Instead, its entries
exhibit an exponential decay away from the off-diagonal according to $\omega_{jk}^* = (-1)^{I(j \neq k)} \exp(-|j-k| \cdot 6/5)$, 
$j,k=1,\ldots,p=500$.\\
\\
Each setup is run for five different values of $n$ (with the exception of 'star', where $d$ varies while
$n$ is fixed to $500$, see Figure \ref{fig:synthetic_benchmark}). For each setup and each value of $n$, 50 replications 
are considered and performance is measured in the form of averages over these replications. To assess performance with regard to the recovery of the graph structure associated with $\Omega_*$ 
(or, equivalently, recovery of the set $\mc{E}^* = \{(j,k):\; \omega_{jk}^* < 0 \} $), we compute Matthew's correlation coefficient (MCC) defined 
by 
\begin{equation*}
\text{MCC} = (\text{TP} \cdot \text{TN} - \text{FP} \cdot \text{FN})/\left\{(\text{TP} + \text{FP})(\text{TP} + \text{FN})(\text{TN} + \text{FP})(\text{TN} + \text{FN}) \right\}^{1/2}, 
\end{equation*}
with $\text{TP}$,$\text{FN}$ etc.~denoting true positives, false negatives etc.~The larger
the criterion, the better the performance. Estimation of $\Omega_*$ is evaluated by the 
error in spectral norm $\nnorm{\Omega_* - \wh{\Omega}}$. We also report the Kullback-Leibler (KL) divergence
$D(\Omega_* \parallel \wh{\Omega}) = \log \det(\Omega_*) - \log \det(\wh{\Omega}) + \tr(\wh{\Omega} \Sigma_*) - p$.

\begin{figure}
\begin{center}
\begin{tabular}{cc} 
\includegraphics[height = 0.20\textheight]{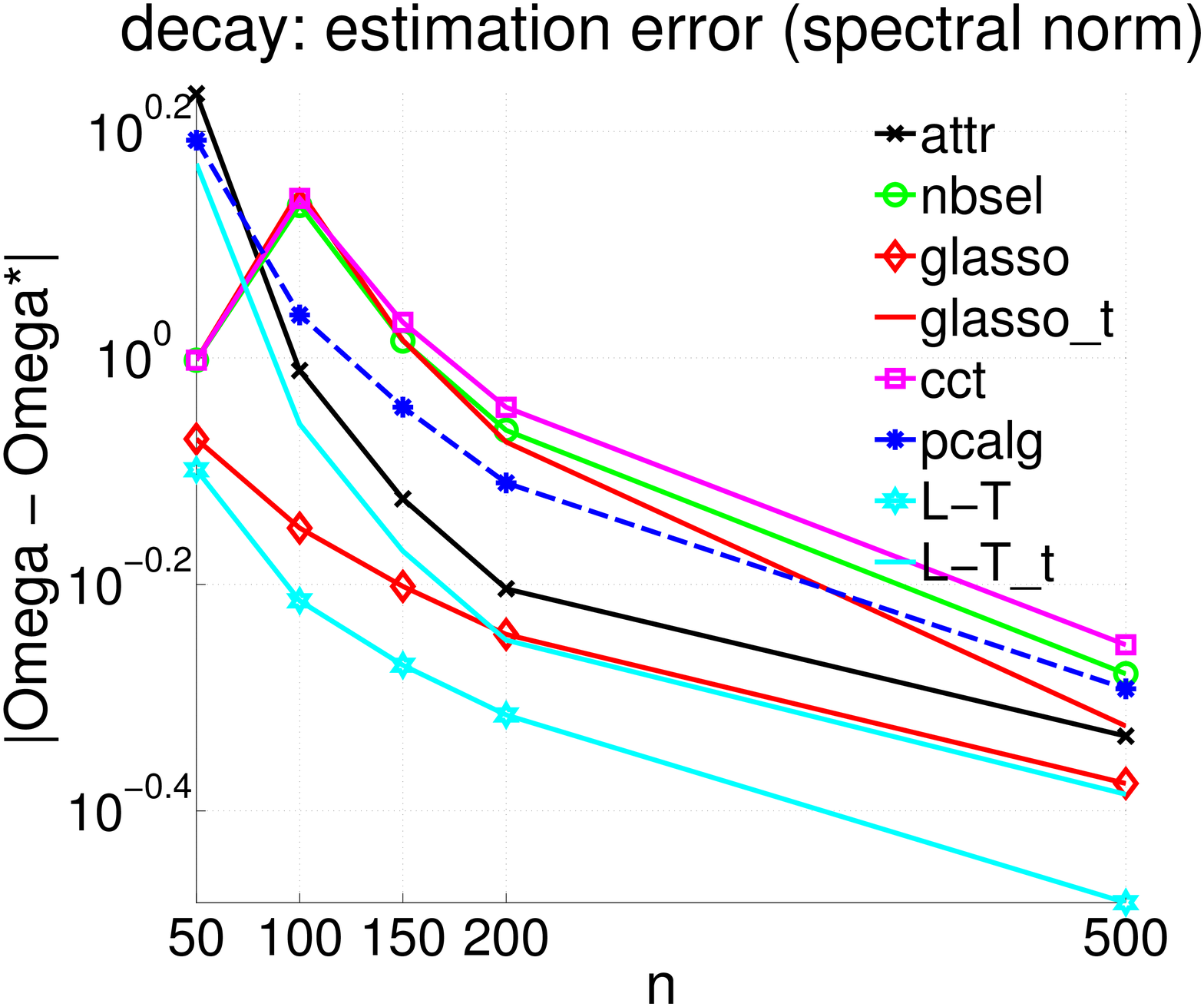} &
\includegraphics[height = 0.20\textheight]{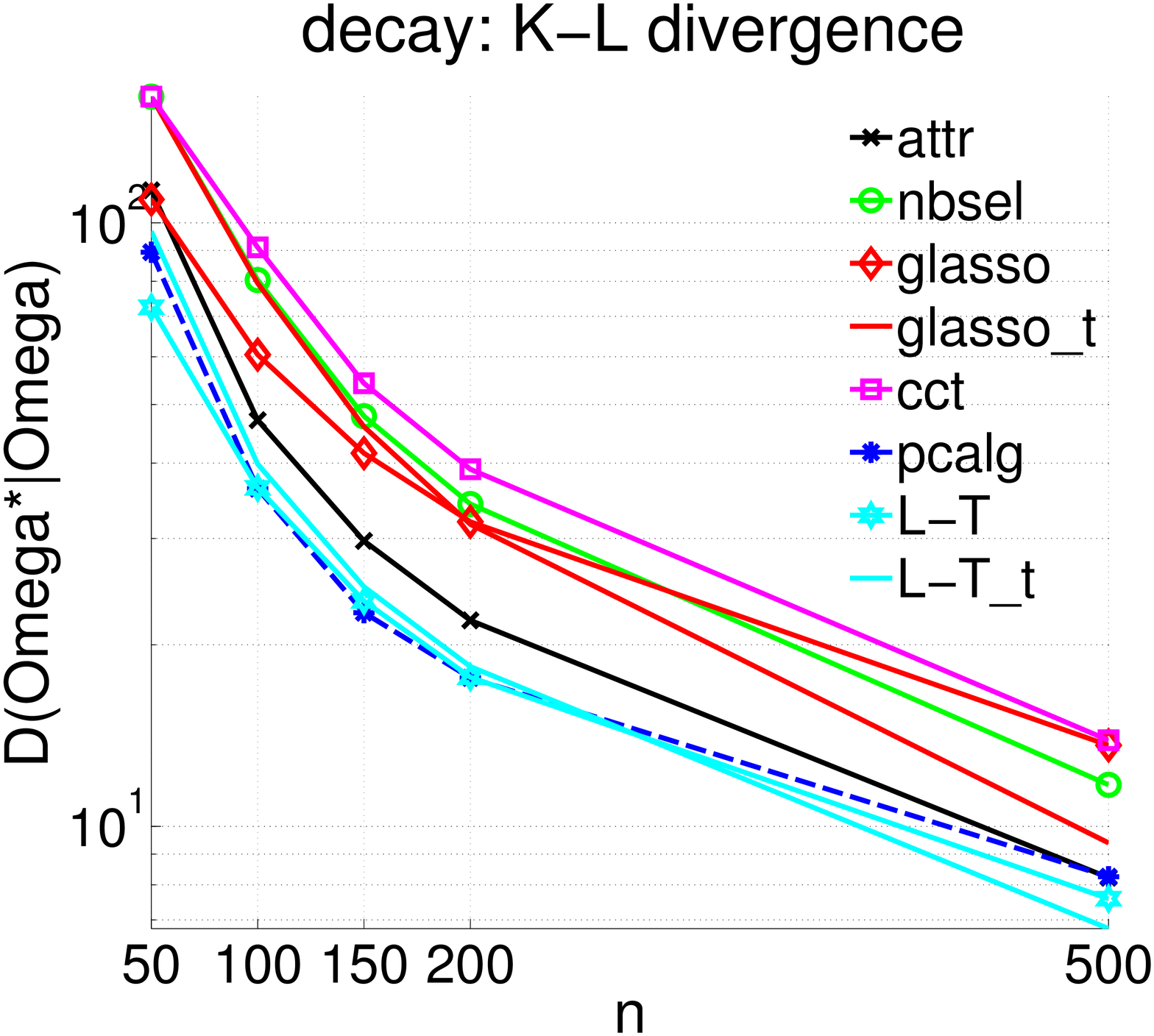} \\
\end{tabular} \hfill \\
\vspace{-.3cm}
\caption{Note that for the setting 'decay' $\Omega_*$ is not 
sparse; hence, we do not report MCC.}\label{fig:synthetic_benchmark_decay}  
\end{center}
\end{figure}

\paragraph{Methods compared} In order to thoroughly benchmark our approach, our
comparison includes various state-of-the-art methods in sparse precision matrix
estimation and structure learning of GMRFs.\\
\\
\emph{attr.} We use the abbreviation 'attr' (mnemonic for attractive random field) to
refer to the sign-constrained log-determinant divergence minimization plus thresholding approach
as described in Section \ref{sec:sparsification}. The threshold $t$ is chosen with the help of the separate
validation set $\{ \wt{x}_i \}_{i=1}^n$ in the following manner. Given the initial estimate 
$\wh{\Omega}$, we compute the $q$-quantiles of its off-diagonal entries, $q \in \{0.7,0.8,0.9,0.95, 0.99, 1, q^* \}$,
where $q^* = 1 - |\mc{E}^*|/(p \cdot (p-1)/2)$ is the quantile corresponding to the smallest non-zero off-diagonal
entry in $\Omega_*$. We then compute the loss on the validation set
\begin{equation}\label{eq:validation_threshold}
-\log \det(\wh{\wh{\Omega}}(t)) + \tr(\wh{\wh{\Omega}}(t) \wt{S}), \quad \; \, \wt{S} = \frac{1}{n} \su x_i x_i^{\T}. 
\end{equation}
based on the re-fitted estimator \eqref{eq:refitting}, with $t$ taken from the above quantiles, and pick 
the value of $t$ for which \eqref{eq:validation_threshold} is minimized.\\
\emph{glasso.} The g(raphical) lasso is defined as minimizer of the $\ell_1$-penalized log-determinant divergence           
\begin{equation*}
\min_{\Omega \in \psd^p} -\log \det(\Omega) + \tr(\Omega S) + \lambda \sum_{(j,k): j \neq k} |\omega_{jk}|, \quad \lambda  \geq 0.
\end{equation*}
Following \cite{Zhou2011}, the parameter $\lambda$ is chosen from the set $\{0.01, 0.05, 0.1, 0.3, 0.5,$
$1, 2, 4, 8, 16 \} \cdot \sqrt{(\log(p)/n)}$   
such that the loss on the validation set as in \eqref{eq:validation_threshold} is minimized. We also consider a thresholding plus refitting variant
of the glasso, denoted by \emph{glasso-t}. Thresholding and refitting proceeds as for 'attr' with $\wh{\Omega}(\lambda^*)$ as initial estimator, 
where $\wh{\Omega}(\lambda^*)$ denotes the glasso estimator with $\lambda = \lambda^*$ chosen as described above.\\
\emph{L-T.} \textsf{L}ake and \textsf{T}enenbaum \cite{LakeTenenbaum2010} consider the class of precision matrices $\mc{L}_I^{p}$ \eqref{eq:LpI} that can be written as a Laplacian matrix plus a multiple of the identity. Sparsity is promoted via an $\ell_1$-penalty
as for the 'glasso', which yields the following optimization problem. 
\begin{equation*}
  \min_{\Omega \in \mc{L}_I^p} -\log \det(\Omega) + \tr(\Omega S) + \lambda \sum_{(j,k): j \neq k} (-\omega_{jk}), \quad \lambda  \geq 0.
\end{equation*}
Recall that $\mc{L}_I^{p} \subset \mc{M}^p$ so that $\Omega \in \mc{L}_I^p \, \Rightarrow \omega_{jk} \leq 0, j \neq k$. The 
parameter $\lambda$ is chosen in the same manner as for the 'glasso'. Likewise, we consider a version with thresholding an refitting, denoted
by \emph{L-T-t}.\\
\\
The following three approaches only try to infer the graph structure of the conditional independence graph, i.e.
their output is an estimate $\wh{\mc{E}}$ of $\mc{E}^* = \{(j,k):\; \omega_{jk}^* < 0 \}$. The precision matrix
is estimated as
\begin{equation*}
\wh{\Omega} = \argmin_{\Omega \in \overline{\psd^p}, \; \omega_{jk} = 0 \;
  \forall (j,k) \notin \wh{\mc{E}}, \, j \neq k} -\log \det(\Omega) + \tr(\Omega S),
\end{equation*}
provided the minimizer exists, see e.g.~\cite{Uhler2012} for sufficient conditions.\\
\emph{nbsel.} \textsf{N}eig\textsf{b}orhood \textsf{sel}ection as proposed in \cite{Mei2006} tries
to infer the graph structure by node-wise $\ell_1$-penalized linear regression in which one
variable is regressed on the remaining ones. Pairs of variables are connected by an edge whenever
at least one of the two associated regression coefficients is nonzero. We use a refined version 
of neighborhood selection \cite{Zhou2011}, in which hard thresholding is applied to the node-wise 
regression coefficients. Following \cite{Zhou2011}, the regularization parameter for the node-wise
$\ell_1$-penalized regressions is chosen from the grid $\{0.01, 0.05, 0.1, 0.3, 0.5, 1, 2, 4, 8, 16 \} \cdot \sqrt{(\log(p)/n)}$ 
to minimize the mean squared prediction error on the validation set over all $p$ regression problems. Subsequent
thresholding is performed according to the scheme used for 'attr'.\\
\emph{cct.} In the \textsf{c}onditional \textsf{c}ovariance \textsf{t}esting approach of \cite{Anandkumar2012} 
one computes for all pairs $(j,k), \, j \neq k$, empirical conditional covariances        
$\wh{\sigma}_{jk|C} = s_{jk} - S_{jC}^{\T} S_{CC}^{-1} S_{Ck}$, where the conditioning set
$C$ ranges over all subsets $\mc{C}(\eta)$ of cardinality at most $\eta \geq 0$, which is the tuning parameter
of the approach. One then obtains $\theta_{jk} = \min_{C \subseteq \mc{C}(\eta)} \wh{\sigma}_{jk|C}$ and connects $j$ and
$k$ by an edge if $\theta_{jk}$ exceeds a suitable threshold (regarding the choice of threshold, we proceed as
for 'attr'). The computational complexity of the procedure is $O(p^{\eta + 2})$ and recovery of $\mc{E}^*$ according to the analysis in
\cite{Anandkumar2012} requires $\eta$ to be chosen as least as large as the size of the minimum vertex separator
in the conditional independence graph over all pair of edges. For the chain and star graph, $\eta = 1$, for the grid
$\eta = 2$ and for grid(3), $\eta = 3$ (which is not considered anymore for computational reasons), and 'cct' is
run with the correct choice of $\eta$ with knowledge about the underlying graph. For 'random' and 'decay', $\eta$ is set
to one.\\
\emph{pcalg.} The PC algorithm \cite{Spirtes2000} is an iterative procedure for inferring pairs of
variables of zero partial correlation. The approach has been further developed and analyzed in the
context of high-dimensional data in \cite{Kalisch2007}. In each iteration, a series of tests for
zero conditional covariances of increasing order, starting from marginal covariances, is performed.                    
Structural consistency requires \emph{faithfulness} of the underlying distribution \cite{Spirtes2000}. 
In Appendix G, we prove that if $\Omega_* \in \mc{M}^p$, faithfulness holds, which justifies the
use of the PC algorithm for the problem at hand. The significance level of the conditional independence
tests is chosen from the grid $\{ 0.001, 0.005, 0.01, 0.02, 0.05, 0.1 \}$ to minimize the loss on
the validation set. For the star graph, the PC algorithm is not run for computational reasons (in fact, its
complexity depends on the maximum vertex degree of the conditional independence graph).

\paragraph{Discussion: statistical performance} Inspecting Figures \ref{fig:synthetic_benchmark} and \ref{fig:synthetic_benchmark_decay},    
we find that the proposed approach performs competitively throughout; for two settings ('chain' and 'random'), it is among the
top competitors. This is rather remarkable, given the fact that our approach is applied in high-dimensional
settings without any explicit form of regularization employed in the first stage of estimation.  
A drop in performance is observed for 'star' as the vertex degree $d$ becomes larger, for 'decay', and
grid/grid(3) at the bottom end of the range considered for the sample size. The comparatively weak performance for 'star'
indicates a sub-optimal dependence of the approach on the maximum vertex degree, which is, besides the 
overall sparsity of $\Omega_*$, a second parameter known to affect performance of sparse precision matrix
estimation methods, see e.g.~\cite{Ravikumar2011, Zhou2011, Cai2011}.\\ 
A general conclusion one can draw from the figures is that all two-stage estimation procedures perform better
than the 'glasso' and 'L-T' (each without thresholding), excluding the non-sparse setting 'decay'.

\begin{table}
\hspace{-.4cm}
{\footnotesize
\begin{tabular}{|l||c|c||c|c||c|c||c|c||c|c||c|c||}
\hline
     & \multicolumn{2}{|c||}{chain} & \multicolumn{2}{|c||}{grid} & \multicolumn{2}{|c||}{grid(3)} & \multicolumn{2}{|c||}{random} & \multicolumn{2}{|c||}{star} & \multicolumn{2}{|c||}{decay} \\
\hline
$n/d$     & $20$ & $40$ &  $75$ & $250$ & $200$ & $500$ & $300$ & $1000$ & $30$ & $10$ & $50$ & $200$ \\
\hline 
attr & 9.4 & 6.4 & 2.7 & 2.6 & 1.3 & 1.3 & 1.5 & 1.2 & 0.6 & 0.6 & 0.8 & 1 \\
nbsel & 6   & 6.2 & 5.8 & 4.7 & 4.4 & 4.1 & 4.1 & 2.6 & 4 & 4.1 & 7.6 &7.1 \\
glasso & 8.2 & 10.1 & 12 & 11.9 & 11.4 & 4.4 & 10.9 & 3.9 & 3 & 3 & 12 & 14.6 \\
L-T & 6.7 & 8.5  & 7.4 & 5.2& 6.8 & 5.5 & 6.4 & 5.1 & 8.8 & 9.7 & 7.9 & 5.5 \\
cct & 16.3 & 16.5 & 37.5 & 38.7 & $-$ & $-$ & 31.1 & 32 & 21.8 & 21.8 & 34.9 &35.4 \\
pcalg & 0.6  &  0.2 & 0.2  & 0.2 & 0.2& 0.4 & 0.2 & 1.4 & 0.2 & 0.2 & 0.2 & 0.2 \\
 \hline
\end{tabular}
}
\vspace{-.3cm}
\caption{Median running times (in 100s) for solving one problem instance for the lower and upper ends of the
range for $n$ (resp.~$d$).}\label{tab:runtimes}
\end{table}

\begin{figure}[t!]
\begin{center}
\begin{tabular}{cc}
\includegraphics[height = 0.2\textheight]{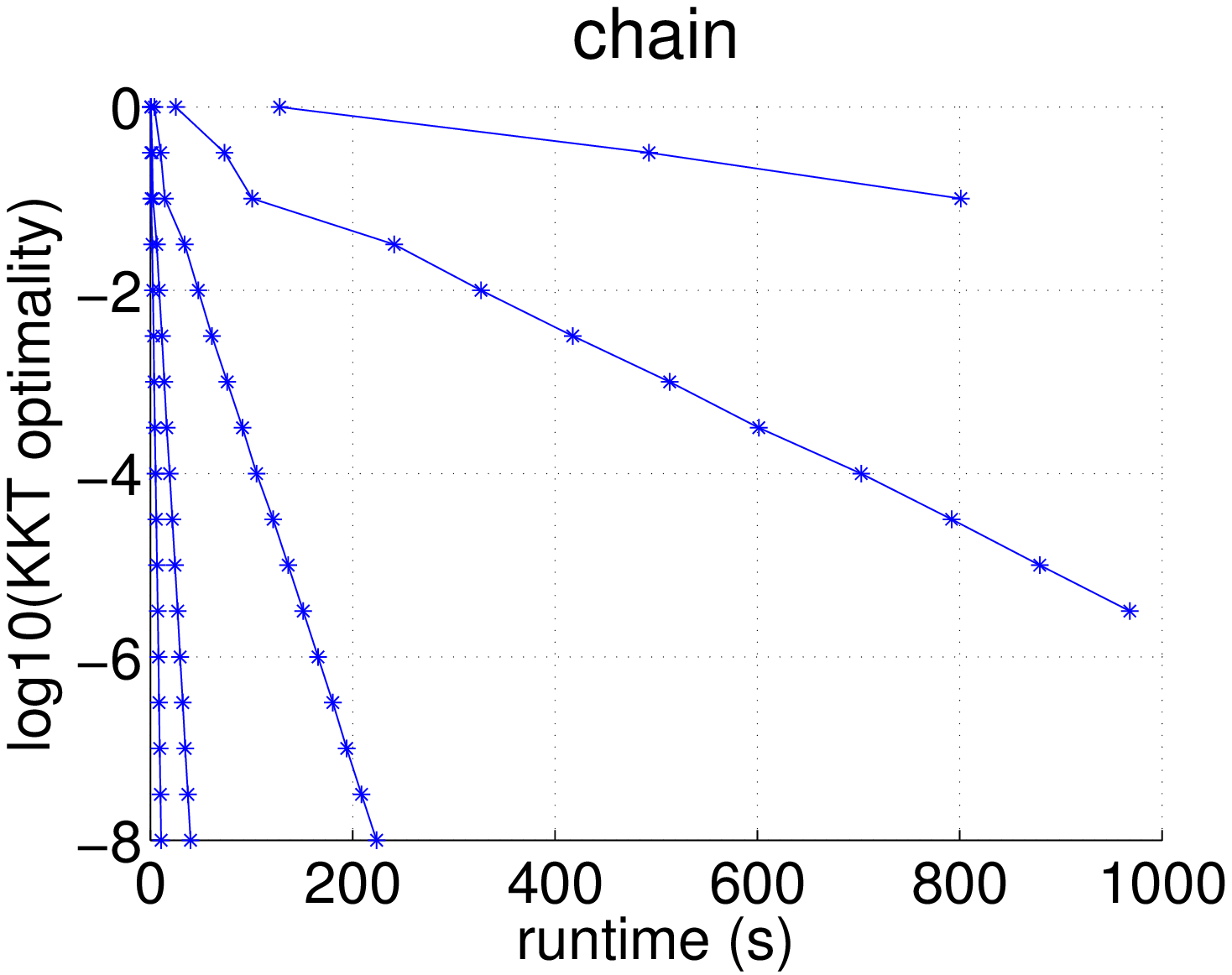} & 
\includegraphics[height = 0.2\textheight]{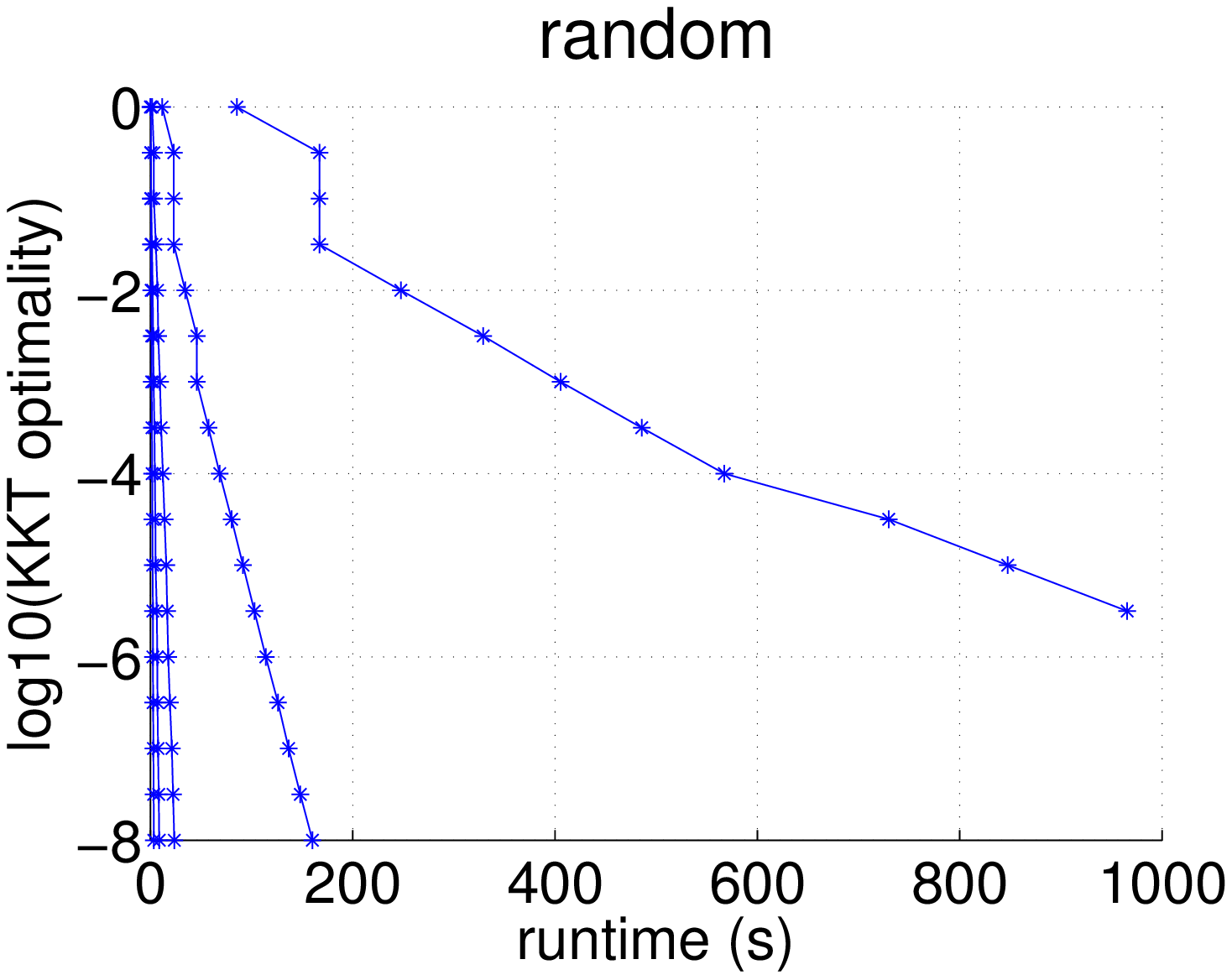} \\
\includegraphics[height = 0.2\textheight]{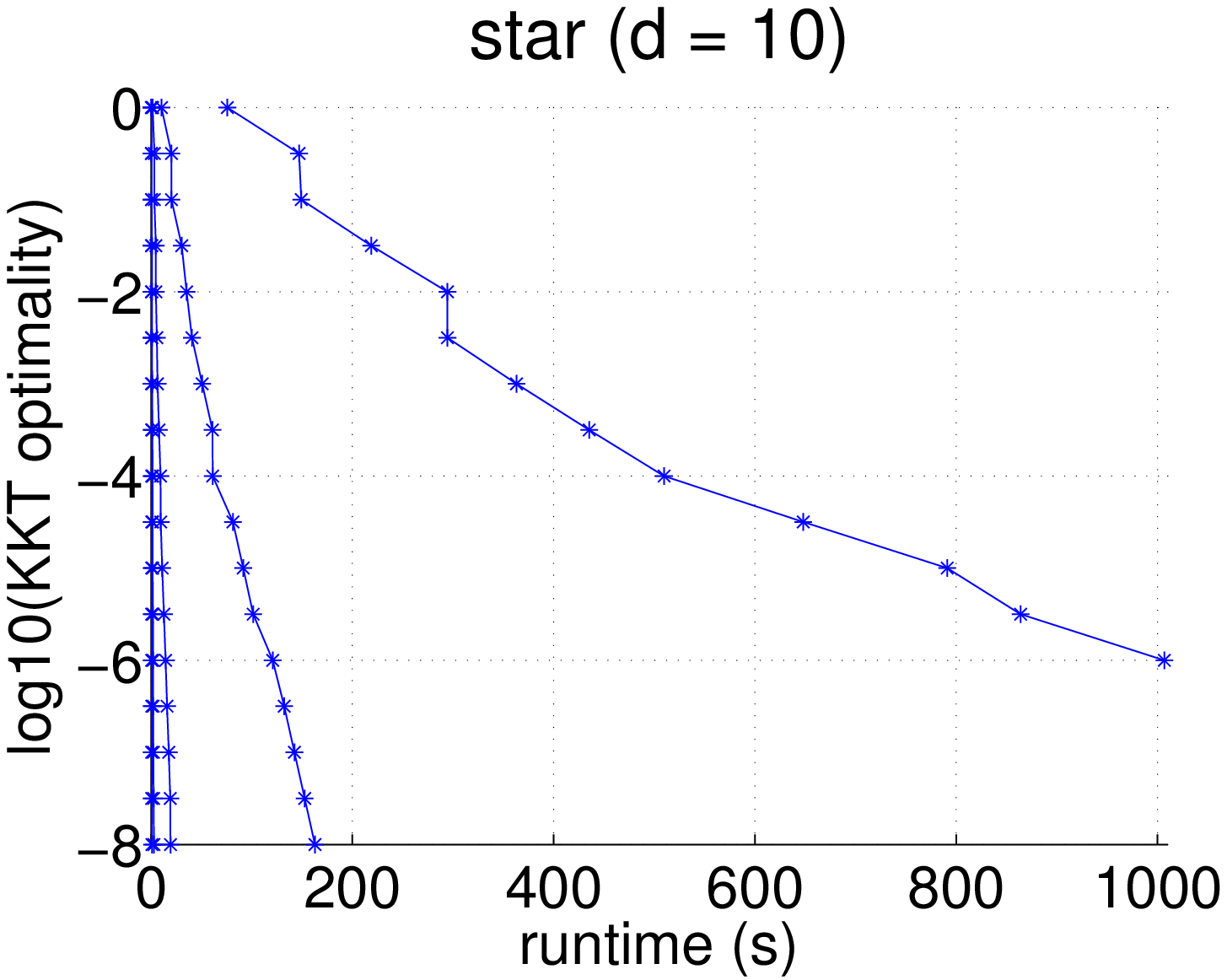} & 
\includegraphics[height = 0.2\textheight]{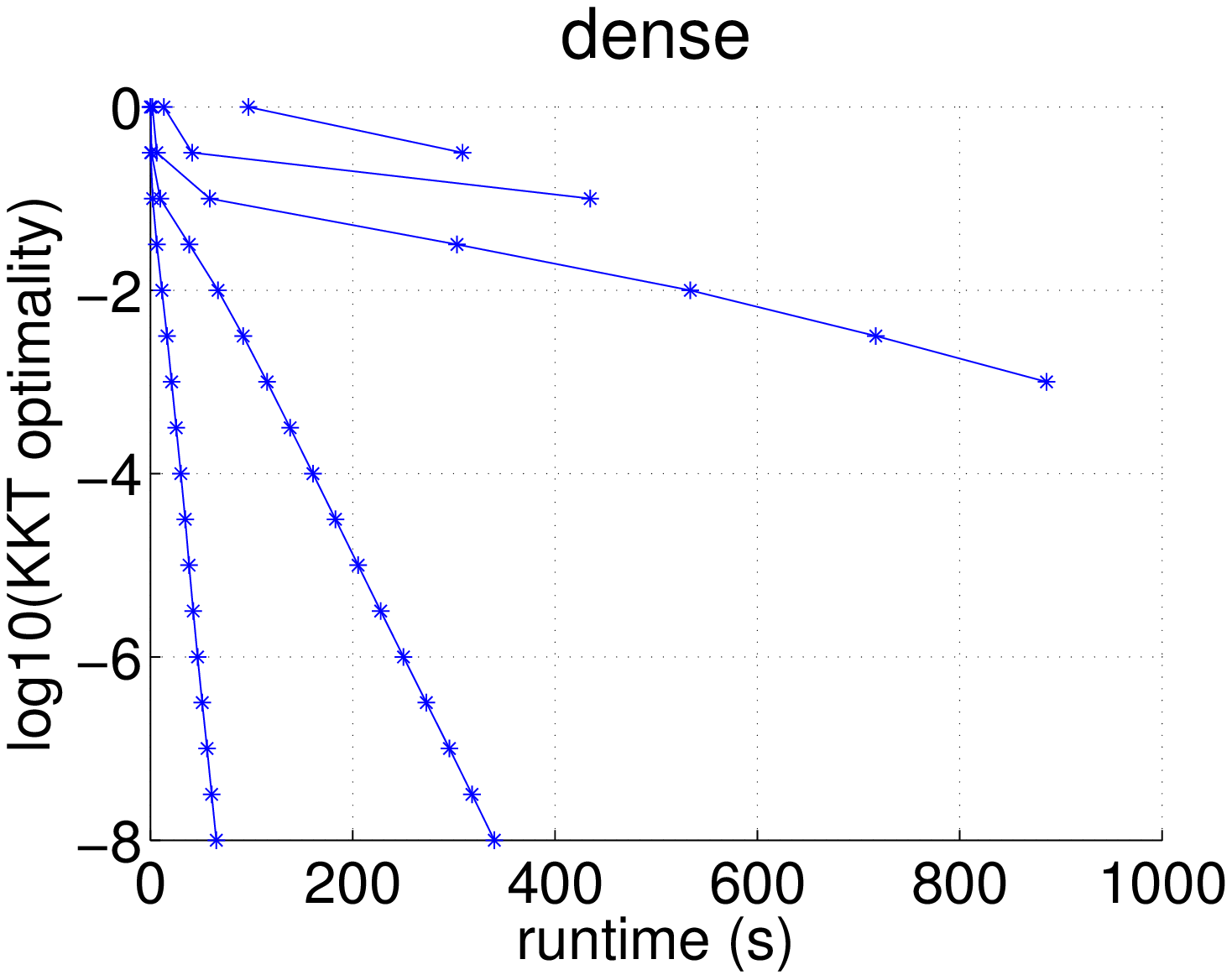} 
\end{tabular}
\end{center}
\vspace{-.5cm}
\caption{Empirical assessment of the speed of convergence of the block coordinate descent algorithm 
  of Section \ref{sec:optimization}. The five different trajectories per plot represent different $p$'s ranging from $2^5$ to
$2^{10}$, with faster convergence corresponding to smaller $p$.}\label{fig:runtimes}
\end{figure}
\paragraph{Discussion: running times} The method of choice should not only have good statistical properties, 
but should as well be favourable in terms of computation. For this reason, we report the results of 
a small runtime comparison in Table \ref{tab:runtimes}. The running times reported there only refer
to the first stage of estimation (including possible hyperparameter tuning), but excluding possible re-estimation steps
given an estimate of the graph structure. Apart from 'pcalg', for which we use the \textsf{R} implementation \cite{Kalisch2012},
all methods are run under \textsf{MATLAB}. For 'glasso' and 'cct' publicly available code is used \cite{glasso_matlab, Anandkumar2012_supp}.
We use code from \cite{L1general} to solve the $\ell_1$-penalized regression problems arising for 'nbsel'. For
'attr', we have implemented the coordinate descent approach of Section \ref{sec:optimization}. An own algorithmic approach based on
projected gradient is used for L-T, which will be reported elsewhere. Table \ref{tab:runtimes} reveals 
that while 'cct' is a theoretically sound approach which has also been seen to perform well empirically, it
falls short in terms of runtime. The comparison is headed by 'pcalg', which outperforms all competitors
by one order of magnitude; however, its computational complexity depends on the maximum degree of the graph,
which becomes a severe issue for the setting 'star'. On average, 'attr' has smaller running times than the 
methods employing regularization ('nbsel', 'glasso' and 'L-T'), since no hyperparameter has to be set. 
The speed of convergence of our computational approach has been investigated in a separate series of
experiments whose results are displayed in Figure \ref{fig:runtimes}. For each of the four settings with $p$ ranging
from $2^6$ to $2^{10}$, we measure the time required to achieve a certain level of KKT optimality as given in \eqref{eq:kkt_practice}.
The times reported in the figure are medians over ten replications obtained per setting. We use the following
sample sizes: $n = 40$ for 'chain', $n = 1000$ for 'random', and $n = 500$ for 'star' and 'dense'. For 'dense',
the only of the four settings which is exclusively considered for the runtime analysis, $\Omega_* = (0.3 \cdot I + 0.7 \bm{1} \bm{1}^{\T})^{-1}$.  
Figure \ref{fig:runtimes} suggest that our approach exhibits a linear rate of convergence, but does not
scale well with $p$ as already indicated in Section \ref{sec:optimization}. Moreover, the speed of convergence visibly depends 
on the structure of $\Omega_*$, with 'dense' being the most difficult setting and 'star' as well as 'random'
being easier than 'chain'.    
 
\subsection{Learning Taxonomies}

\begin{figure}
\begin{flushleft}
\begin{tabular}{c}
\includegraphics[height = 0.425\textheight]{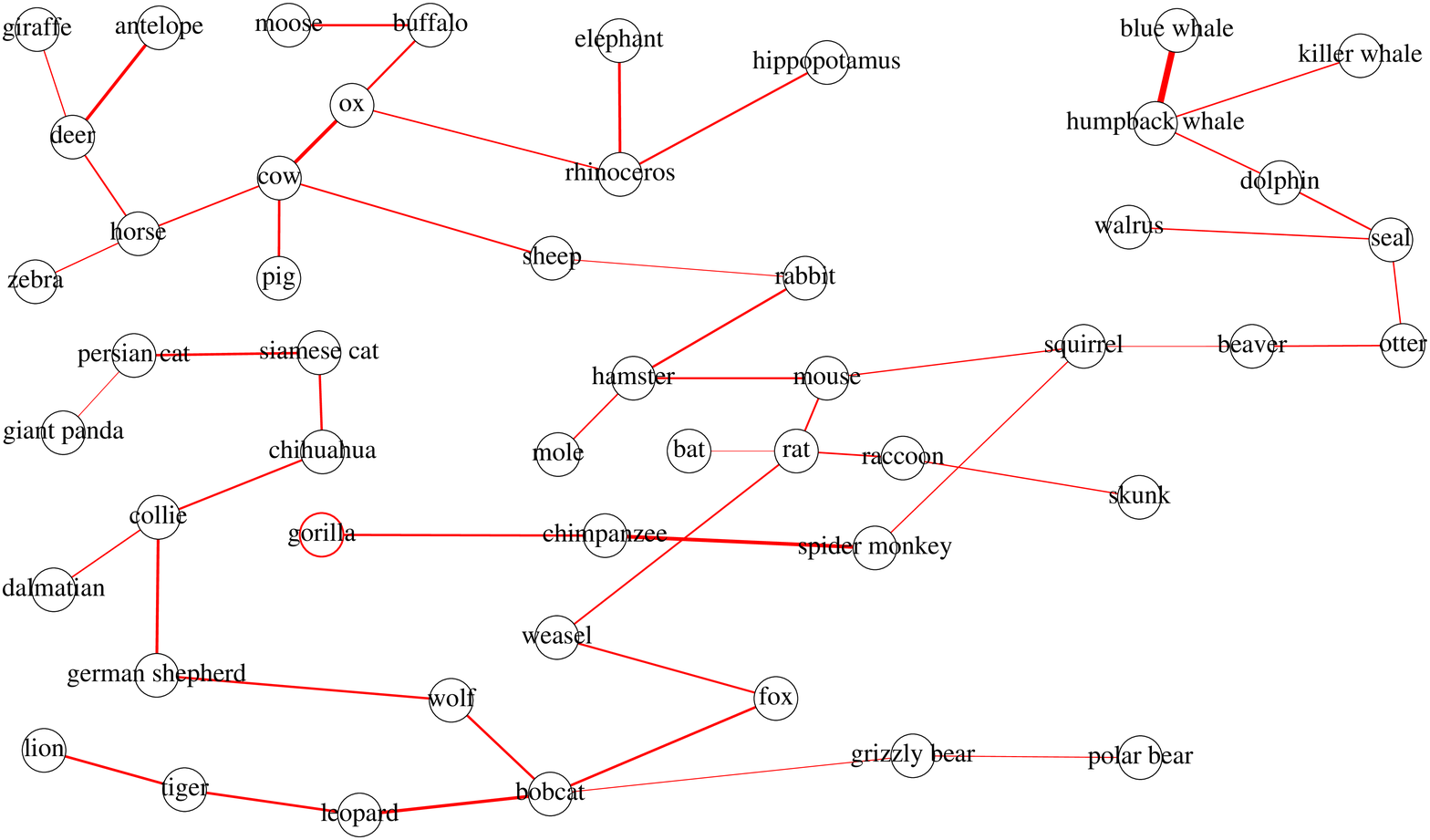} \\
\hline
\includegraphics[height = 0.425\textheight]{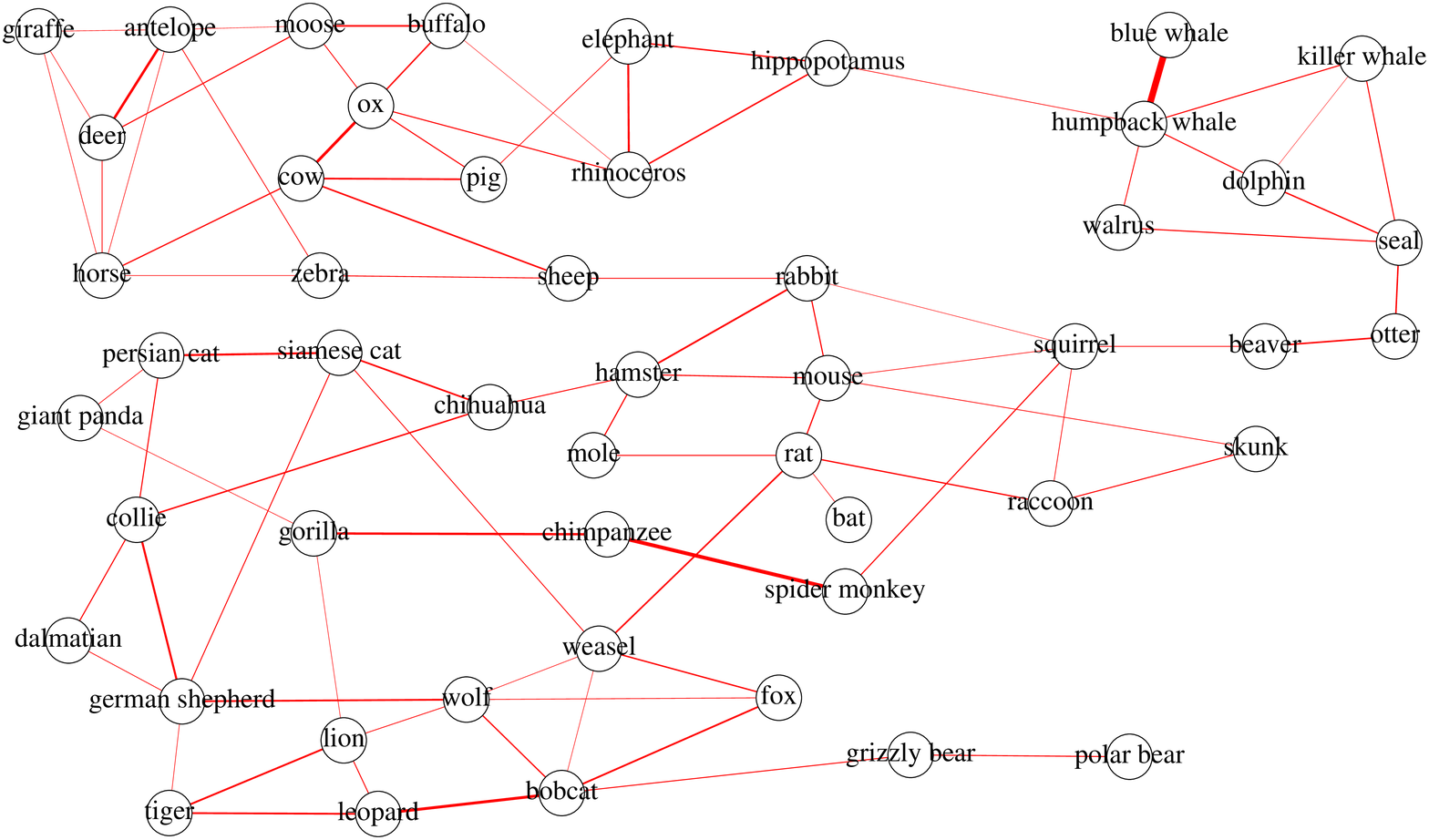}
\end{tabular}
\end{flushleft}

\caption{Graphs for the mammals dataset. Top: graph of the tree model having $p-1 = 49$ edges (all of which are associated with
positive partial correlations) and a median cross-validated loss of $15.63$. Bottom: graph of the attractive model 
after thresholding having $83$ edges and a median
cross-validated loss of $10.88$, with the threshold  set to the 0.47-quantile of the absolute values of the non-zero off-diagonal entries in $\wh{\Omega}$. Edge widths are proportional to the absolute values of the corresponding entries
in the precision matrices.}\label{fig:graphs_animals}
\end{figure}

In the next two paragraphs, we conduct an analysis performed in \cite{LakeTenenbaum2010} 
where precision matrices of the form \eqref{eq:LpI} are considered. 
\vspace{-.25cm}
\paragraph{Mammals dataset} The dataset contains $n = 85$ biological properties of $p = 50$ mammals. The dataset 
is the outcome of a study \cite{Osherson1991}, in which participants were asked to provide scores quantifying the 
relative strength of association between each mammal and the set of biological properties concerning
anatomy, behaviour and living conditions on a scale ranging from 0 (no association) to 100 (perfect association). 
This yields a data matrix $X = (x_{ij})$ with $x_{ij}$ as the mean relative strength of association between
property $i$ and mammal $j$, and in turn a sample covariance $S = \wt{X}^{\T} \wt{X} / n$, where $\wt{X}$ results
from $X$ by centering its columns. The goal is to use the given data to infer a graphical
representation of the $50$ mammals serving as a taxonomy. We here
compare approaches to sparse precision matrix estimation and the resulting graphs associated
with the \emph{negative} off-diagonal elements. Note that it is not meaningful to include edges
corresponding to positive off-diagonal entries, since these edges would not be interpretable
in the context of taxonomic reasoning. It still may make sense not to impose sign constraints 
in estimation, and to construct the graph only from the negative entries, because as discussed
in Section \ref{sec:misspecification}, sign-constrained estimation may lead to a bias even for the underlying negative
entries. We here compare the graphical lasso, Tikhonov regularization, a tree model, thresholding of $S^{-1}$, the method suggested in \cite{LakeTenenbaum2010} ('L-T'), and our sign-constrained approach, among which only the latter two impose sign constraints. 
Tikhonov regularization provides an estimate $\wh{\Omega} = ((1-\alpha) S + \alpha I)^{-1}$, where
$\alpha \in [0,1]$ is a tuning parameter. For the tree model, we restrict the graph associated
with the precision matrix to be a tree, and the estimate  is obtained with the help of the
Chow-Liu algorithm \cite{ChowLiu1968}. In order to judge the usefulness of the sign-constraints, we include thresholding and
re-fitting as described in Section \ref{sec:sparsification}, with the difference that $S^{-1}$ is used as initial
estimator. 
\begin{figure}
\leftline{\includegraphics[height = 0.25\textheight]{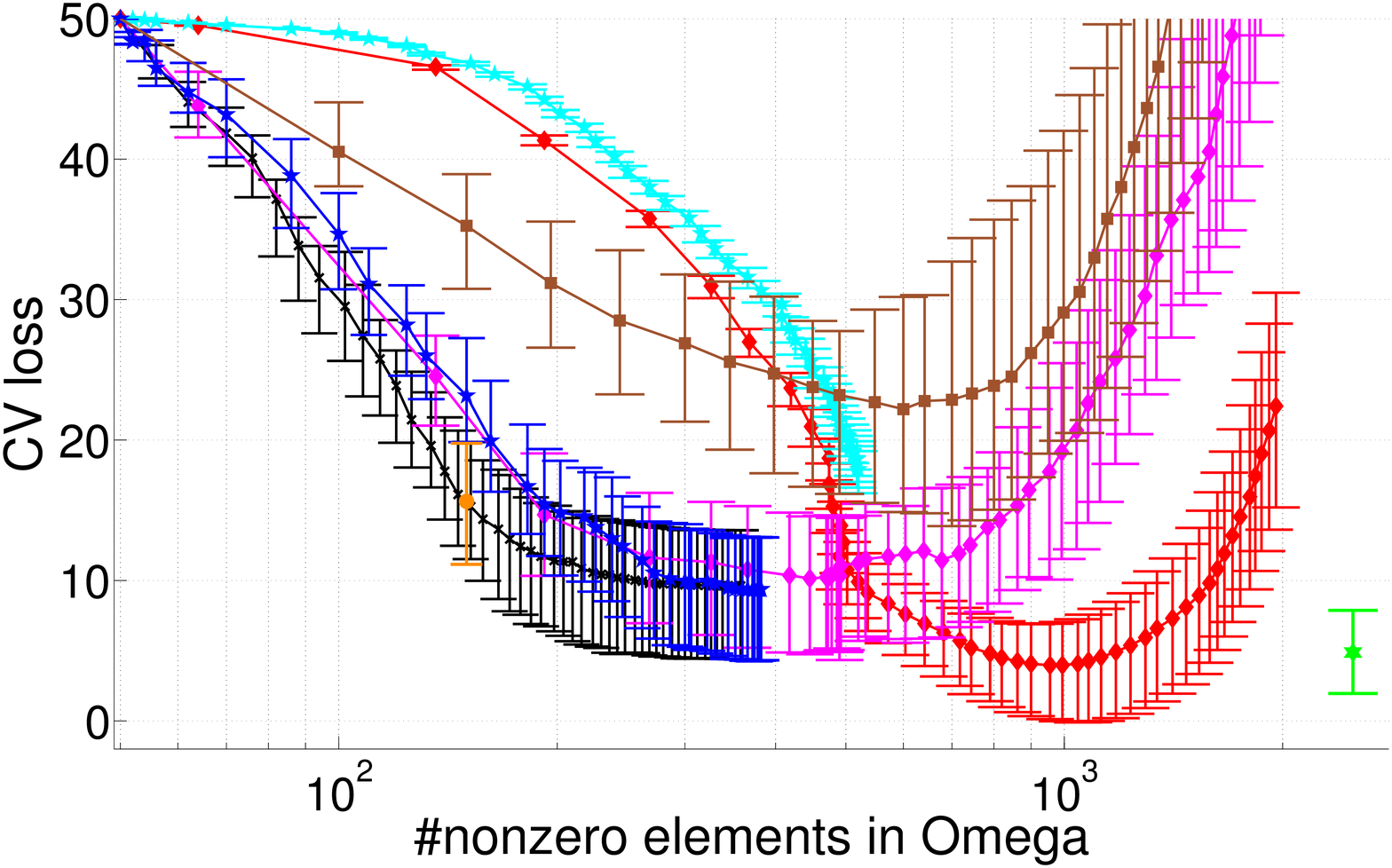}} 
\vspace{-4.5cm}\rightline{\hspace{-0.4cm}\includegraphics[height = 0.17\textheight]{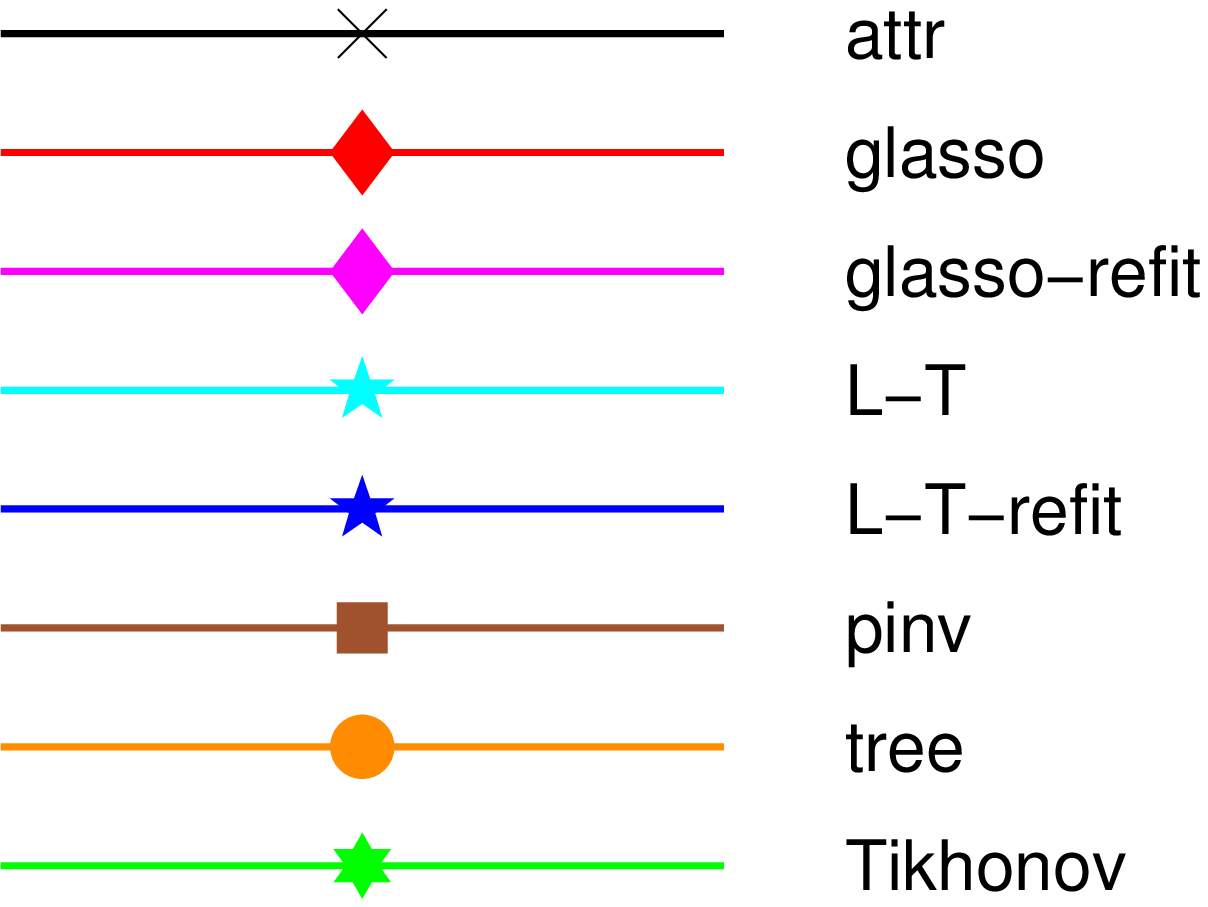} \vspace{2cm}}
\vspace{-1.1cm}
\caption{Cross-validated loss on the mammals dataset in dependency of the sparsity of the estimates. Upper and
lower bars indicate the maximum respectively minimum loss over the ten random partitionings of the set
of observations into ten folds, while points indicate medians. For the tree model ('tree') and Tikhonov regularization ('Tikhonov') the number
of non-zeroes is constant, and three lines/points depict minimum, median (thick) and maximum. 'attr': attractive (sign-constrained) estimation; 
'glasso(-refit)': graphical lasso (with re-fitting given the zero pattern); 'L-T(-refit)': method of Lake and Tenenbaum (with re-fitting given the zero pattern); 
'pinv': thresholding of the (pseudo)-inverse.}\label{fig:cv_animals}
\end{figure} 
Performance is quantified by computing the cross-validated loss defined analogously to \eqref{eq:validation_threshold}.
Regarding cross-validation, we consider ten folds and ten random partitionings into folds.  
We consider 50 different values of the respective tuning parameters corresponding to varying
levels of sparsity of the estimates (except for Tikhonov regularization, which yields fully dense
estimates, and we report the minimum loss over all choices of $\alpha$) as displayed in Figure \ref{fig:cv_animals}.
Our sign-constrained approach performs best for high levels of sparsity, which is the regime
of interest here. Interestingly, in the range of 150 to 200 non-zero entries, both cross-validated
loss and graph come rather close to the tree model (see Figure \ref{fig:graphs_animals}), which
is conventionally used for depicting taxonomic relationships. In general, taxonomies are not required to have a tree structure. As pointed out in \cite{LakeTenenbaum2010}, the data set under consideration contains several features concerning e.g.~habitat or appearance that may induce associations between species that would not be linked in an evolutionary tree. Unlike the tree model, an attractive model is not constrained in the number of edges and may hence constitute a more flexible alternative; this has to be kept in mind when comparing the cross-validated loss of the two approaches. All edges of the Chow-Liu tree are associated with positive partial correlations, which supports
the use of an attractive model. The latter yields several extra edges relative to the tree, some of which are less intuitive   
(elephant and pig, hamster and chihuahua, weasel and wolf). The sign constraints appear to be limiting in the sense that a visibly lower
cross-validated loss can be achieved by methods not imposing sign constraints. However, this
concerns a regime in which the estimates are no longer sparse and thus less interpretable.  
\vspace{-.25cm}

\paragraph{Concepts dataset} The analysis in the present paragraph is of the same spirit as
the previous one, with the important difference that the dataset under consideration is
high-dimensional ($n = 218$, $p = 1000$). The goal is to infer a semantically meaningful graph-based
representation of $1000$ concepts falling into a diverse set of categories such
as food, buildings, animals, clothes or other consumer goods, with the help of answers given
to 218 questions concerning various attributes. The answers are on a five-point scale from
'clearly no' to 'clearly yes', obtained from Amazon Mechanical Turk \cite{mturk}. For simplicity, we treat the data as if the scale
level were metric. The analysis is conducted in the same way as the previous one, with two minor
modifications. In order to reduce computing times, five-fold cross-validation is used and
only 20 different values are considered for the tuning parameters.
\begin{figure}
\leftline{\hspace{-0.6cm}\includegraphics[height = 0.27\textheight]{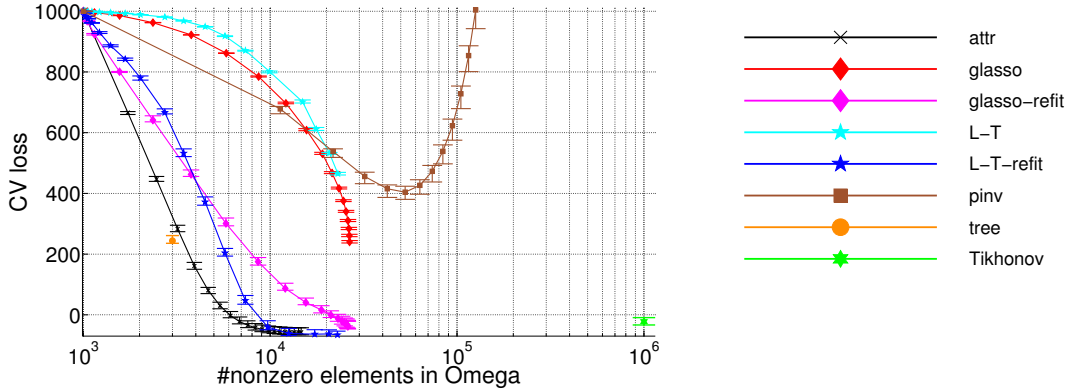}} 
\vspace{-5cm}\rightline{\includegraphics[height = 0.18\textheight]{animals_cv_legend.eps} \vspace{2cm}}
\vspace{-.5cm}
\caption{Cross-validated loss on the concepts dataset in dependency of the sparsity of the estimates. The annotation is as for Figure \ref{fig:cv_animals} above.}\label{fig:cv_concepts}
\end{figure}     
\begin{figure}
\vspace{-4.1cm}
\begin{center}
\includegraphics[height = 0.60\textheight, angle = -90]{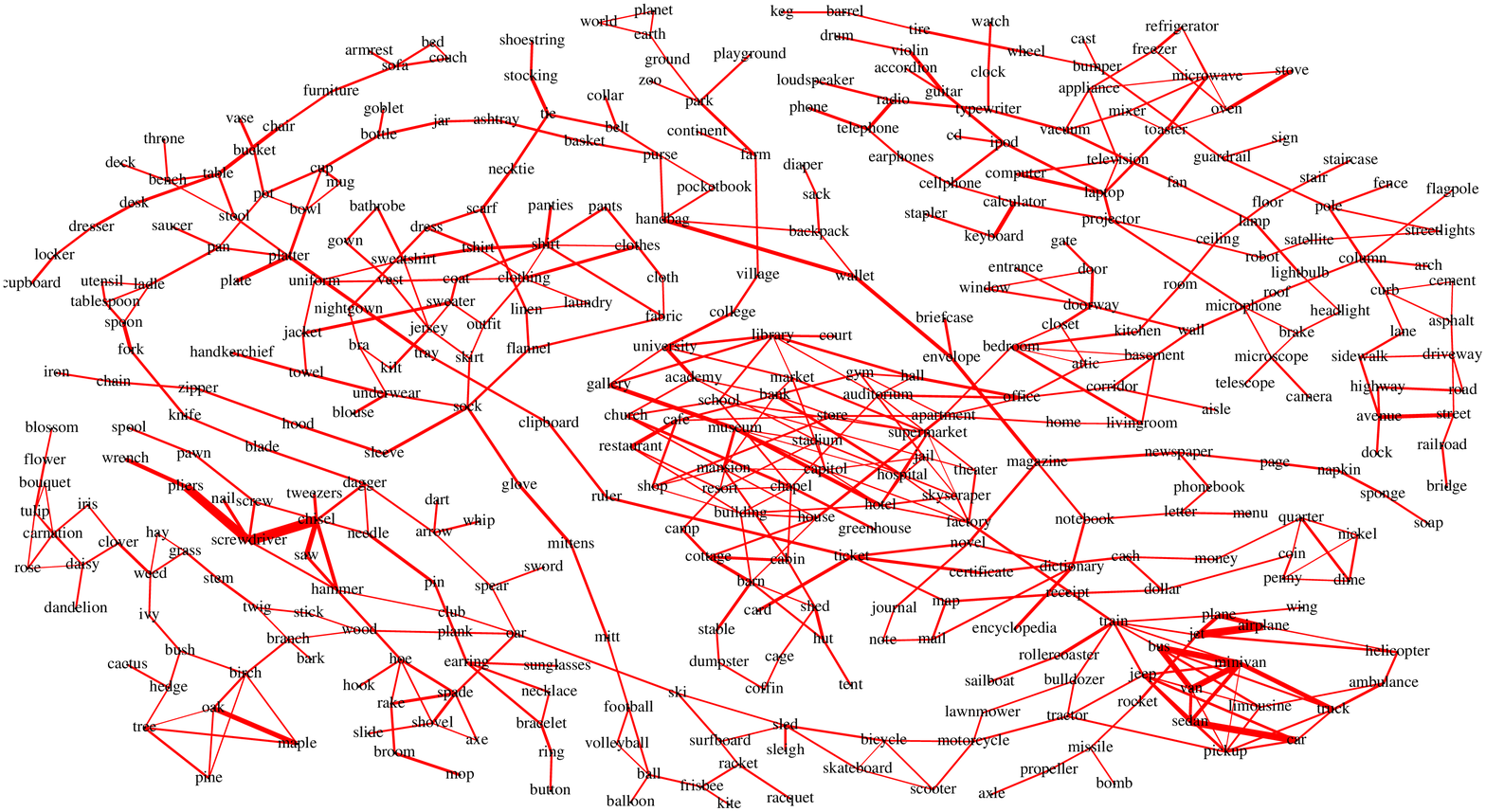}
\end{center}
\vspace{-.7cm}
\caption{Largest connected component (377 vertices) of the attractive model (after thresholding with the 0.84-quantile  of the absolute values of the non-zero off-diagonal entries) for the concepts dataset. 
Edge widths correspond to the absolute values in $\wh{\Omega}$.}\label{fig:graph_concepts}
\end{figure}

\clearpage

\subsection{Covariance modeling of landmark data}

\begin{figure}
\hspace{-0.4cm}
\begin{minipage}[t]{0.25\textwidth}
\vspace{-2cm}\includegraphics[height = 0.18\textheight]{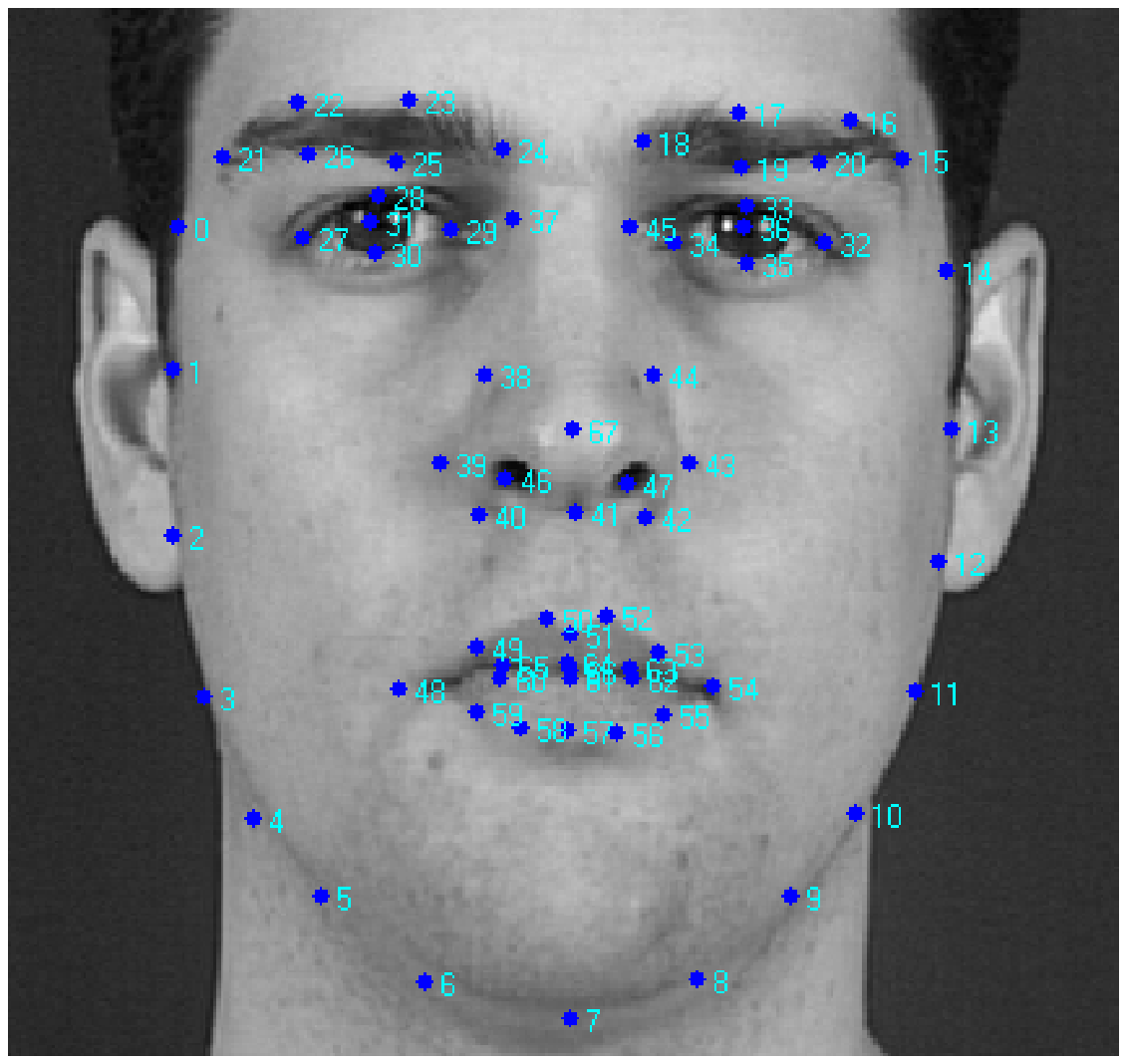}
\end{minipage}
\begin{minipage}{0.7\textwidth}
\begin{tabular}{ll}
\hspace{-0.15cm}\includegraphics[height = 0.23\textheight]{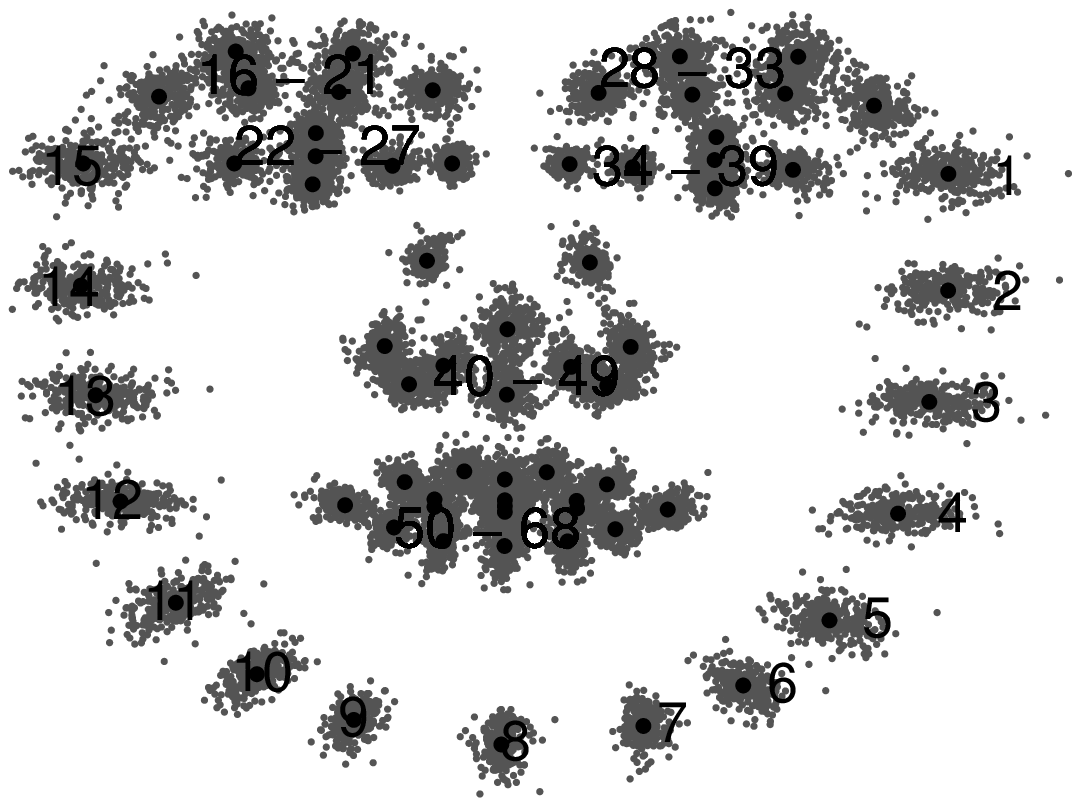} &
\hspace{-0.75cm}\includegraphics[height = 0.23\textheight]{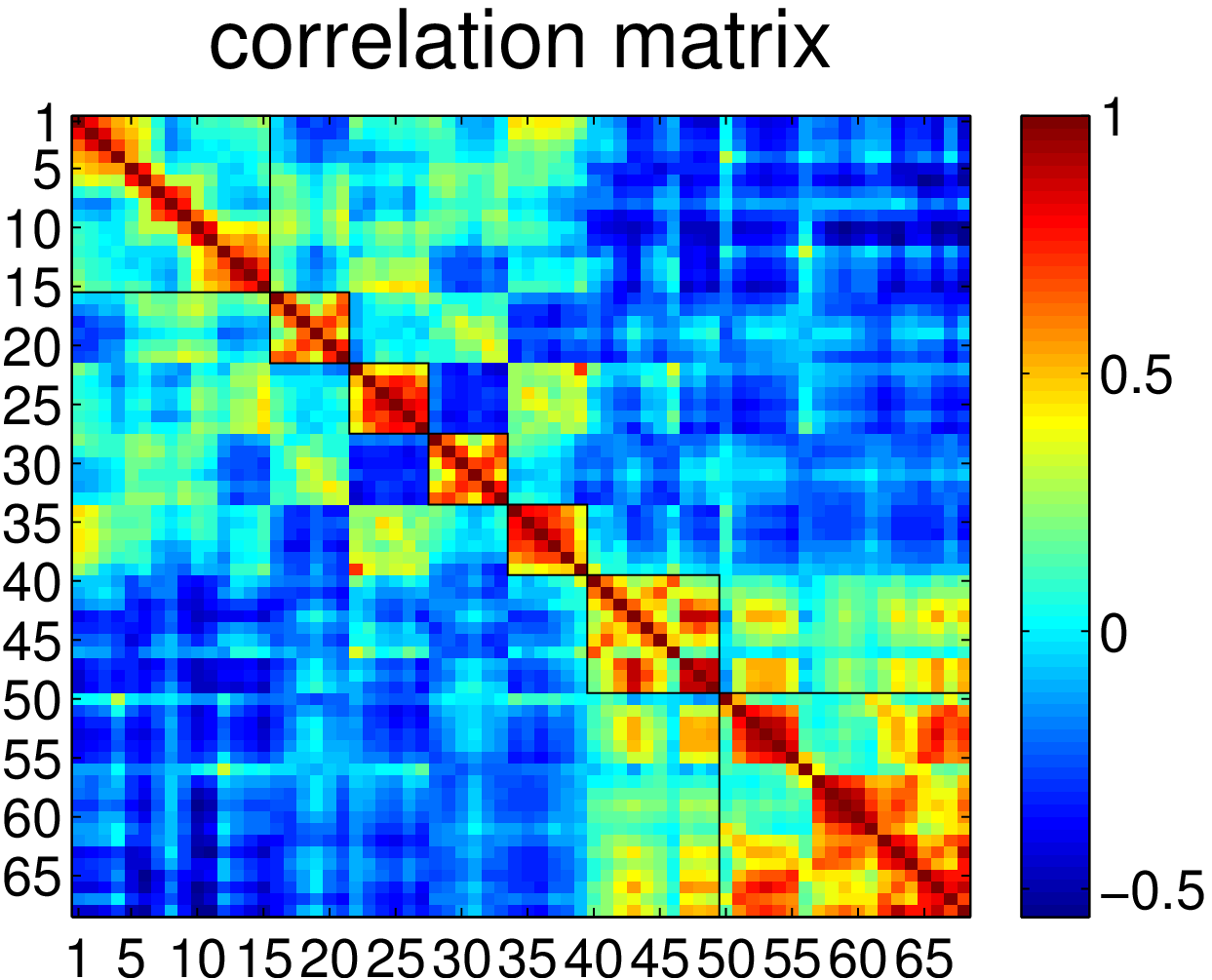}
\end{tabular}
\end{minipage}
\caption{Landmarking scheme (left), scatterplot of landmark points and their numerical labels (middle; contour: landmarks 1-15, eyebrows: landmarks 16-21 and 28-33, eyes: landmarks 22-27 and 34-39, nose: landmarks 40-49, mouth: landmarks 50-68) and correlation
matrix corresponding to the covariance matrix $S^{\textsf{lm}}$ of the training set. The squares indicate
the blocks corresponding to the landmark groups for contour, eyebrows, eyes, nose and mouth.}\label{fig:landmarks}
\end{figure}

\begin{figure}[t!]
\begin{tabular}{ccc}
\includegraphics[height = 0.15\textheight]{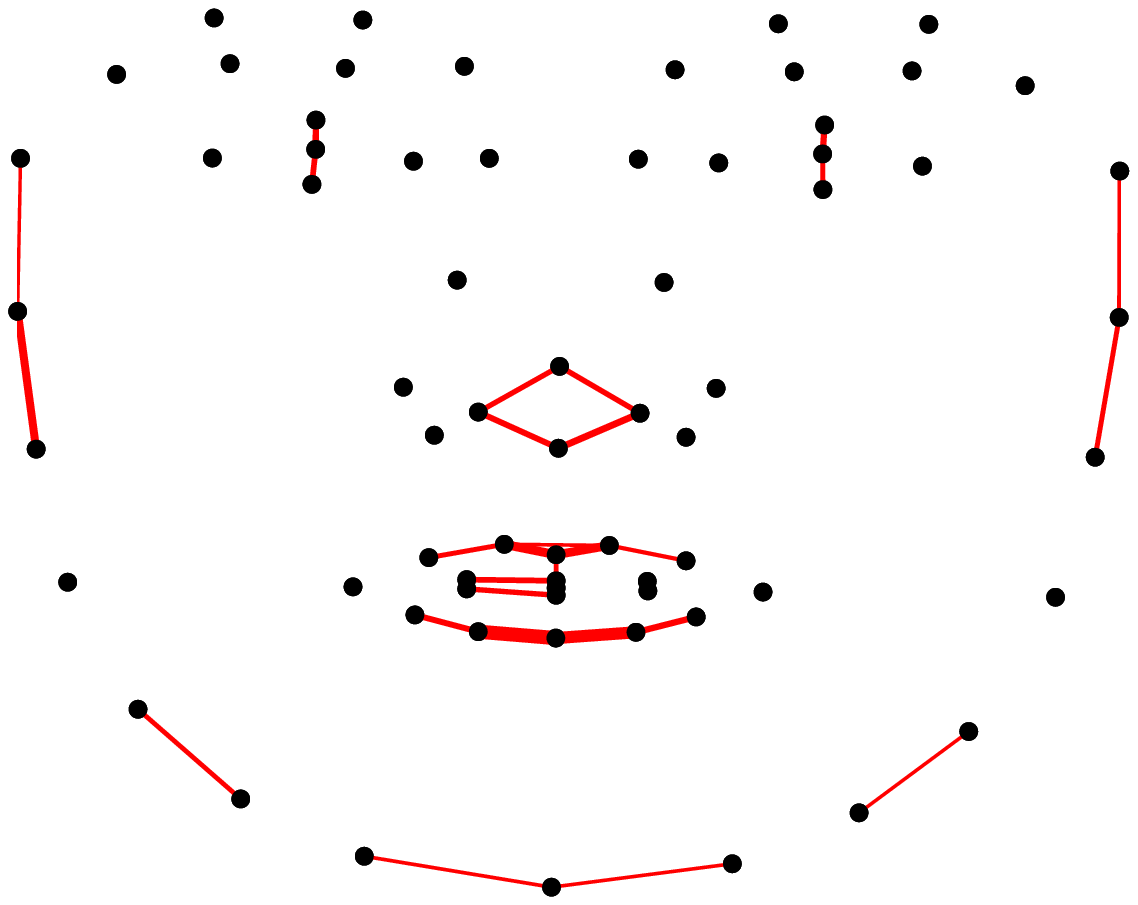} &  \includegraphics[height = 0.15\textheight]{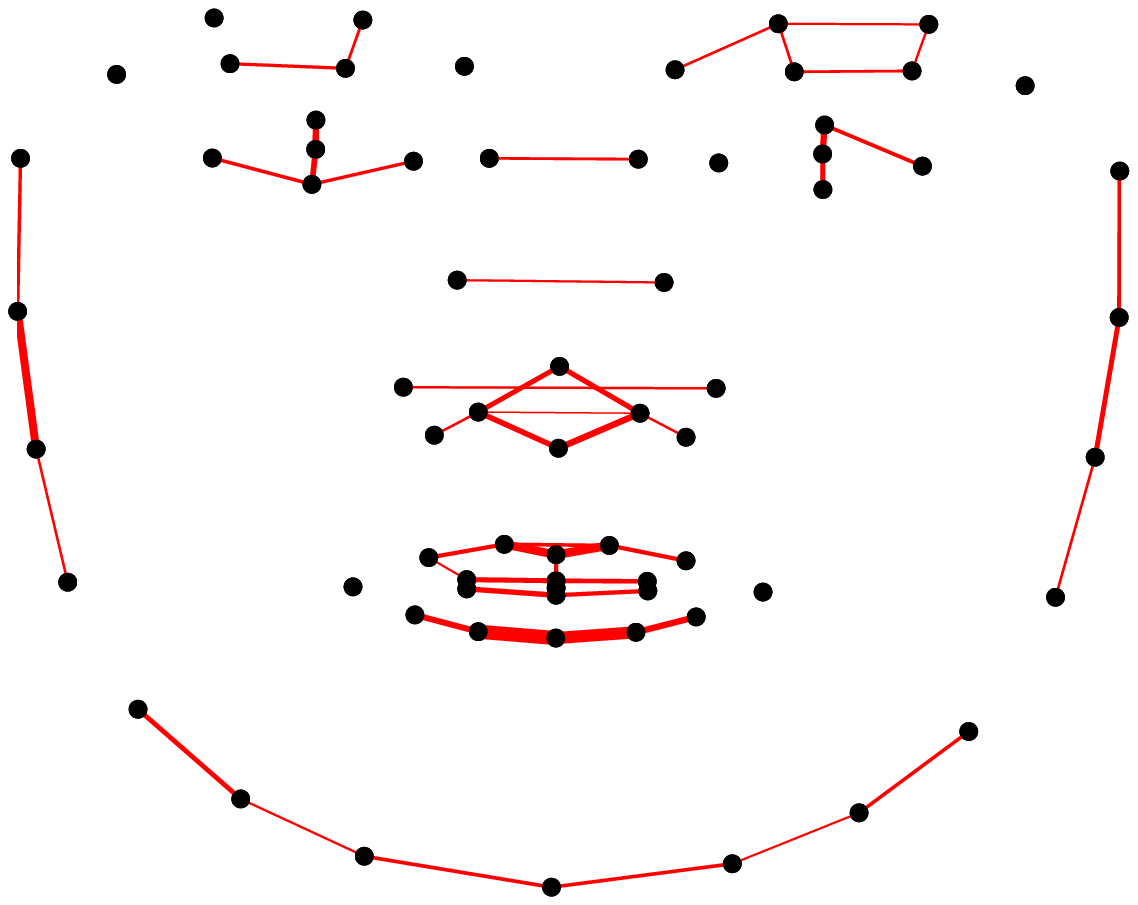}    & \includegraphics[height = 0.15\textheight]{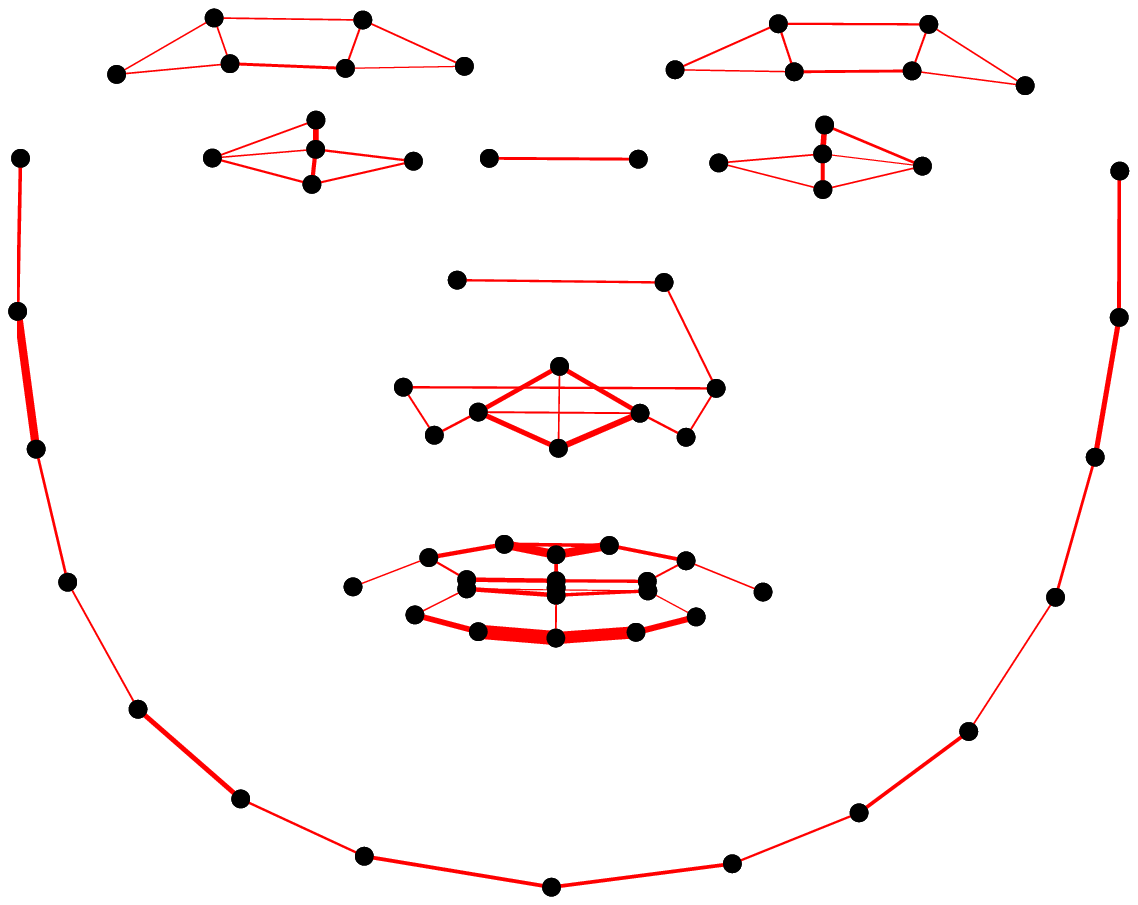}      \\
\includegraphics[height = 0.15\textheight]{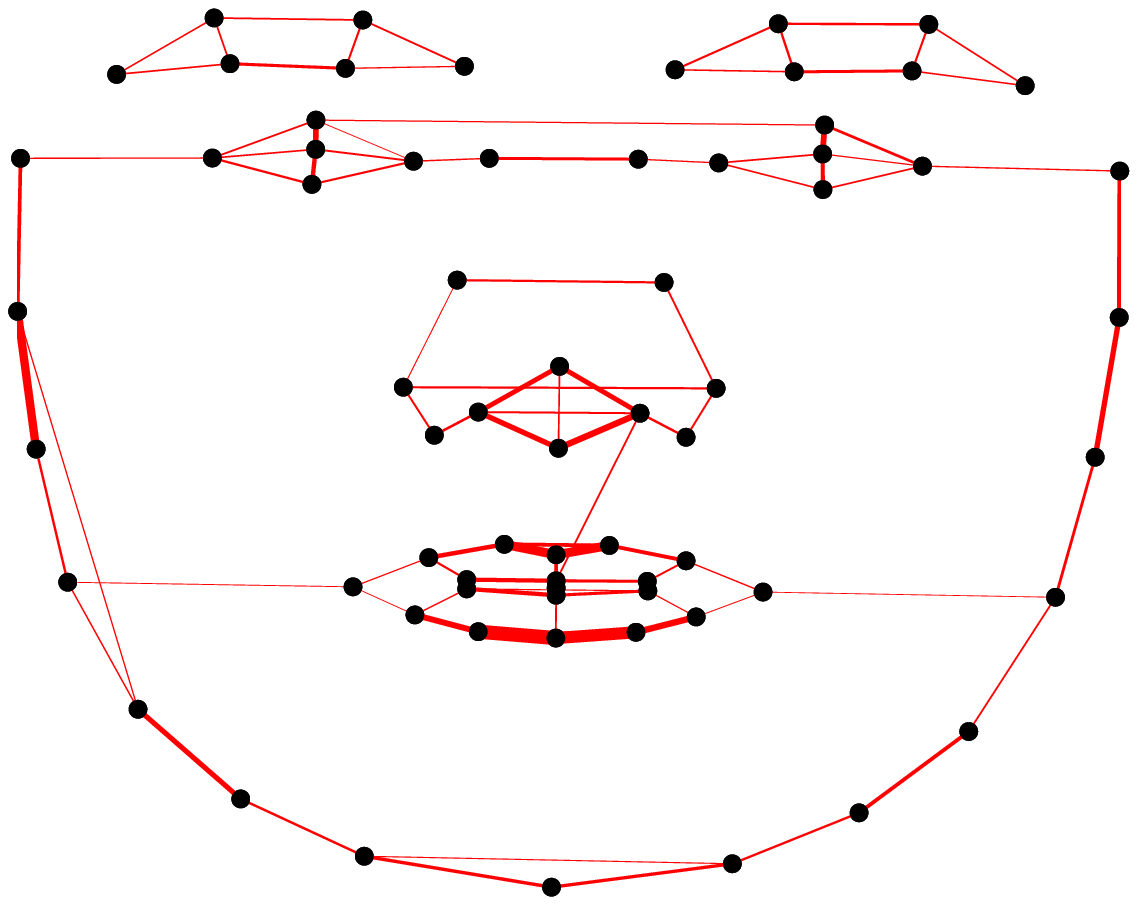}  &  \includegraphics[height = 0.15\textheight]{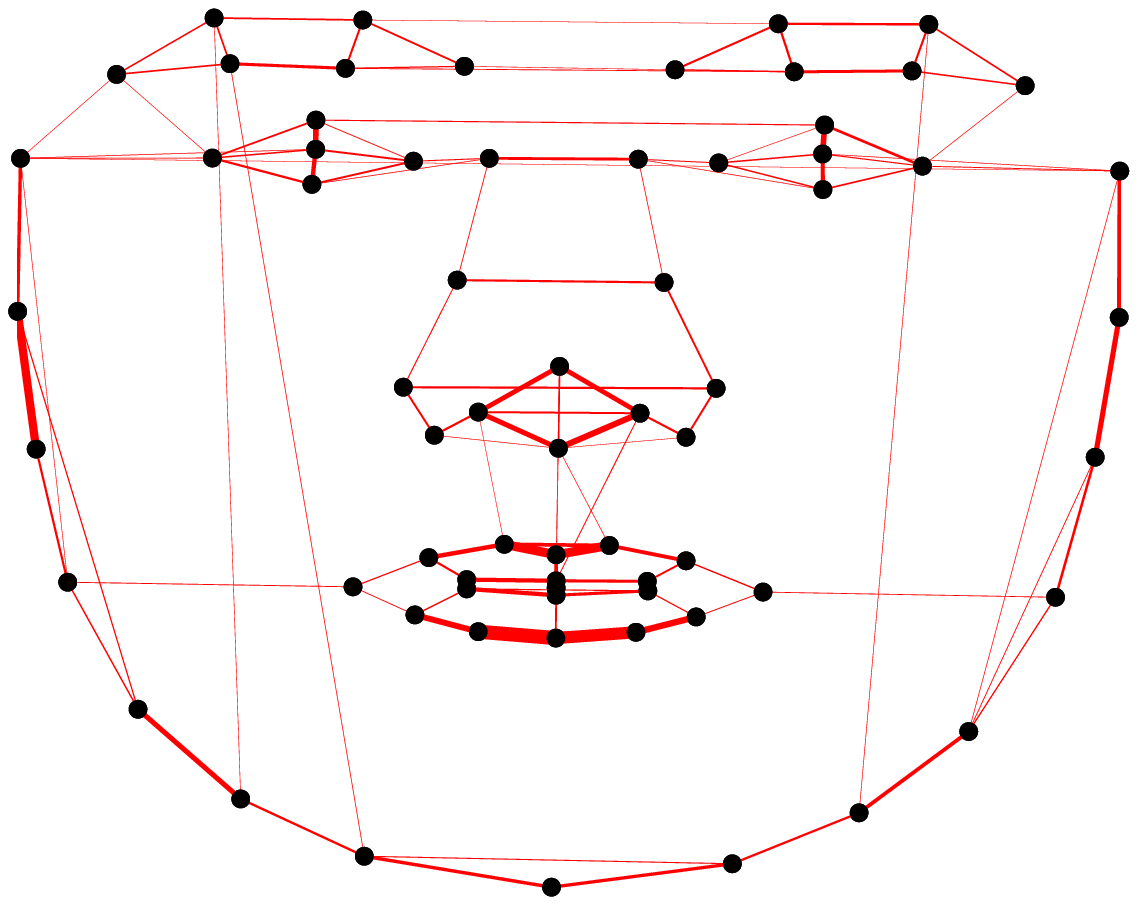}  &
\includegraphics[height = 0.15\textheight]{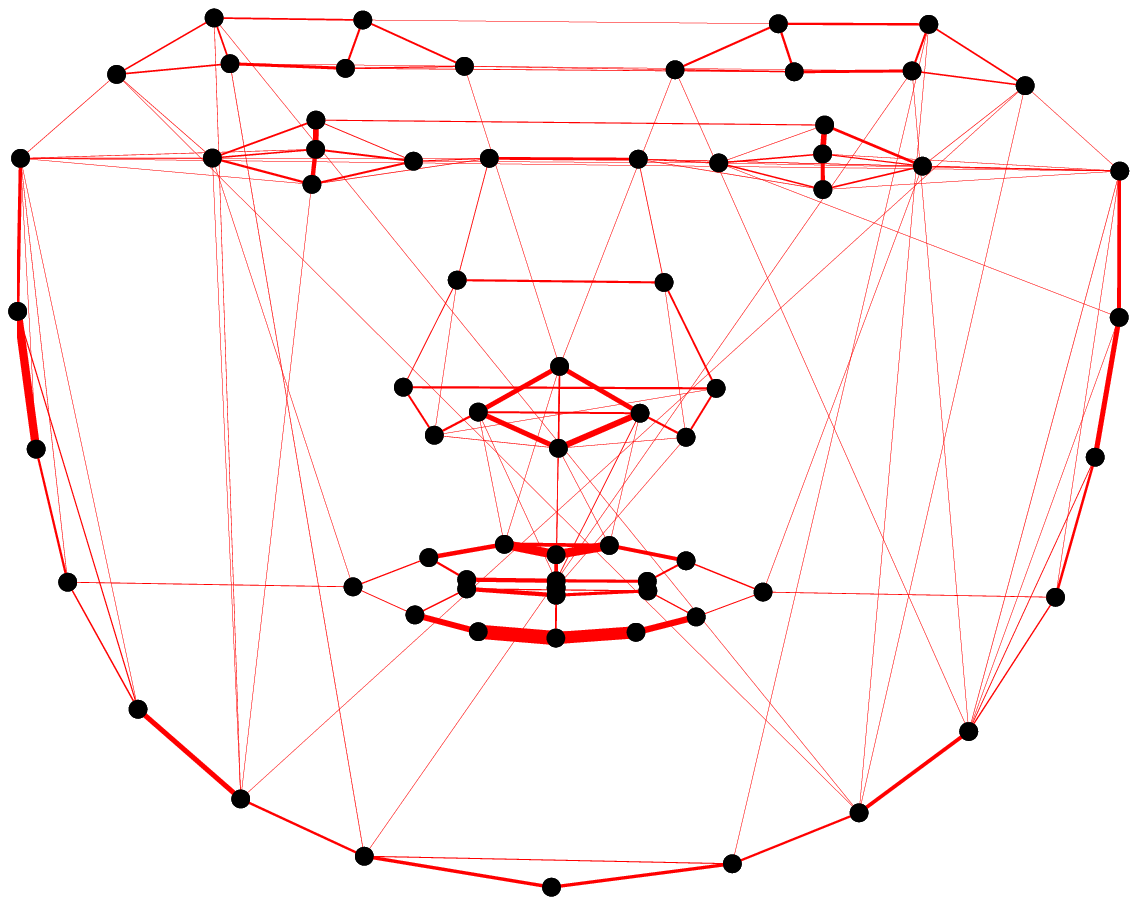}     
\end{tabular}
\rule{0.91\textwidth}{0.5pt}
\begin{tabular}{ccc}
\includegraphics[height = 0.15\textheight]{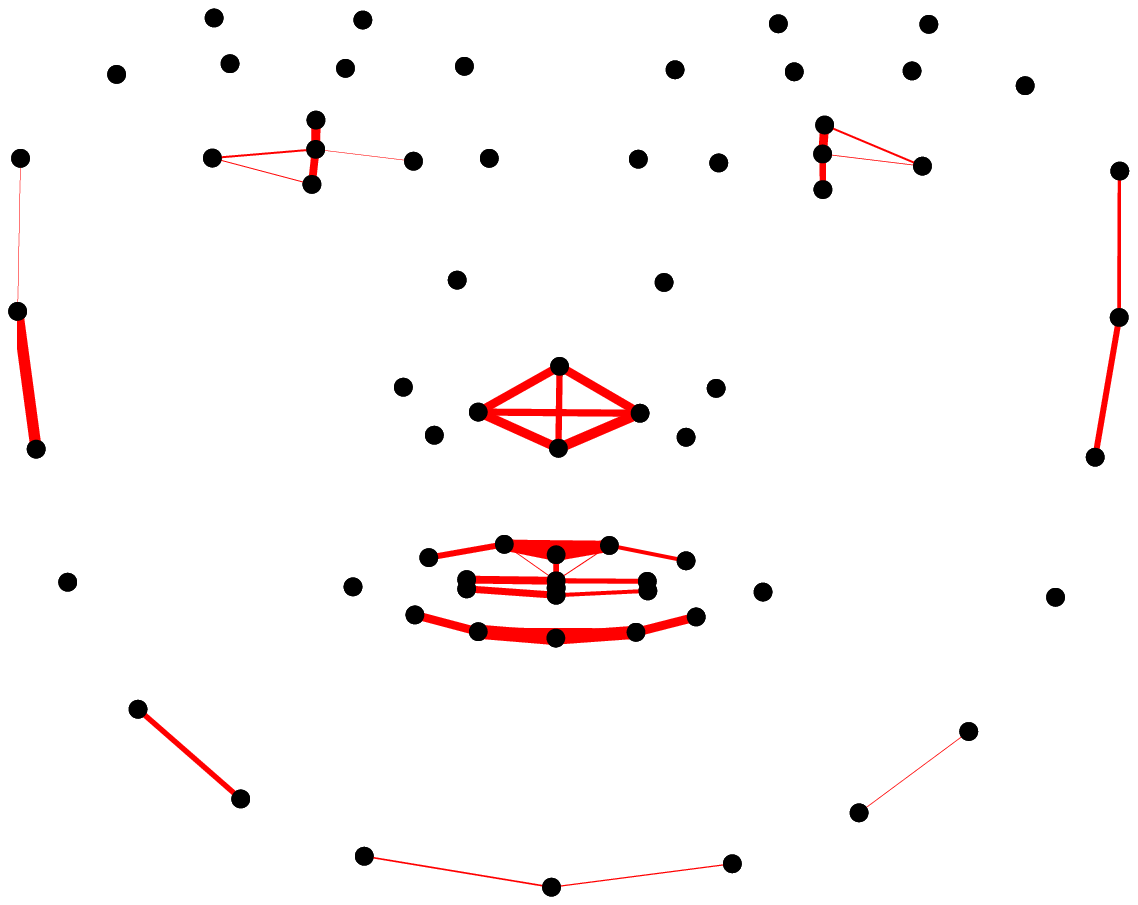} &  \includegraphics[height = 0.15\textheight]{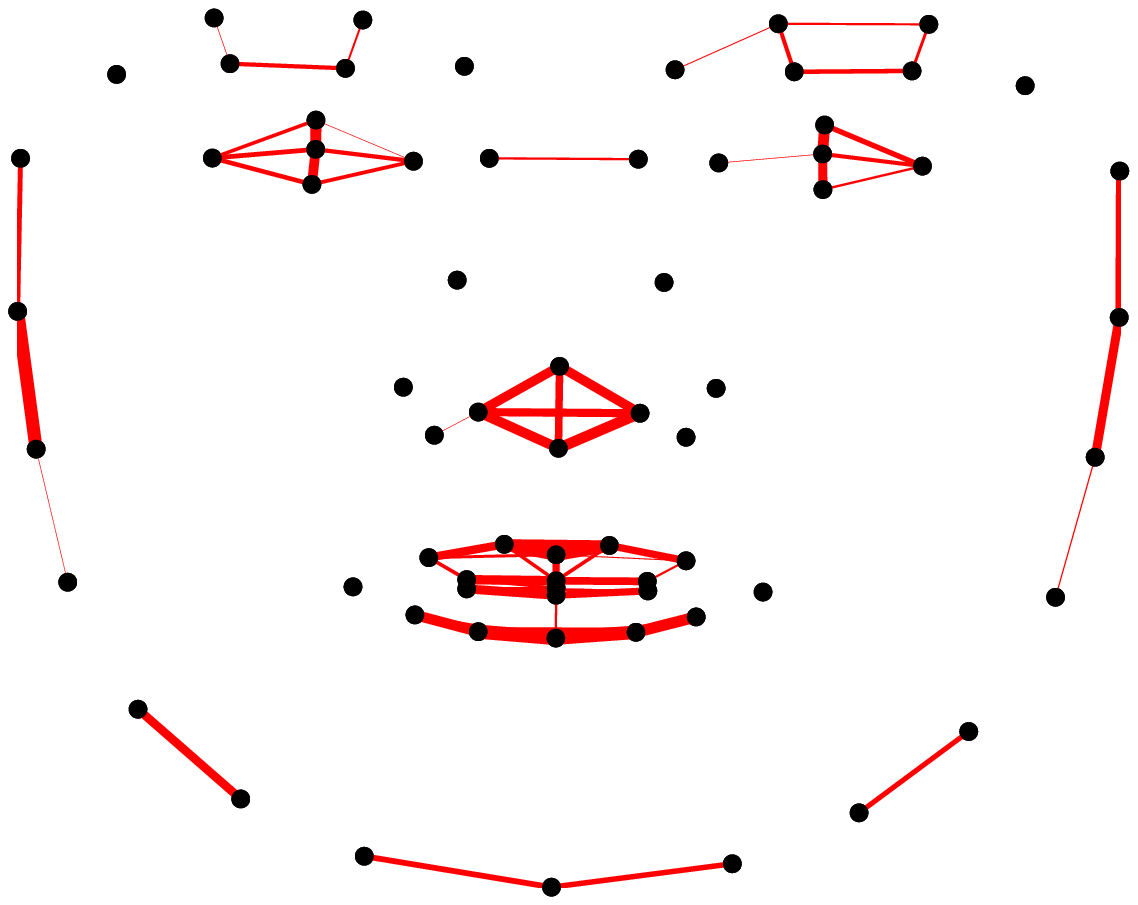}    & \includegraphics[height = 0.15\textheight]{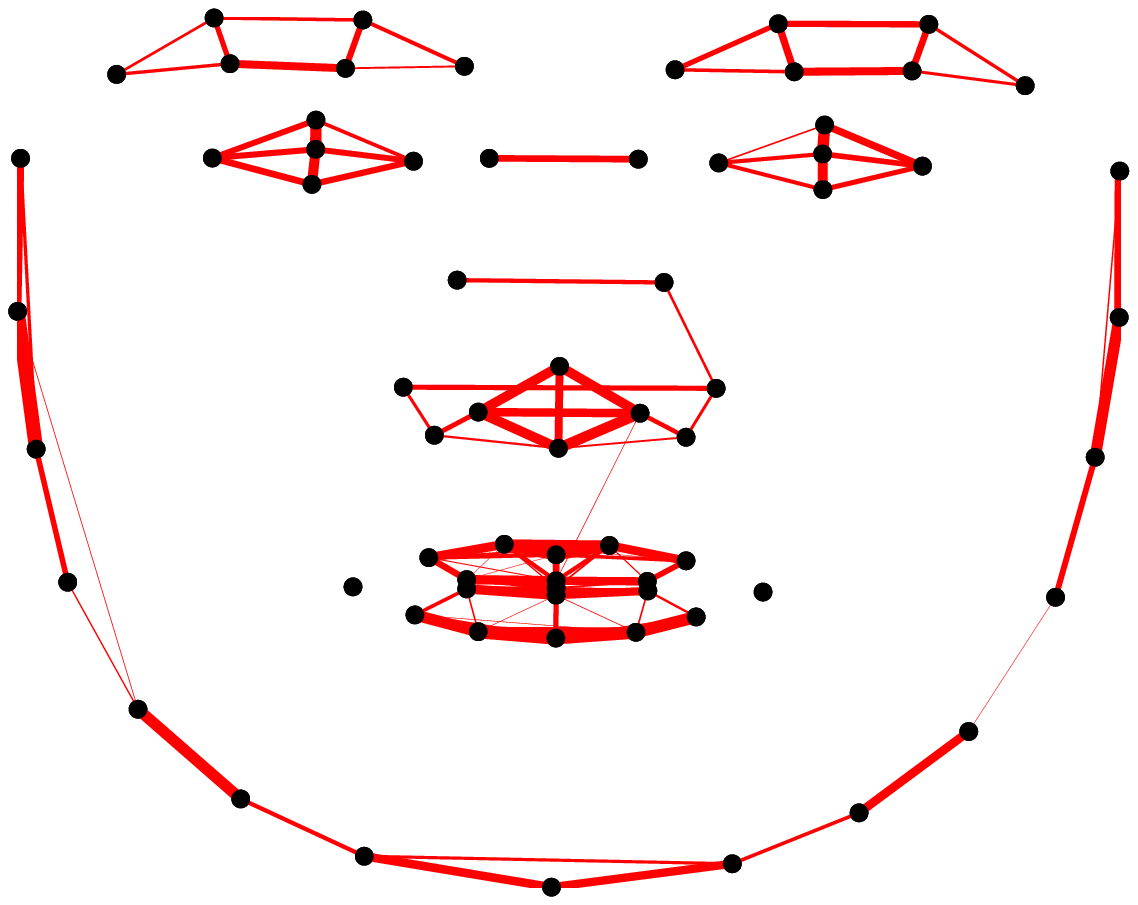}      \\
\includegraphics[height = 0.15\textheight]{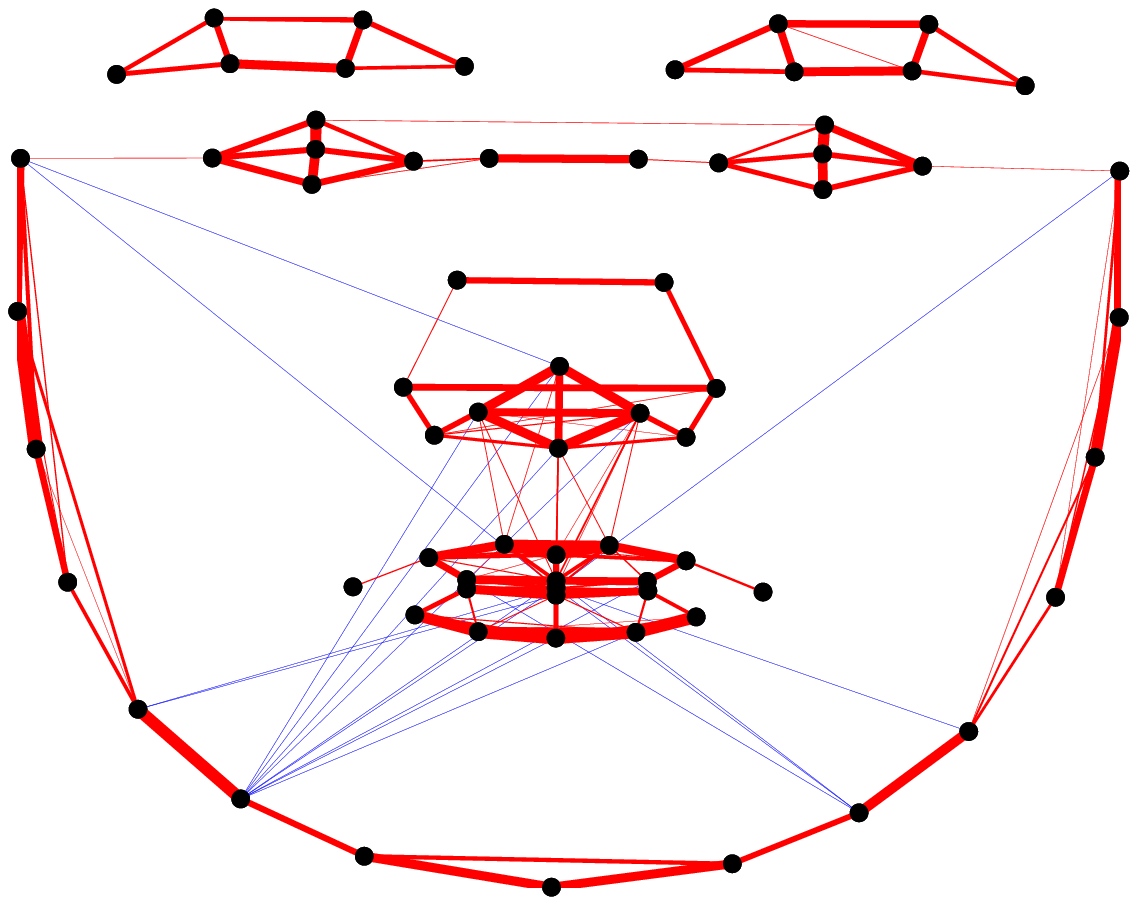}  &  \includegraphics[height = 0.15\textheight]{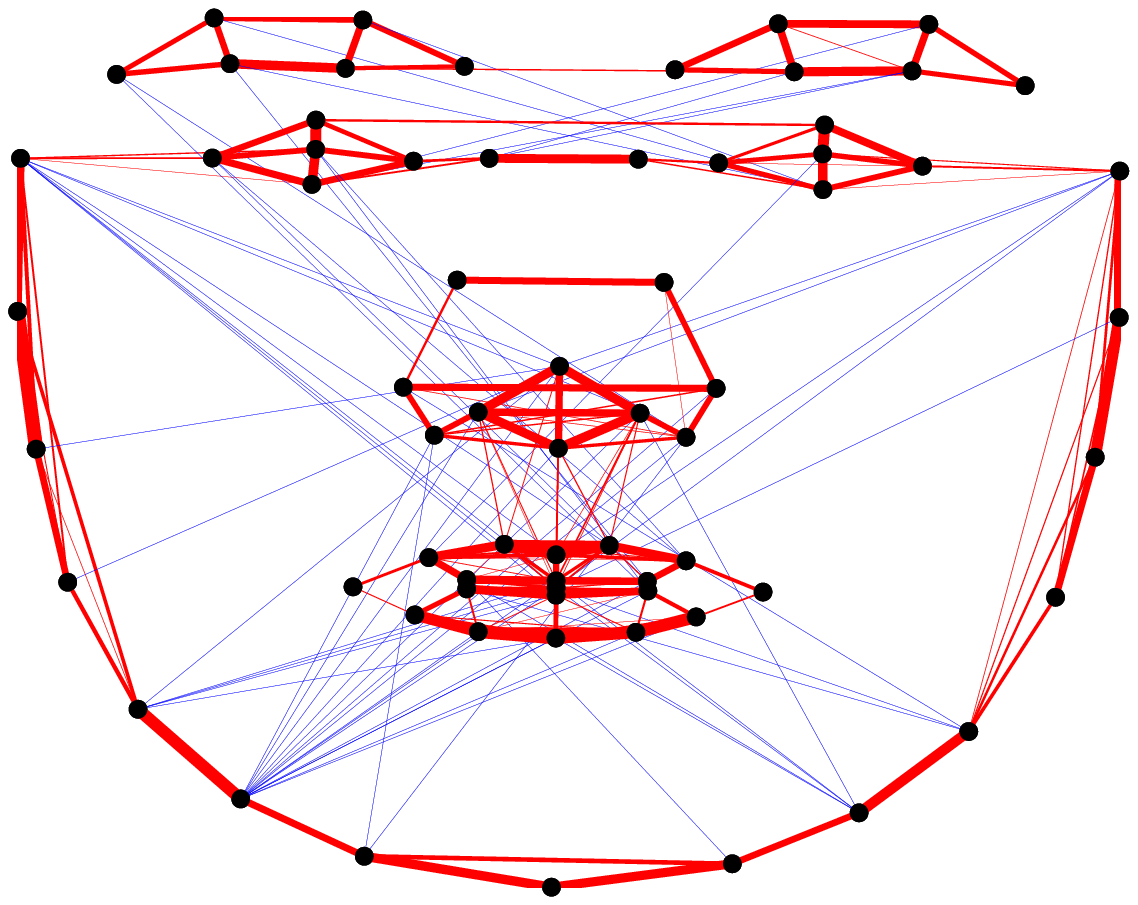}  &
\includegraphics[height = 0.15\textheight]{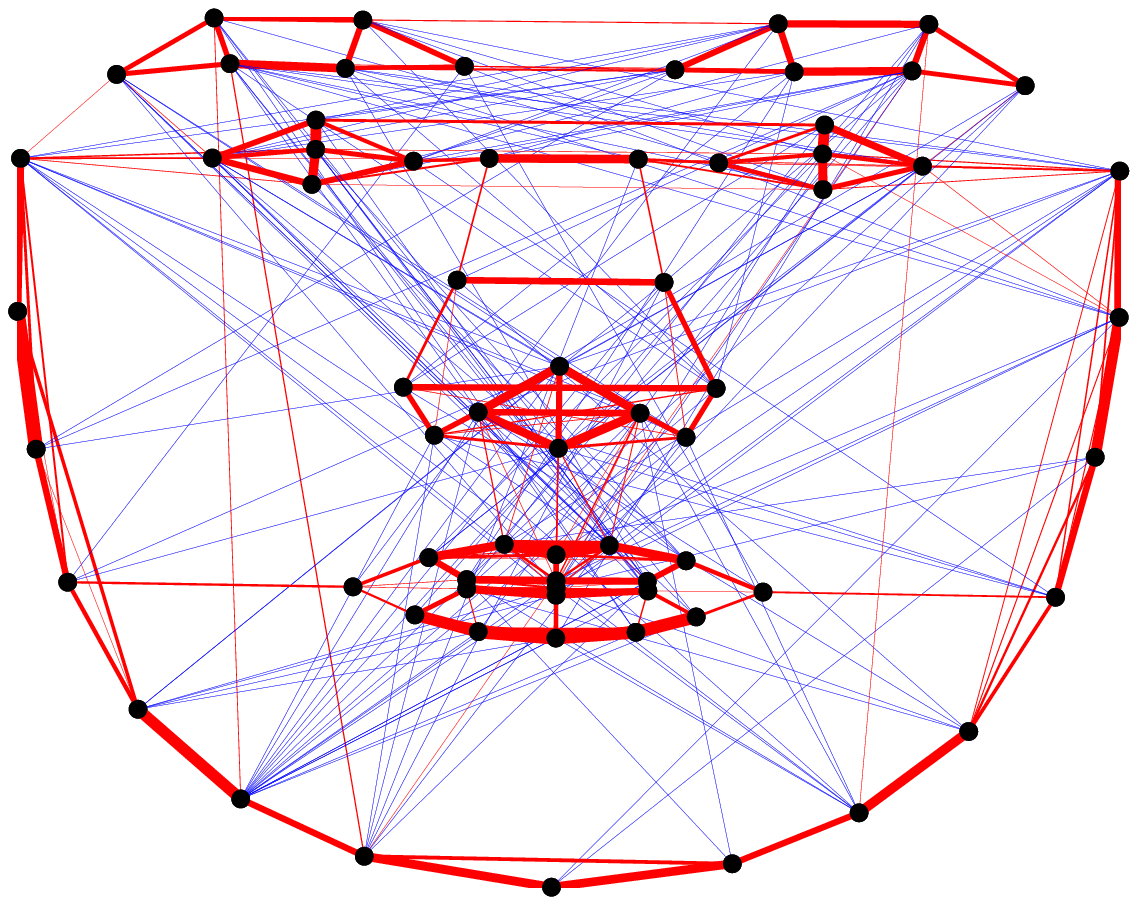}     
\end{tabular}
\caption{Graphs for the landmark dataset for six different levels of sparsity. Top: attractive model. Bottom: graphical lasso. Edge widths
are proportional to the absolute values of the corresponding entries in the estimated precision matrix, and the edge colour represents
their signs (red: negative sign, blue: positive sign).}\label{fig:face_graphs}
\end{figure}

In the present section, we consider the use of positive definite $M$-matrices for modeling
the interdependence of planar landmarks contained in the face image dataset \textsf{XM2VTS} available
from \cite{Cootes}. The concept of our analysis has been inspired by a similar one in \cite{Gu07}, where, however, different
datasets and methods are employed. The full dataset comprises frontal photos of 295 individuals collected over four
sessions, with two shots per individual in each session. A set of $p = 68$ landmarks is collected for each
photo (cf.~left panel of Figure \ref{fig:landmarks}). We here restrict ourselves to the first session, using the firsts of the two shots as the training
set and the second ones as test set. Landmark data falls into the domain of statistical shape analysis
\cite{DrydenMardia2002}, and we apply the usual steps developed in this area (centering, scaling, and generalized
Procrustes analysis) to process the raw data. The thus pre-processed landmark data is given by
$(x_{ij}, y_{ij})$, $i=1,\ldots,n$, $j=1,\ldots,p$, where $(x_{ij}, y_{ij})$ is the coordinate pair of the 
$j$-th landmark of the $i$-th observation (cf.~middle panel of Figure \ref{fig:landmarks}). As starting
point, we consider a factored covariance model of the form $\Sigma_{*} = \Sigma_{*}^{\textsf{xy}} \otimes \Sigma_*^{\textsf{lm}}$,
where $\Sigma_{*}^{\textsf{xy}} \in \psd^2$ represents the covariance between x and y coordinates and $\Sigma_*^{\textsf{lm}} \in \psd^{68}$ 
models the covariance between landmarks (cf.\cite{DrydenMardia2002}, p.167). While this model is restrictive,
since it requires the structure of variability to be the same at each landmark, it is convenient 
for interpretation at the level of landmarks. Instead of fitting such a model directly from the
joint sample covariance $S = \begin{pmatrix} S_{xx} & S_{xy} \\ S_{yx} & S_{yy} \end{pmatrix}$, we consider
a simpler approach that is sufficient for our purpose here. We first determine
\begin{equation}\label{eq:factored_frob}
(S^{\textsf{xy}}, S^{\textsf{lm}}) = \argmin_{\Sigma^{\textsf{xy}} \gec 0, \, \Sigma^{\textsf{lm}} \gec 0} \norm{S - \Sigma^{\textsf{xy}} \otimes \Sigma^{\textsf{lm}}}_F 
\end{equation}
and work subsequently only with $S^{\textsf{lm}}$ in the sequel as input for precision matrix estimation methods. 
Omitting further details here, solving \eqref{eq:factored_frob} does not constitute 
an obstacle, as it turns out that this can be done essentially in closed from. 
The mode of evaluation is as in the previous subsection, but we report the loss on the test set
instead of a cross-validated loss. In Figure \ref{fig:face_graphs}, we compare the graphical lasso
and our approach in terms of the graph structure they provide for several levels of sparsity. Interestingly, 
the graphs are rather similar as long as both of them contain only few edges connecting landmarks belonging
to the distinct parts of the face (eyes, eyebrows, nose, mouth) and a chain along the boundary of the face. As
the graphs become denser, edges appear between the different parts. For the graphical lasso, a considerable fraction
of these edges correspond to negative partial correlations.       
\begin{figure}[h!]
\hspace{-0.8cm}
\begin{minipage}[t]{0.75\textwidth}
\includegraphics[height = 0.31\textheight]{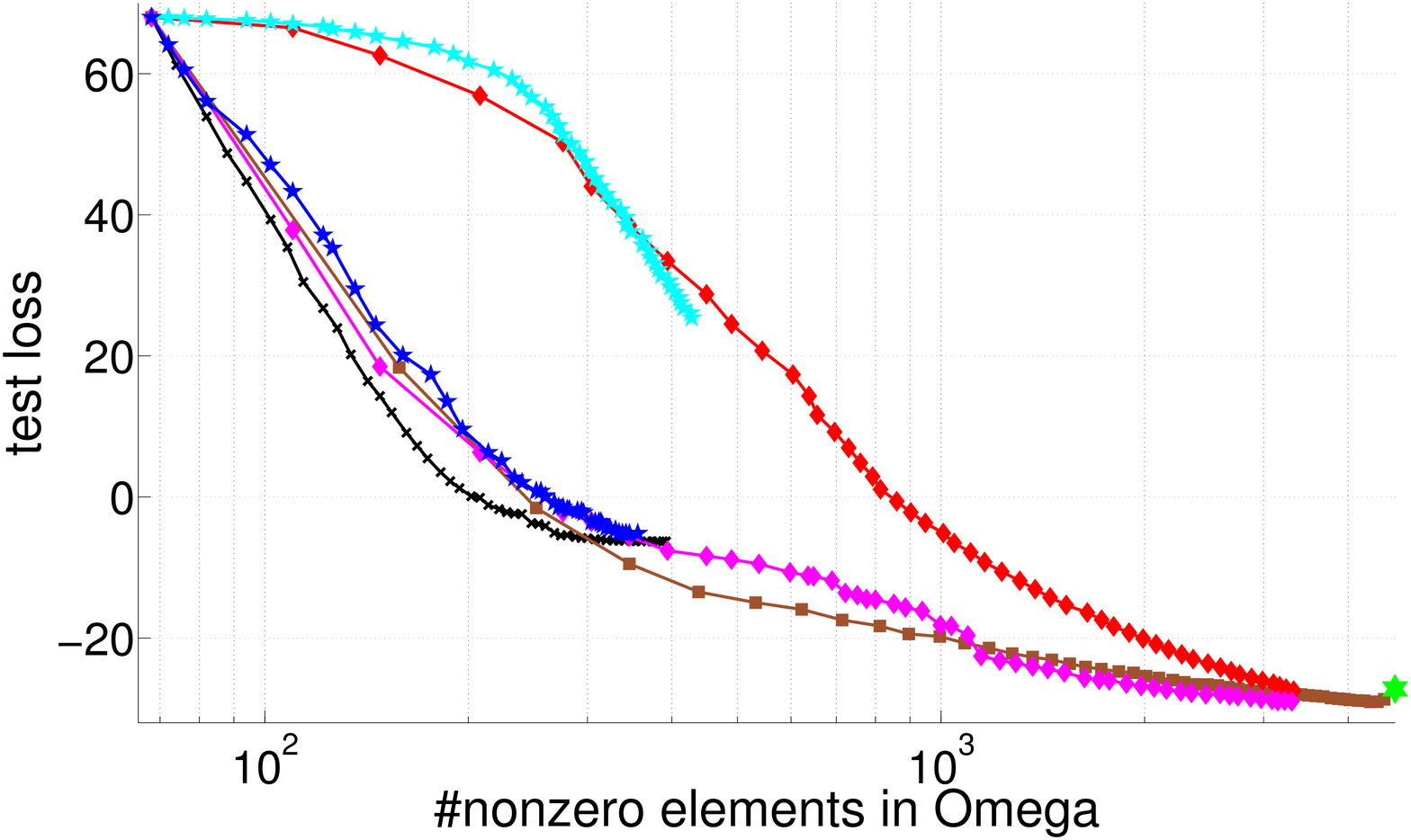}
\end{minipage}
\hspace{-0.5cm}
\begin{minipage}{0.23\textwidth}
\vspace{-6.5cm}
\includegraphics[height = 0.18\textheight]{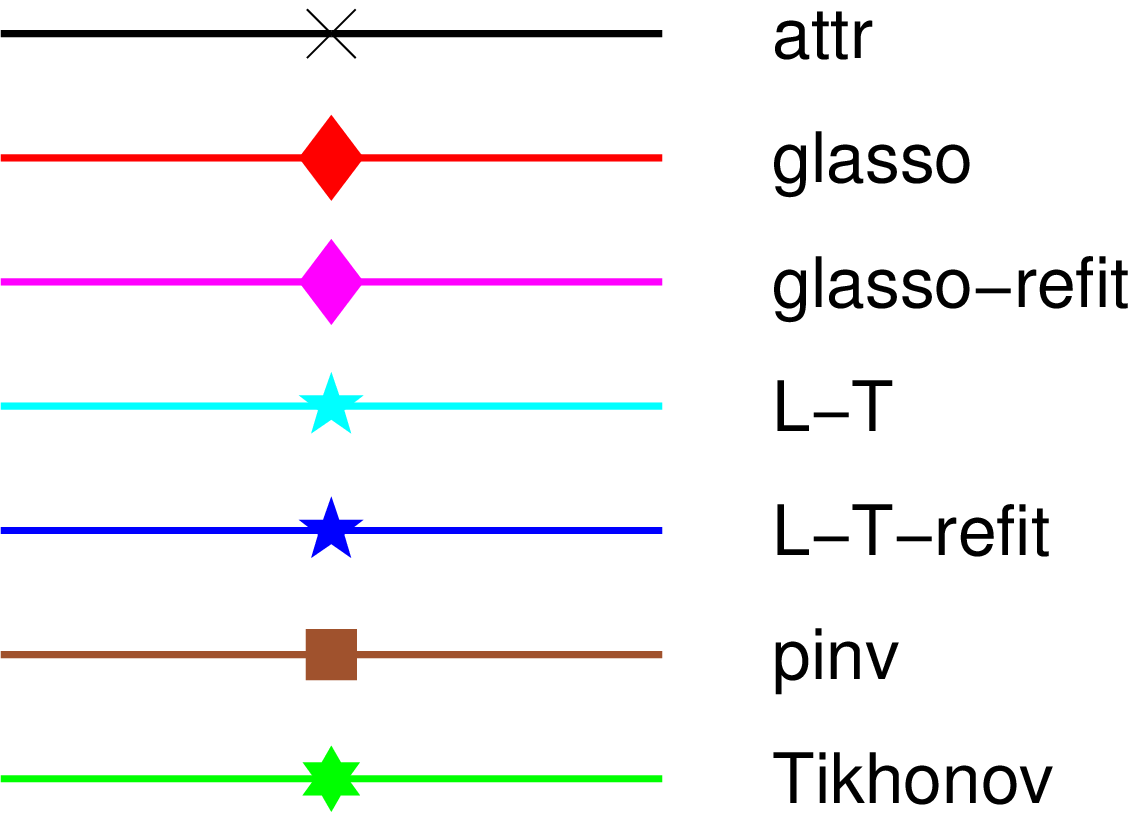}
\end{minipage}
\vspace{-.25cm}
\caption{Test loss on the landmark dataset in dependency of the sparsity of the estimates. The annotation is as for Figure \ref{fig:cv_animals} above.}\label{fig:face_test}
\end{figure}

 \hfill \\

\section{Summary}\label{sec:summary}
The paper has addressed the use of non-positivity constraints on the off-diagonal part in precision
matrix estimation. We have pointed out that the sign constraints constitute a severe restriction
from the perspective of modelling, but also a blessing from the perspective of estimation. Specifically,
we have provided evidence indicating that regularization is no longer compulsory in a high-dimensional regime, and that a simple
post-processing (thresholding) gives rise to well-interpretable results on real world data sets. The theoretical results
of this paper fall behind the empirical findings. Proving structural consistency is left for future research, as is
the development of a faster algorithm that scales more favourably with respect to the number of variables.

\subsection*{Acknowledgments.} We would like to thank Brendan Lake for providing 
us the datasets analyzed in Section \ref{sec:experiments}. 

\appendix

\section{Proof of Theorem \ref{theo:existenceanduniqueness}}\label{proof:existenceanduniqueness}

The proof relies on a characterization of the spectrum of a positive definite $M$-matrix, which
is based on the following results extracted from \cite{Berman2003} and \cite{BermanPlemmons}. In the sequel,
we write $r(A) = \max\{\lambda_1(A), -\lambda_p(A) \}$ for the spectral radius of a symmetric matrix $A \in \R^{p \times p}$
with spectrum $\lambda_1(A) \geq \ldots \geq \lambda_p(A)$. 
\renewcommand{\thesection}{\Alph{section}}

\begin{lemmaApp}\label{lem:characterizationMmatrix} $\Omega \in \overline{\mc{M}^p}$ (resp.~$\Omega \in \mc{M}^p$) if and only if
there exists $B \geq 0$, $B = B^{\T}$ and $\delta \in [r(B), \infty)$ (resp.~$\delta \in (r(B), \infty)$) such that $\Omega = \delta I - B$.
\end{lemmaApp}

\begin{defnApp}  
A matrix $A$ is said to be irreducible if there exists a permutation
matrix $\Pi$ such that 
\begin{equation*}
\Pi A \Pi^{\T} = \left( \begin{array}{cc}
    A_1 & A_2 \\
    0 & A_3 \\
\end{array} \right). 
\end{equation*}
\end{defnApp}

\begin{remApp}\label{rem:reducibility}
Note that if $A$ is symmetric, $\Pi A \Pi^{\T}$ is symmetric as well.
Consequently, if $A$ is symmetric and reducible, it must be a block diagonal
matrix.
\end{remApp}  
 
\begin{theoremApp}\label{theo:perronfrobenius}(Perron$-$Frobenius) Let $A \geq 0$ be a square irreducible matrix. Then
$r(A)$ is a positive, algebraically simple eigenvalue of $A$, and $A$ has a corresponding positive eigenvector.
\end{theoremApp}

\begin{remApp} If $A \geq 0$ is symmetric and reducible, then, after a
suitable permutation, $A$ is block diagonal (cf.~Remark \ref{rem:reducibility}) with blocks $A_1,\ldots,A_K$, say, and the Perron-Frobenius Theorem applies to each
block. The positive eigenvectors corresponding to $r(A_1),\ldots,r(A_K)$
have disjoint supports given by the row/column indices of the blocks.
\end{remApp}

The following lemma is an immediate consequence of Lemma \ref{lem:characterizationMmatrix}. 

\begin{lemmaApp}\label{lem:spectrumMmatrix}
Let $\Omega \in \mc{M}^p$. Consider $B$ and $\delta$ according to Lemma 
\ref{lem:characterizationMmatrix} and let
\begin{equation*}
\gamma \coloneq r(B),  \quad \sigma_j \coloneq \lambda_j(B/\gamma), \; \, j=1,\ldots,p, \; \; \; \, \eps \coloneq \delta - \gamma. 
\end{equation*}
Then, $v \in \R^p$ is an eigenvector of $\Omega$ if and only if it is an eigenvector of $B$, and the spectrum of $\Omega$ is given by
$\lambda_j(\Omega) = \eps + \gamma (1 - \sigma_{p - j + 1}) > 0$, $j=1,\ldots,p$.
\end{lemmaApp}

\paragraph{Proof of Theorem \ref{theo:existenceanduniqueness}}          

It will be shown that under the stated conditions,
there exists $R < \infty$ such that 
\begin{equation}\label{eq:min_spectralball}
\min_{\Omega \in \overline{\mc{M}^p}} -\log \det(\Omega) + \tr(\Omega S) = \min_{\Omega \in \overline{\mc{M}^p}:\,\norm{\Omega} \leq R} -\log \det(\Omega) + \tr(\Omega S), 
\end{equation}
which implies existence of a minimizer. Uniqueness readily follows from existence in view of the strict convexity of the negativ log-determinant on 
$\overline{\psd^p}$. To establish \eqref{eq:min_spectralball}, it suffices to show that the objective is bounded
from below. Let $\Omega \in \mc{M}^p$ be arbirtrary. Expanding $\Omega = \sum_{j = 1}^p \lambda_j u_j u_j^{\T}$, 
where $\lambda_j \coloneq \lambda_j(\Omega), \, j=1,\ldots,p$, and $\{ u_j \}_{j=1}^p$ are the corresponding eigenvectors, the 
objective evaluated at $\Omega$ can be written as
\begin{align*}
-\log \det(\Omega) + \tr(\Omega S) &= -\sum_{j = 1}^p \log(\lambda_j) + \sum_{j = 1}^p \lambda_j u_j^{\T} S u_{j} \\
&=-\sum_{j = 1}^p \log \left\{ \eps +  \gamma (1 - \sigma_{p - j + 1}) \right\} +\\ 
&\;\;\;\,+\sum_{j = 1}^p\left\{ \eps +  \gamma (1 - \sigma_{p - j + 1})  \right\} u_j^{\T} S u_{j}, 
\end{align*}
where the second equality follows from $\Omega \in \mc{M}^p$ and 
Lemma \ref{lem:spectrumMmatrix}. We first note that the objective can be unbounded from below only if $\gamma$ is unbounded and $\eps$ is bounded from above, respectively. To show the boundedness of $\eps$, note that 
\begin{align}\label{eq:epsbounded}
&-\sum_{j = 1}^p \log \left\{ \eps +  \gamma (1 - \sigma_{p - j + 1}) \right\} + 
\sum_{j = 1}^p\left\{ \eps +  \gamma (1 - \sigma_{p - j + 1})  \right\} u_j^{\T} S u_{j} \notag \\
&\geq -\sum_{j = 1}^p \log \left\{ \eps +  \gamma (1 - \sigma_{p - j + 1}) \right\}
+ \sum_{j = 1}^p  \eps u_j^{\T} S u_j \\
&\geq  -\sum_{j = 1}^p \log \left\{ \eps +  \gamma (1 - \sigma_{p - j + 1}) \right\} 
+ \eps \, \max_{1 \leq j \leq p}  u_j^{\T} S u_j. \notag
\end{align}
If the objective is unbounded from below, $\eps$ can hence be unbounded only if 
$S u_j = 0$, $j = 1,\ldots,p$, because the first term decreases
logarithmically in $\eps$, whereas the second term increases linearly in $\eps$. The condition $S u_j = 0$, $j = 1,\ldots,p$, however, implies that 
$S$ is the zero matrix, which contradicts the assumption that $S$ has
positive diagonal entries. Using a similar argument as in \eqref{eq:epsbounded} while noting that $|\sigma_k| \leq 1$, $k=1,\ldots,p$, if $\gamma$ is unbounded, the objective can be unbounded from below only if 
\begin{align*}
\bigwedge_{j=1}^p  (S u_{j} = 0 \, \vee \,  \gamma (1 - \sigma_{p - j + 1}) = 0) \; \Longleftrightarrow \, 
\bigwedge_{j=1}^{p-1}  (S u_{j} = 0 \, \vee \,  \gamma (1 - \sigma_{p - j + 1}) = 0),
\end{align*}
where the equivalence follows from $\sigma_1 = \lambda_1(B/\gamma) = 1$ according to 
Lemma \ref{lem:spectrumMmatrix} and the Perron-Frobenius Theorem. Consider now the following cases.\\
\emph{Case 1.} If $1 - \sigma_{p - j + 1} > 0, \, j=1,\ldots,p-1$, we must have $S u_j = 0$, $j = 1,\ldots,p-1$,  
which implies that $S$ has rank one.\\
\emph{Case 2.} If $1 - \sigma_{k} = 0$ for some $k \in \{2,\ldots,p\}$, we have
$\sigma_k = \sigma_1 = r(B)$. Consequently, the eigenvalue $\sigma_1$ of $B$ has multiplicity greater than $1$.
According to Theorem \ref{theo:perronfrobenius}, the symmetric non-negative matrix $B$
must be reducible, as must be $\Omega$. As a result, there exists a 
partitioning $I_1,\ldots,I_K$ ($K \geq 2$) of $\{1,\ldots,p\}$ and
a permutation matrix $\Pi$ so that $\Pi \Omega \Pi^{\T} = \text{bdiag}( \Omega_{11},\ldots,\Omega_{KK})$,
where $\Omega_{kk}$ is the principal submatrix corresponding to
index set $I_k$, $k=1,\ldots,K$. We may assume that the $\{ \Omega_{kk} \}_{k = 1}^K$ are irreducible; otherwise, 
we could simply apply additional permutations to end up with $K' > K$ blocks that are irreducible. The objective can now be decoupled as follows.
\begin{align*}
-\log \det(\Omega) + \tr(\Omega S) &=  -\log \det(\Pi \Omega \Pi^{\T}) + \tr\left\{(\Pi \Omega \Pi^{\T}) (\Pi S \Pi^{\T})\right\} \\
                                             &= \sum_{k = 1}^K -\log \det(\Omega_{kk}) + \tr(\Omega_{kk} S_{kk}),
\end{align*}
where $S_{kk}$, $k=1,\ldots,K$, are the principal submatrices of $\Pi S \Pi^{\T}$ corresponding
to $\Omega_{kk}$, $k=1,\ldots,K$. From the last display, we conclude that the objective
is unbounded from below only if there exists $k \in \{1,\ldots,K\}$ such that 
$-\log \det(\Omega_{kk}) + \tr(\Omega_{kk} S_{kk})$ is unbounded from below. At this 
point, for each $k=1,\ldots,K$, we may recur to the reasoning that has led us to the case distinction above. 
As the $\{ \Omega_{kk} \}_{k = 1}^K$ are irreducible, the second case cannot occur any longer. Consequently,
the objective is unbounded from below only if one of the $\{ S_{kk} \}_{k = 1}^K$ has rank one, i.e.~for
some $k \in \{1,\ldots,K\}$, we have that $S_{kk} = v v^{\T}$ with $v \in \R^{|I_k|}$. Note that if
$|I_k| = 1$, $S_{kk}$ is a scalar, i.e.~a diagonal entry of $S$ which is assumed to be positive, 
so that $-\log \det(\Omega_{kk}) + \tr(\Omega_{kk} S_{kk})$ cannot be unbounded from below.
If $|I_k| \geq 3$, there must exist $\ell,m \in \{1,\ldots,|I_k| \}$, $\ell \neq m$, so that
\begin{equation*} 
0 < (S_{kk})_{\ell m} = v_{\ell} v_{m} = \sqrt{v_{\ell}^2 \, v_{m}^2} =  \sqrt{(S_{kk})_{\ell \ell} (S_{kk})_{mm}}.
\end{equation*}
This can be seen by enumerating all possible sign patterns of $v$ for $|I_k| = 3$, i.e.~schematically $(+,+,+),(-,+,+),\ldots,(-,-,+),(-,-,-)$ (note that all entries of $v$ must be non-zero, since the diagonal entries of $S$ are positive by assumption) and verifying that at least
one of $v_1 \cdot v_2$, $v_1 \cdot v_3$, $v_2 \cdot v_3$ must be positive. To finish the proof of the theorem, it remains to consider the case $|I_k| = 2$, where
\begin{equation*}
S_{kk} = \begin{pmatrix}
        v_{1}^2 & v_{1} v_{2} \\
        v_{1} v_{2} & v_{2}^2 
        \end{pmatrix}. 
\end{equation*}
Denote the eigenvectors of $\Omega_{kk} = \delta_{kk} I - B_{kk}$, say, by $z_1$ and $z_2$, $z_1^{\T} z_2 = 0$. Invoking
Theorem \ref{theo:perronfrobenius}, we may choose $z_1$ such that its two entries are of the same sign
and those of $z_2$ have a different sign. Note that since Case 1 (see above) occurs for $S_{kk}$,
the eigenvector $z_2$ not corresponding to the largest eigenvalue of $B$ must
satisfy $S_{kk} z_2 = 0$. This can hold only if $v_1 v_2 > 0$. $\qed$

\renewcommand{\thesection}{Appendix \Alph{section}}

\section{Proof of Theorem \ref{theo:characterization}}\label{proof:characterization}

Let $\Omega = \Sigma^{-1}$ and consider the Bregman divergence between $\Omega$ and $\Omega' \in \overline{\psd^p}$ which is induced by the
log-determinant  
\begin{equation*}
D(\Omega \parallel \Omega') = - \log \det(\Omega') + \log \det(\Omega) + \tr(\Sigma(\Omega' - \Omega)).
\end{equation*}
In virtue of properties of Bregman divergences, we have
\begin{equation}\label{eq:Bregmanproperty}
D(\Omega \parallel \Omega') \geq 0, \; \text{with equality holding if and only if} \;\, \Omega' = \Omega. 
\end{equation}
If $\Sigma \in  \psd^p$ is an inverse $M$-matrix, then
$\Omega \in \mc{M}^p$. By \eqref{eq:Bregmanproperty} and the definition of $D$  
\begin{equation*}
\Omega \in \mc{M}^p  \, \Longrightarrow \, \Omega = \argmin_{\Omega' \in \overline{\mc{M}^p}} D(\Omega \parallel \Omega') 
\, \Longleftrightarrow \,  \Omega = \argmin_{\Omega' \in \overline{\mc{M}^p}} -\log \det(\Omega') + \tr(\Sigma \Omega'). 
\end{equation*}
From the duality relation \eqref{eq:dualproblem}, the last property implies that
\begin{equation*}
\Sigma = \argmax_{\Sigma' \in \overline{\psd}, \; \; \Sigma' \geq \Sigma, \;
  \diag(\Sigma') = \diag(\Sigma)} \log \det(\Sigma') + p 
\end{equation*}
For the opposite direction, suppose that $\Omega \in \psd^p \setminus \mc{M}^p$ and denote \newline $\Omega_{\bullet} = \argmin_{\Omega' \in \overline{\mc{M}^p}} D(\Omega \parallel \Omega')$. 
From \eqref{eq:Bregmanproperty}, we have $D(\Omega \parallel \Omega_{\bullet}) > 0$ and hence also
\begin{equation*}
-\log \det(\Omega_{\bullet}) + \tr(\Sigma \Omega_{\bullet}) > -\log \det(\Omega) + \tr(\Sigma \Omega) = \log \det(\Sigma) + p. 
\end{equation*}
Denote $\Sigma_{\bullet} = \Omega_{\bullet}^{-1}$. By definition of $\Omega_{\bullet}$ and \eqref{eq:dualproblem},
\begin{equation*}
\log \det(\Sigma_{\bullet}) + p = -\log \det(\Omega_{\bullet}) + \tr(\Sigma \Omega_{\bullet}) > -\log \det(\Omega) + \tr(\Sigma \Omega) = \log \det(\Sigma) + p.
\end{equation*}
It follows that $\log \det(\Sigma_{\bullet}) > \log \det(\Sigma)$ with $\Sigma_{\bullet} \geq \Sigma$ and $\text{diag}(\Sigma_{\bullet}) = \text{diag}(\Sigma)$.

\section{Proofs for the Examples in Section \ref{sec:misspecification}}

\paragraph{(3) AR(2)-structure} First note that under the condition $4|\phi_2| < \phi_1$, it holds that $\phi_1 > 0$. 
Hence if $\phi_2$ is also non-negative, the parameters of the corresponding AR(2) process and thus
all partial correlations are non-negative \cite{Rue2001} so that $\Omega_* \in \mc{M}^p$, in which case
property \eqref{eq:misspec_signpreserve} is trivially satisfied. On the other hand, if $\phi_2 < 0$ while $4|\phi_2| < \phi_1$,
the sequence $(\rho_{\ell})$ as given in \eqref{eq:yulewalkerAR2} is monotonically decreasing (cf.~\cite{Chatfield2003}, p.45).
This property will allow us to show that $\Omega_{\bullet}$  \eqref{eq:logdet_population_Omega} is an $M$-matrix corresponding to an AR(1)-structure. To this end, it is established
that the inverse $\Sigma_{\bullet}$ has entries $\sigma_{jk}^{\bullet} =
\rho_{|j-k|}^{\bullet}$, $j,k=1,\ldots,p$, where with $\rho_1$ as in \eqref{eq:yulewalkerAR2},  
\begin{equation*}
\rho_{\ell}^{\bullet} = \rho_1^{\ell} = \left(\frac{\phi_1}{1 - \phi_2}\right)^{\ell}, \quad \ell=1,\ldots,p-1. 
\end{equation*}
In order to verify the optimality conditions of \eqref{eq:logdet_population_Omega}/\eqref{eq:logdet_population_Sigma} according to \eqref{eq:kkt}, it
suffices to show complementarity slackness \eqref{eq:compslackness_componentwise}, that is $\rho_1^{\bullet} = \rho_1$ and further $\rho_{\ell}^{\bullet} >
\rho_{\ell}$ for $\ell=2,\ldots,p-1$. The first claim is immediate from
\eqref{eq:yulewalkerAR2}. The second claim is proved by induction. The base case ($\ell = 2$) follows from
\begin{align*} 
\rho_2 &= \phi_1 \rho_1 + \phi_2 \\
       &=\rho_1 (1 - \phi_2) \rho_1 + \phi_2 = \rho_1^2 + \phi_2 (1 -
       \rho_1^2) < \rho_1^2 = \rho_{2}^{\bullet},
\end{align*}
since $\phi_2 < 0$. Considering $\ell \geq 3$, we have
\begin{align*}
\rho_{\ell} &= \phi_1 \rho_{\ell-1} + \phi_2 \rho_{\ell-2} \\
           &= \rho_1 (1 - \phi_2) \rho_{\ell-1} + \phi_2 \rho_{\ell-2} \\
           &= \rho_1 \rho_{\ell-1} + \phi_2 (\rho_{\ell-2} - \rho_1
           \rho_{\ell-1})
\end{align*}
The second term is negative, because the term inside the
brackets is positive, as $\rho_{\ell-2} > \rho_{\ell-1}$, recalling that $(\rho_{\ell})$ is monotonically decreasing. By the induction
hypothesis, $\rho_{\ell-1} < \rho_{\ell-1}^{\bullet} = \rho_{1}^{\ell-1}$ so that 
$\rho_1 \rho_{\ell-1} < \rho_1^{\ell}$ and consequently $\rho_{\ell} <
\rho_1^{\ell} = \rho_{\ell}^{\bullet}$, $\ell=2,\ldots,p-1$, as claimed.$\qed$
\paragraph{(4) Star structure} We verify the expressions for $\Omega_{\bullet}$ and $\Sigma_{\bullet}$ as given in the main text. 
This is done by checking the KKT optimality conditions of \eqref{eq:logdet_population_Omega}, cf.\eqref{eq:kkt}. We set
\begin{equation*}
\Gamma_{\bullet} = \begin{pmatrix}
          0 & \delta^{\T} \\
          \delta & - \rho \wt{\rho}^{\T} -
              \wt{\rho} \rho^{\T} + 2 \wt{\rho} \wt{\rho}^{\T} 
          \end{pmatrix},
\end{equation*}
with $\rho$, $\wt{\rho}$ as defined in the main text and $\delta = \rho - \wt{\rho}$, so that $\Sigma_{\bullet} = \Sigma_{*} + \Gamma_{\bullet}$. Because of
dual feasibility and complementarity slackness, the following has to hold for the entries ($\gamma_{jk}^{\bullet}$) of $\Gamma_{\bullet}$: 
\begin{align*}
(1):\,  &\quad \gamma_{jj}^{\bullet} = 0, \; j=1,\ldots,p,\\
(2):\,  &\quad \gamma_{j1}^{\bullet} = \gamma_{1j}^{\bullet} = 0 \,\; \text{for all} \; j \neq 1 \; \text{such that} \; \omega_{j1}^{\bullet} < 0,\\
(3):\,  &\quad \gamma_{j1}^{\bullet} = \gamma_{1j}^{\bullet} > 0 \,\; \text{for all} \; j \neq 1 \; \text{such that} \; \omega_{j1}^{\bullet} = 0,\\
(4):\,  &\quad \gamma_{jk}^{\bullet} = 0 \,\; \text{for all} \; (j,k), j,k \neq 1,  \; \text{and} \; \omega_{jk}^{\bullet} < 0,\\
(5):\,  &\quad \gamma_{jk}^{\bullet} > 0 \,\; \text{for all} \; (j,k), j,k \neq 1,  \; \text{and} \; \omega_{jk}^{\bullet} = 0.
\end{align*}
The requirements (2) and (3) follow immediately from the definition of $\delta$ and $\wt{\rho}$. For the remainder,  
let $j,k > 1$ be arbitrary and set $\ell = j - 1$, $m = k-1$. We have that
\begin{equation*}
\gamma_{jk}^{\bullet} = -\rho_{\ell} \wt{\rho}_{m} - \wt{\rho}_{\ell} \rho_{m} + 2\wt{\rho}_{\ell}\wt{\rho}_{m}\begin{cases}
    =0 &\quad \; \; \rho_{\ell} < 0, \; \rho_m < 0,\\
    >0 &\quad \; \; \rho_{\ell} > 0, \; \rho_m < 0,\\
    >0 &\quad \; \; \rho_{\ell} < 0, \; \rho_m > 0,\\
    =0  &\quad \;\; \rho_{\ell} > 0, \; \rho_m > 0,\\
    =0  &\quad \;\; \rho_{\ell} = 0 \; \text{or} \; \rho_m = 0.                    
\end{cases}
\end{equation*}
The diagonal entries $\gamma_{jj}^{\bullet}, j > 1$, equal zero since for $\ell = m$, $\rho_{\ell} = \rho_m$
and consequently $\text{sign}(\rho_{\ell}) = \text{sign}(\rho_m)$, which confirms (1). 
Concerning (4), note that for $j \neq k$, $\omega_{jk}^{\bullet} < 0$ if and
only if $\delta_{\ell}=\rho_{\ell} - \wt{\rho}_{\ell} > 0$ and $\delta_m = \rho_{m} - \wt{\rho}_{m} > 0$ if and only if 
$\rho_{\ell} > 0$ and $\rho_m > 0$ so that the corresponding $\gamma_{jk}^{\bullet}
= 0$. Point (5) follows with the converse argument, starting from 
$\omega_{jk}^{\bullet} = 0$ if and only if $\delta_{\ell} = 0$ or $\delta_m = 0$.
We now check the stationarity condition
$\Omega_{\bullet} = (\Sigma_* + \Gamma_{\bullet})^{-1}$.
We first verify the expression $I - \delta \delta^{\T}/(1 +
\nnorm{\delta}_2^2)$ for the bottom right block of $\Omega_{\bullet}$. Using
Schur complements, we get that the block is given by
\begin{align*}
\left(I + \rho \rho^{\T} - \rho \wt{\rho}^{\T} -
              \wt{\rho} \rho^{\T} + 2 \wt{\rho} \wt{\rho}^{\T} - \wt{\rho}
              \wt{\rho}^{\T} \right)^{-1}&= \left( I + (\rho - \wt{\rho})(\rho - \wt{\rho})^{\T} \right)^{-1} \\
&= (I + \delta \delta^{\T})^{-1} = I - \frac{\delta \delta^{\T}}{1 + \nnorm{\delta}_2^2}
\end{align*}
by the Sherman-Woodbury-Morrison formula. For the remaining  blocks, 
one verifies that $\Omega_{\bullet} \Sigma_{\bullet} = I$ directly with matrix multiplication, noting that 
\begin{equation*}
\wt{\rho} - (\wt{\rho} - \delta \delta^{\T} \wt{\rho}/(1 +
\nnorm{\delta}_2^2)) = 0, 
\end{equation*}
because $\delta^{\T} \wt{\rho} = 0$. To conclude the proof, it remains to check that 
$\Omega_{\bullet} \succ 0$. We show that all principal minors of $\Omega_{\bullet}$
are positive, which requires the following conditions to hold:
\begin{equation*}
I - \frac{\delta \delta^{\T}}{1 + \nnorm{\delta}_2^2} \succ 0, \; \, \text{and} \; \, ( 1 +
\nnorm{\wt{\rho}}_2^2 - \wt{\rho}^{\T} (I - \delta \delta^{\T}/(1 +
\nnorm{\delta}_2^2))^{-1} \wt{\rho}) > 0.
\end{equation*}
The first part follows by noting that the smallest eigenvalue of the matrix
is given by $1 - \nnorm{\delta}_2^2/(1 + \nnorm{\delta}_2^2) > 0$. The second
part results from $\delta^{\T} \wt{\rho} = 0$.$\qed$

\section{Proof of Theorem \ref{theo:convergence}}\label{proof:convergence}
It is not hard to verify that the conditions of Proposition 2.7.1 in \cite{Bertsekas1999}, a general
result concerning convergence of block coordinate descent, are satisfied. In particular, as 
discussed in Section \ref{sec:optimization}, the problems associated with each coordinate block are strictly convex
and hence have a unique minimizer. Proposition 2.7.1 in \cite{Bertsekas1999} then yields that each limit point of the 
sequence of iterates $\{ \Omega^t \}$  is a stationary point and thus the unique minimizer $\wh{\Omega}$. Existence
of a limit point requires the $\{ \Omega^t \}$ be contained in a compact set. This follows from 
\begin{equation*}
\lim_{\norm{\Omega} \rightarrow \infty} -\log \det(\Omega) + \tr(\Omega S) = +\infty,
\end{equation*}
which can be established using the reasoning leading to the proof of Theorem \ref{theo:existenceanduniqueness}, 
and the fact that for all $t$
\begin{equation*}
 -\log \det(\Omega^{t+1}) + \tr(\Omega^{t+1} S) \leq  -\log \det(\Omega^{t}) + \tr(\Omega^t S). 
\end{equation*}  

\section{Proof of Proposition \ref{prop:fixedp}}\label{proof:fixedp}
The proof will be reduced to a general scheme for establishing consistency
of $M$-estimators. An $M$-estimator is defined as maximizer $\wh{\theta}$ of a function 
of the form
\begin{equation*}
\theta \mapsto M_n(\theta) \coloneq \frac{1}{n} \su m_{\theta}(X_i), 
\end{equation*}
over some metric space $(\Theta, d)$, where the random variables $\{ X_i \}_{i = 1}^n$ represent the samples
drawn i.i.d.~according to a certain probability measure.
\renewcommand{\thesection}{\Alph{section}}
\begin{theoremApp}\label{theo:consistencyM}[from Theorem 5.14 in \cite{vanderVaart1998}]
Let the following conditions be fulfilled.
\begin{enumerate}
\item The map $\theta \mapsto m_{\theta}(x)$ is upper-semicontinuous for almost all $x$.
\item For every ball $U \subset \Theta$, $\E \left[\sup_{\theta \in U} m_{\theta}(X) \right] < \infty$.
\item  There exists a compact set $K \subset \Theta$ so that $\p(\wh{\theta} \in K) \rightarrow 1$ as $n \rightarrow \infty$.  
\end{enumerate}
Then: $d(\wh{\theta}, \Theta_{\bullet}) \rightarrow 0$ as $n \rightarrow \infty$, where $\Theta_{\bullet} = \argmax_{\theta \in \Theta} \E[m_{\theta}(X)]$.  
\end{theoremApp}
Before applying the above theorem, we first state and prove the
following lemma.
\begin{lemmaApp}\label{lem:consistencyspectral} 
If the random vector $X$ has finite fourth moments, that is 
$\E[X_j^4] < \infty$ for all $j=1,\ldots,p$, then the spectrum of the sample covariance
matrix $S = \frac{1}{n} \sum_{i = 1} (x_i - \mu_*) (x_i - \mu_*)^{\T}$ satisfies
\begin{equation*}
\lambda_j(S) = \lambda_j(\Sigma_*) + o_{\p}(1), \;j=1,\ldots,p,
\end{equation*}
as $n \rightarrow \infty$ (and $p$ stays fixed). 
\end{lemmaApp}
\begin{bew} The assumption of having finite fourth moments implies that \\ 
$s_{jk} = \sigma_{jk}^* + o_{\p}(1) \; \forall (j,k)$ and hence also that
\begin{equation*} 
\norm{S - \Sigma_{*}}_F \leq p^{1/2} \, \max_{(j,k)} |s_{jk} - \sigma_{jk}^*| = o_{\p}(1).  
\end{equation*}
The claim then follows from the Hoffmann-Wielandt Theorem. 
\end{bew}  

\paragraph{Proof of Proposition \ref{prop:fixedp}} 
First note that the sign-constrained log-determinant divergence minimization
\eqref{eq:logdetext} falls under the framework of M-estimation with 
$\Theta = \overline{\mc{M}^p}$, $\wh{\theta} = \wh{\Omega}$, and we may take $d$ as the metric that is induced by the spectral norm. For the 
function $m_{\theta}$, we have $m_{\theta}(X) = \log \det(\theta) - \tr(\theta X X^{\T})$. We now verify all three conditions of Theorem \ref{theo:consistencyM}. The first condition obviously holds true. Regarding 2., let $U = \{\theta \in \Theta: \nnorm{\theta - \theta_0} \leq r \}$
for some $\theta_0 \in \Theta$ a ball of radius $r > 0$, where $\nnorm{\cdot}$ denotes the spectral norm. We have   
\begin{align*}
\E \left[\sup_{\theta \in U} \log \det(\theta) - \tr(\theta X X^{\T}) \right] &\leq \sup_{\theta \in U} \log \det(\theta) < \infty,
\end{align*}
since $\E[\tr(\theta X X^{\T})] \geq 0$ in view of $\theta \in \overline{\mathbb{S}_+^p}$, and the first term is bounded from above, because so is $\sup_{\theta \in U} \nnorm{\theta}$. We finally turn to 3. Using an eigen-expansion of $\wh{\Omega}$, we have
\begin{align*}
\log \det(\wh{\Omega}) - \tr(\wh{\Omega} S) &\leq p \log(\lambda_1(\wh{\Omega})) - \lambda_1(\wh{\Omega}) \lambda_p(S) \\
&\leq \underbrace{p \log(\lambda_1(\wh{\Omega})) - \lambda_1(\wh{\Omega})  \{ \lambda_p(\Sigma_*)}_{\invcoloneq U(\lambda_1(\wh{\Omega}))} + o_{\p}(1)  \},
\end{align*}
using Lemma \ref{lem:consistencyspectral}. On the other hand, since $\wh{\Omega}$ is a minimizer and $\Omega_{\bullet}$ \eqref{eq:logdet_population_Omega} is feasible,   
\begin{align*}
\log \det(\wh{\Omega}) - \tr(\wh{\Omega} S) &\geq  \log \det(\Omega_{\bullet}) - \tr(\Omega_{\bullet} S) \\
                                            &\geq  p \log (\lambda_p(\Omega_{\bullet}))  -  \lambda_1(S) \sum_{j = 1}^p \lambda_j(\Omega_{\bullet})\\
                                            &=\underbrace{p \log (\lambda_p(\Omega_{\bullet})) -  \lambda_1(\Sigma_*) \sum_{j = 1}^p \lambda_j(\Omega_{\bullet})}_{\invcoloneq L} + o_{\p}(1)
\end{align*}
Therefore, with probability tending to one as $n \rightarrow \infty$, $\wh{\Omega}$ is contained in the compact set 
\begin{equation*}
K = \{\Omega:\, \norm{\Omega} \leq B(U, L) \}, \quad \text{where} \; B(U,L) = \sup \{b \geq 0:\, U(b) \geq L \}.
\end{equation*}  
Note that $b$ must be bounded from above as $\lim_{b \rightarrow  \infty} U(b) = -\infty$. $\qed$  

\renewcommand{\thesection}{Appendix \Alph{section}}

\section{Proof of Proposition \ref{prop:singleedge_recovery}}\label{app:singleedge_recovery}

\renewcommand{\thesection}{\Alph{section}}
Our proof depends on the following lemma. 
\begin{lemmaApp}\label{lem:supbound_nnlinearsystem} 
Consider the system of linear equations $A x = b$, where $A \in \R_{+}^{d
  \times d}$, $x \in \R_{+}^d$ and $b \in \R_+^d$ have only non-negative
entries. Then $x_j \leq b_j/a_{jj}$ for all $j=1,\ldots,d$.
\end{lemmaApp}

\begin{bew} For any $j \in \{1,\ldots,d \}$, we have that
$a_{jj} x_j \leq \sum_{k = 1}^d a_{jk} x_k = b_{j}$, using the non-negativity of all entries.
\end{bew}

\paragraph{Proof of Proposition \ref{prop:singleedge_recovery}}

We first prove that $\max_{(j,k) \notin \mc{E}^*, \, j \neq k} (-\wh{\omega}_{jk}) \leq c_1B$.  
From the KKT optimality conditions \eqref{eq:kkt}, we have that
$\wh{\Sigma} = S + \wh{\Gamma}$ with $\wh{\gamma}_{jj} = 0$ for all $j$ and $\wh{\gamma}_{jk} = 0$ whenever
\begin{equation*}
(j,k) \in \wh{\mc{E}}, \quad \text{where} \; \, \wh{\mc{E}} = \{(j,k): \; \; \wh{\omega}_{jk} < 0\},
\end{equation*}
and consequently
\begin{equation}\label{eq:noslack}
(j,k) \in \wh{\mc{E}} \; \Longrightarrow \; \wh{\sigma}_{jk} = s_{jk}.  
\end{equation}
Now choose $(\bar{j}, \bar{k}) \in \triangle = \wh{\mc{E}} \setminus \mc{E}^*$ 
(if $\triangle = \emptyset$, the claim would follow trivially) 
such that
\begin{equation*}
\wh{\omega}_{\bar{j} \bar{k}} = \min_{(j,k) \in \triangle} \wh{\omega}_{jk}. 
\end{equation*}
Using the partitioning scheme \eqref{eq:partitioning} with $\wh{\Sigma}$
respectively $\wh{\Omega}$ and $j = \bar{j}$, and using Schur complements, we obtain that
\begin{equation}\label{eq:statcol}
\wh{\sigma}_{\bar{j}} = \frac{ \wh{\Omega}_{\bar{j} \bar{j}}^{-1} \left(-\wh{\omega}_{\bar{j}}
  \right)}{\wh{\omega}_{\bar{j} \bar{j}}- \wh{\omega}_{\bar{j}}^{\T}
  \wh{\Omega}_{\bar{j} \bar{j}}^{-1}
  \wh{\omega}_{\bar{j}}} =  \wh{\sigma}_{\bar{j} \bar{j}}\wh{\Omega}_{\bar{j} \bar{j}}^{-1} \left(-\wh{\omega}_{\bar{j}}
  \right) = s_{\bar{j} \bar{j}} \wh{\Omega}_{\bar{j} \bar{j}}^{-1} \left(-\wh{\omega}_{\bar{j}} \right).
\end{equation}
Using Schur complements again, 
\begin{equation}\label{eq:condcov}
\wh{\Omega}_{\bar{j} \bar{j}}^{-1} = \wh{\Sigma}_{\bar{j} \bar{j}} - \frac{\wh{\sigma}_{\bar{j}}
  \wh{\sigma}_{\bar{j}}^{\T}}{\wh{\sigma}_{\bar{j} \bar{j}}} = 
\wh{\Sigma}_{\bar{j} \bar{j}} - \frac{\wh{\sigma}_{\bar{j}}
  \wh{\sigma}_{\bar{j}}^{\T}}{s_{\bar{j} \bar{j}}}
\end{equation}
Combining \eqref{eq:statcol} and \eqref{eq:condcov}, we obtain
\begin{equation}\label{eq:statcol_comb}
\underbrace{\left(\wh{\Sigma}_{\bar{j} \bar{j}} - \frac{\wh{\sigma}_{\bar{j}}
  \wh{\sigma}_{\bar{j}}^{\T}}{\wh{\sigma}_{\bar{j} \bar{j}}}
\right)}_{\invcoloneq \wt{\Sigma}_{\bar{j} \bar{j}}}
(-\wh{\omega}_{\bar{j}}) = \frac{\wh{\sigma}_{\bar{j}}}{s_{\bar{j} \bar{j}}}.
\end{equation}
Let $\mc{A}_{\bar{j}} = \{l \in \{1,\ldots,p-1 \}:\, (-\wh{\omega}_{\bar{j}})_l
> 0\}$. Then, \eqref{eq:statcol_comb} can equivalently be written as 
\begin{align}\label{eq:activepart}
&\left(\wt{\Sigma}_{\bar{j} \bar{j}} \right)_{\mc{A}_{\bar{j}} \,
  \mc{A}_{\bar{j}}} (-\wh{\omega}_{\bar{j}})_{\mc{A}_{\bar{j}}} = \left(
\frac{s_{\bar{j}}}{s_{\bar{j} \bar{j}}} \right)_{\mc{A}_{\bar{j}}},\\
&\left(\wt{\Sigma}_{\bar{j} \bar{j}} \right)_{\mc{A}_{\bar{j}}^c \,
  \mc{A}_{\bar{j}}} (-\wh{\omega}_{\bar{j}})_{\mc{A}_{\bar{j}}} = \left(
\frac{s_{\bar{j}} + \gamma_{\bar{j}}}{s_{\bar{j} \bar{j}}} \right)_{\mc{A}_{\bar{j}}^c} \notag,
\end{align}
where we have used \eqref{eq:noslack}. In order to upper bound $(-\wh{\omega}_{\bar{j} \bar{k}})$, we
consider \eqref{eq:activepart}. Since $\wh{\Omega}_{\bar{j} \bar{j}}$ is an $M$-matrix, its inverse $\wt{\Sigma}_{\bar{j} \bar{j}}$ has only non-negative
entries, and we are hence in position to apply Lemma \ref{lem:supbound_nnlinearsystem}.  We obtain
\begin{equation*}
(-\wh{\omega}_{\bar{j} \bar{k}}) \leq \frac{s_{\bar{j} \bar{k}} / s_{\bar{j}
    \bar{j}}}{\wh{\sigma}_{\bar{k} \bar{k}} - \wh{\sigma}_{\bar{j}
        \bar{k}}^2 / \wh{\sigma}_{\bar{j} \bar{j}}} = \frac{s_{\bar{j} \bar{k}} / s_{\bar{j}
    \bar{j}}}{s_{\bar{k} \bar{k}} - s_{\bar{j}
        \bar{k}}^2/s_{\bar{j} \bar{j}}},
\end{equation*}
where the second equality is again a consequence of \eqref{eq:noslack}. Using
that $\sigma_{\bar{j} \bar{k}}^{*} = 0$, $\sigma_{\bar{j} \bar{j}}^{*} =
\sigma_{\bar{k} \bar{k}}^{*} = 1$, the bound \eqref{eq:supbound} yields
\begin{equation}\label{eq:supbound_final}
(-\wh{\omega}_{\bar{j} \bar{k}}) \leq \frac{s_{\bar{j} \bar{k}}}{s_{\bar{j} \bar{j}} 
  s_{\bar{k} \bar{k}} - s_{\bar{j}
        \bar{k}}^2} \leq \frac{B}{1 - 2 B - B^2} \leq C_0 B \; \, \text{as} \; n \rightarrow \infty,
\end{equation}
since $B = o_{\p}(1)$ as $n \rightarrow \infty$. In the sequel, we derive a lower bound on the entry of $(-\wh{\omega}_{12})$
corresponding to $\mc{E}^*$. For this purpose, let us re-consider
Eq.~\eqref{eq:statcol_comb} for $\bar{j} = 1$, that is 
$\wt{\Sigma}_{11} (-\wh{\omega}_1) = \wh{\sigma}_1/s_{11}$. Expanding this
equation entry-wise, we get    
\begin{align*} 
\underbrace{\begin{pmatrix}
\wt{\sigma}_{22} & \wt{\sigma}_{23} & \ldots & \wt{\sigma}_{2p} \\
\wt{\sigma}_{23} & \wt{\sigma}_{33} & \ldots & \wt{\sigma}_{3p} \\
 \vdots               &    \vdots              & \ddots       & \vdots                  \\
\wt{\sigma}_{2p} & \wt{\sigma}_{3p} & \ldots & \wt{\sigma}_{pp}
\end{pmatrix}}_{\wt{\Sigma}_{11}} 
\begin{pmatrix}
-\wh{\omega}_{12} \\
-\wh{\omega}_{13} \\
\vdots \\
-\wh{\omega}_{1p}
\end{pmatrix} = \begin{pmatrix}
\wh{\sigma}_{12}/s_{11} \\
\wh{\sigma}_{13}/s_{11} \\
\vdots \\
\wh{\sigma}_{1p}/s_{11}
\end{pmatrix},  
\end{align*}
with $\wt{\sigma}_{jk} = \wh{\sigma}_{jk} - \wh{\sigma}_{1j} \wh{\sigma}_{1k} / s_{11}, \; j,k=2,\ldots,p$.
Consider now the top equation
\begin{equation}\label{eq:topequation}
\wt{\sigma}_{22} (-\wh{\omega}_{12}) + \underbrace{\sum_{l = 3}^p
  \wt{\sigma}_{2l} (-\wh{\omega}_{1l})}_{\invcoloneq \delta} =
\wh{\sigma}_{12}/s_{11}. 
\end{equation}
The order of the term $\delta$, which is the inner product of the first row of
$\wt{\Sigma}_{11}$ (excluding the diagonal element $\wt{\sigma}_{22}$) and
$(-\wh{\omega}_{13},\ldots,-\wh{\omega}_{1p})$, can be upper bounded by taking
the corresponding inner products associated with the remaining rows of
$\wt{\Sigma}_{11}$ as a reference, noting that the $\{\wt{\sigma}_{jk}, \, j \neq
k, \, (j,k) \notin \{1,2\} \}$ are of the same order, since the $\{\wh{\sigma}_{jk}, \, j \neq
k, \, (j,k) \notin \{1,2\} \}$ are exchangeable. Accordingly, the right hand
sides $\wh{\sigma}_{13}/s_{11},\ldots,\wh{\sigma}_{1p}/s_{11}$ are also of the same
order, which is at most $O_{\p}(B)$ in view of the complementarity slackness
condition \eqref{eq:noslack} and the scaling of the $\{ s_{jk}/s_{jj},\, j \neq
k, \, (j,k) \notin \{1,2\} \}$ (if
$\wh{\omega}_{13}=\ldots=\wh{\omega}_{1p}=0$, we would have $\delta =
0$). Formally, the argument reads
\begin{align}\label{eq:Delta_order}
\delta = \sum_{l=3}^p \wt{\sigma}_{2l} (-\wh{\omega}_{1l}) \leq C_1 \max_{3
  \leq j \leq p} \sum_{l = 3}^p \wt{\sigma}_{jl} (-\wh{\omega}_{1l}) &= \max_{3
\leq j \leq p} \wh{\sigma}_{1j}/s_{11} \notag\\
&\leq C_2 \max_{3
  \leq j \leq p} s_{1j}/s_{11} \leq C_{3} B,
\end{align}
as $n \rightarrow \infty$. Suppose for a moment that $(-\wh{\omega}_{12}) > 0$ so that $\wh{\sigma}_{12}
= s_{12}$. Substituting \eqref{eq:Delta_order} back into
\eqref{eq:topequation} and resolving for $(-\wh{\omega}_{12})$, we obtain
\begin{align}\label{eq:lowerbound_final}
(-\wh{\omega}_{12}) = \frac{\wh{\sigma}_{12}/s_{11} - \delta}{\wt{\sigma}_{22}} &=
\frac{s_{12}/s_{11} - \delta}{s_{22} - s_{12}^2/s_{11}} \notag \\ 
&\geq \frac{\rho - C_4 B}{1 - \rho^2 + C_5 B} = (-\omega_{12}^*) - C_6
B,
\end{align}
as  $n \rightarrow \infty$, i.e.~if $\rho > C_4 B$, we verify that indeed $(-\wh{\omega}_{12}) >
0$. Altogether, \eqref{eq:supbound_final} and \eqref{eq:lowerbound_final}
indicate that if $(-\wh{\omega}_{12}) > t,  \; t = (C_0 + C_6) B$, the thresholding procedure
\eqref{eq:Omega_thres} would yield $\wh{\mc{E}}(t) = \mc{E}^*$, i.e.~recovery of
the edge set. $\qed$

\renewcommand{\thesection}{\Alph{section}}

\section{$M$-matrices and Faithfulness}
We show that if the precision matrix of a Gaussian random vector is an $M$-matrix, then faithfulness as defined below holds. 

\renewcommand{\thesection}{\Alph{section}}

\begin{defnApp}\cite{Spirtes2000, Buhlmann2011} Let $X$ be a multivariate Gaussian random vector
with covariance matrix $\Sigma \in \mathbb{S}_+^p$. Then $X$ is said to have a faithful
distribution if for any disjoint triple $A, B, C$ of subsets of
$\{1,\ldots,p \}$,
\begin{align*}
X_A \indep X_B \,| X_C \; \Longrightarrow \; &C \, \textbf{\text{separates}} \, A \, \text{and} \, B \\
&\text{in the conditional independence graph associated with} \, \Sigma^{-1},
\end{align*}
that is, vertices in $A$ are connected with vertices in $B$ at most only via vertices in $C$.
\end{defnApp}
Note that the converse statement, i.e.~separation of $A$ and $B$ by $C$ in the conditional independence graph implies 
conditional independence of $X_A$ and $X_B$ given $X_C$, always holds (global Markov property).   
\begin{propApp} If $\Sigma^{-1} \in \mc{M}^p$, then $X$ has a faithful distribution. 
\end{propApp}
\begin{bew} In the sequel, we will show that 
\begin{align}\label{eq:hierarch}
X_A \indep X_B \,| X_C \; \Longrightarrow X_A \indep X_B \,| X_{C'} \; \,\forall C' \supseteq C, \; \, C' \subseteq \{1,\ldots,p\} \setminus (A \cup B). 
\end{align}
Choosing $C' =  \{1,\ldots,p\} \setminus (A \cup B)$, the claim follows from the global Markov property w.r.t.
the conditional independence graph associated with $\Sigma^{-1}$. Let $\gamma = A \cup B \cup C$, $\gamma' = A \cup B \cup C'$ and
 $\Omega^{\gamma\gamma} = \Sigma_{\gamma \gamma}^{-1}$, $\;\Omega^{\gamma'\gamma'} = \Sigma_{\gamma' \gamma'}^{-1}$. 
We will compute $\Omega^{\gamma' \gamma'}$ incrementally from $\Omega^{\gamma \gamma}$ by using the decomposition    
$C' \setminus C = \{ i_1 \} \cup \ldots \cup \{ i_q \}$ ($q = |C' \setminus C|$), successively obtaining 
$\Omega^{\gamma_1 \gamma_1}, \ldots, \Omega^{\gamma_q \gamma_q} = \Omega^{\gamma' \gamma'}$, where $\gamma_1 = A \cup B \cup C \cup \{ i_1 \},\ldots,
\gamma_q = A \cup B \cup C \cup \{ i_1 \} \cup \ldots \cup \{ i_q \} = \gamma'$. 
Starting from $\Omega^{\gamma_1 \gamma_1}$, we partition its inverse $\Sigma_{\gamma_1 \gamma_1}$ as
\begin{equation*}
\Sigma_{\gamma_1 \gamma_1} =  \begin{pmatrix}
                         \Sigma_{\gamma\gamma} & u_1^{\T} \\
                          u_1 & u_{11} \\     
                         \end{pmatrix},
\end{equation*}
where the vector $u_1$ and the scalar correspond to the added index $i_1$. The partitioned inverse formula yields 
the following for $\Omega^{\gamma_1 \gamma_1}_{\gamma \gamma}$, the principal submatrix of $\Omega^{\gamma_1 \gamma_1}$
associated with the index set $\gamma$:
\begin{align*}
\Omega^{\gamma_1 \gamma_1}_{\gamma \gamma} = \left( \Sigma_{\gamma\gamma}  - \frac{u_1 u_1^{\T}}{u_{11}} \right)^{-1} 
                                  &= \Sigma_{\gamma \gamma}^{-1} + \underbrace{\frac{\Sigma_{\gamma \gamma}^{-1} u_1 u_1^{\T} \Sigma_{\gamma \gamma}^{-1}}{u_{11} - u_1^{\T} \Sigma_{\gamma \gamma}^{-1} u_1}}_{\invcoloneq P} = \Omega^{\gamma \gamma} + P,\\
\end{align*}
where the second identity results from the Sherman-Woodbury-Morrison formula.   
Note that $P$ has only non-negative entries. To verify this, observe that $\Sigma_{\gamma \gamma}^{-1} u_1$ equals
the vector of regression coefficients one obtains when regressing the variable with index $i_1$ on the variables in $\gamma$, 
which must be non-negative because $\Sigma^{-1} \in \mc{M}^p$ and hence also $\Sigma_{\gamma_1 \gamma_1}^{-1} \in \mc{M}^p$ (cf.~\eqref{eq:Omega_Schur} and the comments thereafter). Furthermore,
\begin{equation*}
u_{11} - u_1^{\T} \Sigma_{\gamma \gamma}^{-1} u_1 = \det(u_{11} - u_1^{\T} \Sigma_{\gamma \gamma}^{-1} u_1) = \det(\Sigma_{\gamma_1 \gamma_1})/\det(\Sigma_{\gamma \gamma}) > 0
\end{equation*}
in virtue of the positive definiteness of $\Sigma$. The non-negativity of $P$ implies that 
$\Omega^{\gamma \gamma}\leq \Omega_{\gamma \gamma}^{\gamma_1 \gamma_1}$. We may now repeat the same
argument to obtain successively $\Omega_{\gamma_1 \gamma_1}^{\gamma_2 \gamma_2},\ldots,\Omega_{\gamma_{q-1} \gamma_{q-1}}^{\gamma_q \gamma_q}$. Consequently,
we must have $\Omega^{\gamma \gamma} \leq \Omega_{\gamma \gamma}^{\gamma_q \gamma_q}$ and hence, since $X_A \indep X_B | X_C$, in particular that
$0 = \Omega^{\gamma \gamma}_{AB} \leq \Omega^{\gamma_q \gamma_q}_{AB} = \Omega^{\gamma' \gamma'}_{AB}$ (recall that $\gamma' = A \cup B \cup C'$). The inequality must hold 
with equality, because $\Omega^{\gamma'\gamma'}$ is a positive definite $M$-matrix (cf.~\eqref{eq:Omega_Schur}). We conclude
the assertion from $\Omega^{\gamma' \gamma'}_{AB} = 0 \Longleftrightarrow X_A \indep X_B | X_{C'}$.

\end{bew}

\renewcommand{\thesection}{Appendix \Alph{section}}





\bibliographystyle{plain}
\bibliography{../lit/references.bib}

\begin{thebibliography}{10}

\bibitem{Anandkumar2012_supp}
A.~Anandkumar, V.~Tan, F.~Huang, and A.~Willsky.
\newblock {Supplementary code to 'High-Dimensional Gaussian Graphical Model
  Selection:\\ $\text{Walk-Summability and Local Separation Criterion'}$}.
\newblock
  {http://newport.eecs.uci.edu/anandkumar/pubs/GaussianStructLearning-code.zip}.

\bibitem{Anandkumar2012}
A.~Anandkumar, V.~Tan, F.~Huang, and A.~Willsky.
\newblock High-dimensional graphical model selection: Tractable graph families
  and necessary conditions.
\newblock {\em Journal of Machine Learning Research}, 13:2293--2337, 2012.

\bibitem{Ban2008}
O.~Banerjee, L.~El Ghaoui, and A.~d'Aspremont.
\newblock {Model Selection Through Sparse Maximum Likelihood Estimation for
  Multivariate Gaussian or Binary Data}.
\newblock {\em Journal of Machine Learning Research}, 9:485--516, 2008.

\bibitem{BermanPlemmons}
A.~Berman and R.~Plemmons.
\newblock {\em Nonnegative matrices in the mathematical sciences}.
\newblock SIAM Classics in Applied Mathematics, 1994.

\bibitem{Berman2003}
A.~Berman and N.~Shaked-Monderer.
\newblock {\em Completely positive matrices}.
\newblock World Scientific, 2003.

\bibitem{Bertsekas1999}
D.~Bertsekas.
\newblock {\em {Nonlinear Programming}}.
\newblock Athena Scientific, 1999.

\bibitem{BoydVandenberghe2004}
S.~Boyd and L.~Vandenberghe.
\newblock {\em Convex Optimization}.
\newblock Cambridge University Press, 2004.

\bibitem{Buhlmann2011}
P.~B\"uhlmann and S.~van~de Geer.
\newblock {\em Statistics for high-dimensional data}.
\newblock Springer, 2011.

\bibitem{Cai2011}
T.~Cai, W.~Liu, and X.~Luo.
\newblock {A Constrained $\ell_1$ Minimization Approach to Sparse Precision
  Matrix Estimation}.
\newblock {\em Journal of the American Statistical Association}, 106:594--607,
  2011.

\bibitem{Chatfield2003}
C.~Chatfield.
\newblock {\em {The Analysis of Time Series: an introduction}}.
\newblock Chapmann \& Hall/CRC, 2003.

\bibitem{Drton2007}
S.~Chaudhuri, M.~Drton, and T.~Richardson.
\newblock Estimation of a covariance matrix with zeros.
\newblock {\em Biometrika}, 94:199--216, 2007.

\bibitem{ChowLiu1968}
C.~Chow and C.~Liu.
\newblock {Approximating discrete probability distributions with dependence
  trees}.
\newblock {\em IEEE Transactions on Information Theory}, 14:462--467, 1968.

\bibitem{Cootes}
T.~Cootes.
\newblock {$\text{XM2VTS face images}$}.
\newblock
  {http://personalpages.manchester.ac.uk/\\staff/timothy.f.cootes/data/xm2vts/xm2vts\verb?_?markup.html}.

\bibitem{Dahl2005}
J.~Dahl, V.~Roychowdhury, and L.~Vandenberghe.
\newblock {Maximum likelihood estimation of Gaussian graphical models:
  numerical implementation and topology selection}.
\newblock Technical report, University of California, Los Angeles, 2005.

\bibitem{Dempster1972}
A.~Dempster.
\newblock Covariance selection.
\newblock {\em Biometrics}, 28:157--175, 1972.

\bibitem{Dhillon2007}
I.~Dhillon and J.~Tropp.
\newblock {Matrix nearness problems with Bregman divergences}.
\newblock {\em SIAM Journal on Matrix Analysis and Applications},
  29:1120--1146, 2007.

\bibitem{DrydenMardia2002}
I.~Dryden and K.~Mardia.
\newblock {\em {Statistical Shape Analysis}}.
\newblock Wiley, 2002.

\bibitem{Duchi2009}
J.~Duchi, S.~Gould, and D.~Koller.
\newblock {Projected subgradient methods for learning sparse Gaussians}.
\newblock In {\em Artificial Intelligence and Statistics (AISTATS)}, 2009.

\bibitem{Bolviken1982}
{E. B\o{}lviken}.
\newblock {Probability inequalities for the multivariate normal with
  non-negative partial correlations}.
\newblock {\em Scandinavian Journal of Statistics}, 9:49--58, 1982.

\bibitem{Fan2009}
J.~Fan, Y.~Feng, and Y.~Wu.
\newblock {Network exploration via the adaptive lasso and SCAD penalties}.
\newblock {\em The Annals of Applied Statistics}, 3:521--541, 3.

\bibitem{Foygel2010}
R.~Foygel and M.~Drton.
\newblock {Extended Bayesian information criteria for Gaussian graphical
  models}.
\newblock In {\em Advances in Neural Information Processing Systems 23}, pages
  2020--2028, 2010.

\bibitem{Friedman2008}
J.~Friedman, T.~Hastie, and R.~Tibshirani.
\newblock {Sparse inverse covariance estimation with the graphical lasso}.
\newblock {\em Biostatistics}, 9:432--441, 2008.

\bibitem{cvx}
M.~Grant and S.~Boyd.
\newblock {CVX}: Matlab software for disciplined convex programming, version
  1.21, 2011.

\bibitem{Gu07}
L.~Gu, E.~Xing, and T.~Kanade.
\newblock {Learning GMRF Structures for Spatial Priors}.
\newblock In {\em CVPR}, 2007.

\bibitem{Honorio2009}
J.~Honorio, L.~Ortiz, and D.~Samaras.
\newblock {Sparse and Locally Constant Gaussian Graphical Models}.
\newblock In {\em Advances in Neural Information Processing Systems 22}. 2009.

\bibitem{mturk}
http://www.mturk.com.

\bibitem{Kalisch2007}
M.~Kalisch and P.~B\"uhlmann.
\newblock {Estimating High-Dimensional Directed Acyclic Graphs with the
  PC-Algorithm}.
\newblock {\em Journal of Machine Learning Research}, 8:613--636, 2007.

\bibitem{Kalisch2012}
M.~Kalisch, M.~M\"achler, D.~Colombo, M.~Matthuis, and P.~B\"uhlmann.
\newblock {Causal inference using graphical models with the R package pcalg}.
\newblock {\em Journal of Statistical Software}, 47:1--26, 2012.

\bibitem{Karlin1980}
S.~Karlin and Y.~Rinott.
\newblock {Classes of orderings of measures and related correlation
  inequalities I: multivariate totally positive distributions}.
\newblock {\em Journal of Multivariate Analysis}, 10:467--498, 1980.

\bibitem{Karlin1983}
S.~Karlin and Y.~Rinott.
\newblock M-matrices as covariance matrices of multinormal distributions.
\newblock {\em {Linear Algebra and Its Applications}}, 52:419--438, 1983.

\bibitem{glasso_matlab}
H.~Karshenas.
\newblock {Graphical lasso in \textsf{R} and \textsf{MATLAB}}.
\newblock http://www-stat.stanford.edu/~tibs/glasso/.

\bibitem{LakeTenenbaum2010}
B.~Lake and J.~Tenenbaum.
\newblock {Discovering Structure by Learning Sparse Graphs}.
\newblock In {\em {Proceedings of the 33rd Annual Cognitive Science
  Conference}}, 2010.

\bibitem{Lauritzen1996}
S.~Lauritzen.
\newblock {\em {Graphical Models}}.
\newblock {Oxford University Press}, 1996.

\bibitem{LuoTseng1992}
Z.~Luo and P.~Tseng.
\newblock On the convergence of the coordinate descent method for convex
  differentiable minimization.
\newblock {\em Journal of Optimization Theory and Applications}, 72:7--35,
  1992.

\bibitem{Malioutov2006}
D.~Malioutov, J.~Johnson, and A.~Willsky.
\newblock {Walk-Sums and Belief Propagation}.
\newblock {\em Journal of Machine Learning Research}, 7:2031--2064, 2006.

\bibitem{Mardia1979}
K.~Mardia, J.~Kent, and J.~Bibby.
\newblock {\em Multivariate Analysis}.
\newblock Academic Press, 1979.

\bibitem{Meinshausen2013}
N.~Meinshausen.
\newblock {Sign-constrained least squares estimation for high-dimensional
  regression}.
\newblock {\em The Electronic Journal of Statistics}, 7:1607--1631, 2013.

\bibitem{Mei2006}
N.~Meinshausen and P.~B\"uhlmann.
\newblock High-dimensional graphs and variable selection with the lasso.
\newblock {\em The Annals of Statistics}, 34:1436--1462, 2006.

\bibitem{Osherson1991}
D.~Osherson, J.~Stern, O.~Wilkie, M.~Stob, and E.~Smith.
\newblock Default probability.
\newblock {\em Cognitive Science}, 15:251--269, 1991.

\bibitem{Ostrowski1937}
A.~Ostrowski.
\newblock {\"Uber die Determinanten mit \"uberwiegender Hauptdiagonale}.
\newblock {\em Commentarii Mathematici Helvetici}, 10:69--96, 1937.

\bibitem{Portugal1994}
L.~Portugal, J.~Judice, and L.~Vicente.
\newblock A comparison of block pivoting and interior-point algorithms for
  linear least squares problems with non-negative variables.
\newblock {\em Mathematics of Computation}, 63:625--643, 1994.

\bibitem{Ravikumar2011}
P.~Ravikumar, M.~Wainwright, G.~Raskutti, and B.~Yu.
\newblock {High-dimensional covariance estimation by minimizing
  $\ell_1$-penalized log-determinant divergence}.
\newblock {\em The Electronic Journal of Statistics}, 4:935--980, 2011.

\bibitem{Rinott2004}
Y.~Rinott and M.~Scarsini.
\newblock Total positivity order and the normal distribution.
\newblock {\em Journal of Multivariate Analysis}, 97:1251--1261, 2004.

\bibitem{Rothman2008}
A.~Rothman, P.~Bickel, L.~Levina, and J.~Zhu.
\newblock {Sparse permutation invariant covariance estimation}.
\newblock {\em The Electronic Journal of Statistics}, 2:494--515, 2008.

\bibitem{Rue2001}
H.~Rue and L.~Held.
\newblock {\em {Gaussian Markov Random Fields}}.
\newblock Chapman and Hall/CRC, Boca Raton, 2001.

\bibitem{L1general}
M.~Schmidt.
\newblock {\textsf{L1General}: minimizing differentiable functions
  with\L1-\\regularization}.
\newblock {http://www.di.ens.fr/$\sim$mschmidt/Software/L1General.html}.

\bibitem{Shen2012}
X.~Shen, W.~Pan, and Y.~Zhu.
\newblock {Likelihood-based selection and sharp parameter estimation}.
\newblock {\em Journal of the American Statistical Association}, 107:223--232,
  2012.

\bibitem{Slawski2013nnlsalg}
M.~Slawski.
\newblock Problem-specific peformance analysis of non-negative least squares
  solvers with a focus on instances with sparse solutions.
\newblock Working manuscript, 2013.

\bibitem{SlawskiHein2011nips}
M.~Slawski and M.~Hein.
\newblock Sparse recovery by thresholded non-negative least squares.
\newblock In {\em Advances in Neural Information Processing Systems 24}, pages
  1926--1934. 2011.

\bibitem{Spirtes2000}
P.~Spirtes, C.~Glymour, and R.~Scheines.
\newblock {\em {Causation, Prediction, and Search}}.
\newblock MIT press, 2000.

\bibitem{Uhler2012}
C.~Uhler.
\newblock {Geometry of maximum likelihood estimation in Gaussian graphical
  models}.
\newblock {\em The Annals of Statistics}, 40:238--261, 2012.

\bibitem{vanderVaart1998}
A.~van~der Vaart.
\newblock {\em {Asymptotic Statistics}}.
\newblock Cambridge University Press, 1998.

\bibitem{Whittaker1990}
J.~Whittaker.
\newblock {\em {Graphical Models in Applied Multivariate Statistics}}.
\newblock Wiley, Chichester, 1990.

\bibitem{Yuan2010}
M.~Yuan.
\newblock {High Dimensional Inverse Covariance Matrix Estimation via Linear
  Programming}.
\newblock {\em Journal of Machine Learning Research}, 11:2261--2286, 2010.

\bibitem{Yuan2007}
M.~Yuan and Y.~Lin.
\newblock {Model Selection and Estimation in the Gaussian Graphical Model }.
\newblock {\em Biometrika}, 94:19--35, 2007.

\bibitem{Zhou2011}
S.~Zhou, P.~R\"utimann, M.~Xu, and P.~B\"uhlmann.
\newblock {High-dimensional covariance estimation based on Gaussian graphical
  models}.
\newblock {\em Journal of Machine Learning Research}, 12:2975--3026, 2011.

\end{thebibliography}







\end{document}